\documentclass[12pt]{amsart}

\usepackage{amsmath}
\usepackage{amstext}
\usepackage{amscd}
\usepackage{amssymb}

\usepackage[dvips]{graphicx}

\def\showfig{\iftrue}

\newtheorem{thm}[subsection]{Theorem}
\newtheorem{lem}[subsection]{Lemma}
\newtheorem{cor}[subsection]{Corollary}
\newtheorem{claim}[subsubsection]{Claim}
\newtheorem{prop}[subsection]{Proposition}

\newtheorem{rem}[subsection]{Remark}
\newtheorem{exmp}[subsection]{Example}
\newtheorem{defn}[subsection]{Definition}

\makeatletter
\@addtoreset{subsubsection}{section}
\@addtoreset{subsubsection}{subsection}
\@addtoreset{equation}{subsection}
\@addtoreset{figure}{subsection}
\@addtoreset{table}{subsection}
\makeatother

\newcommand\chEq{%
\setcounter{equation}{\value{subsection}}%
\addtocounter{subsection}{1}
}

 \newcommand\chFig{%
 \protect\setcounter{figure}{\value{subsection}}%
 \addtocounter{subsection}{1}
 }

\newcommand\chTbl{%
\setcounter{table}{\value{subsection}}%
\addtocounter{subsection}{1}
}

\newcommand\SHOWFIG[1]{
 \showfig{#1}\fi
}

\providecommand\eqref[1]{(\ref{#1})}

\newcommand\CCC{{\mathbb C}}

\newcommand\RRR{{\mathbb R}}

\newcommand\ZZZ{{\mathbb Z}}

\newcommand\proofstyle[1]{{\it #1}} 
\providecommand\proof{\proofstyle{Proof.~}}

\newcommand\PathComp[2]{ {#1}_{#2} }
\newcommand\Frechet{Fr\'echet}

\newcommand\id{\mathrm{id}}
\newcommand\IM{\mathrm{Im}}
\newcommand\Int{\mathrm{Int}}
\newcommand\supp{\mathrm{supp\,}}
\newcommand\Fix{\mathrm{Fix\,}}
\newcommand\interior{\mathrm{Int}}

\newcommand\grad{\triangledown} 

\newcommand\rank{rank} 

\newcommand\Orbit{\mathcal{O}}
\newcommand\Stab{\mathcal{S}}
\newcommand\Diff{\mathcal{D}}
\newcommand\Aut{\mathrm{Aut}}

\newcommand\Cinf{C^{\infty}}
\newcommand\Sinf{\mathcal{S}}

\newcommand\Unity[1]{ \PathComp{#1}{\id} }

\newcommand\StabId{\Unity{\Stab}}

\newcommand\grp{G}

\newcommand\pnt{z}

\newcommand\eps{\varepsilon}

\newcommand\aCircle{S^1}
\newcommand\Psp{P}
\newcommand\manif{M}
\newcommand\tmanif{\widetilde{\manif}} 
\newcommand\flow{\Phi}

\newcommand\End{{\mathcal{E}}}

\newcommand\path{\omega}

\newcommand\cfunc{\gamma}

\newcommand\partitf{\Delta_{\mrsfunc}}

\newcommand\Nbh{\mathcal{N}}

\newcommand\nbh{U}

\newcommand\mrsfunc{f}
\newcommand\tmrsfunc{\widetilde{\mrsfunc}}

\newcommand\DiffP{\Diff({\Psp})}

\newcommand\DiffM{\Diff(\manif)}

\newcommand\dimM{m}

\newcommand\difM{h}

\newcommand\Fld{F}
\newcommand\Gfld{G}
\newcommand\Hfld{H}

\newcommand\singf{\Sigma_\mrsfunc}

\newcommand\cmap{h}

\newcommand\smfunc{\Cinf(\manif,\RRR)}
\newcommand\smmap{\Cinf(\manif,\manif)}

\newcommand\Shift{\varphi}

\newcommand\smone{\Cinf(\manif,\Psp)}

\newcommand\smr{\Cinf(\manif,\RRR)}

\newcommand\smrd{\Cinf_{\partial}(\manif,\RRR)}

\newcommand\smoned{\Cinf_{\partial}(\manif,\Psp)}
\newcommand\smm{\Cinf(\manif,\manif)}

\newcommand\afunc{\alpha}

\newcommand\bdifM{\bar{\difM}}

\newcommand\gdifM{\theta}

\newcommand\bnd{B}

\newcommand\gfunc{g} 

\newcommand\bndcnt{b}

\newcommand\genus{g}

\newcommand\fcr{(\mrsfunc,\singf)}
\newcommand\Mcr{(\manif,\singf)}
\newcommand\Orbfcr{\Orbit{\fcr}}

\newcommand\Orbf{\Orbit(\mrsfunc)}
\newcommand\Orbff{\Orbit_{\mrsfunc}(\mrsfunc)}
\newcommand\Orbffcr{\Orbit_{\mrsfunc}(\mrsfunc,\singf)}
\newcommand\DiffMcr{\Diff\Mcr}
\newcommand\DiffIdMcr{\Diff_{\id}\Mcr}
\newcommand\DiffIdM{\Diff_{\id}(\manif)}

\newcommand\Stabfcr{\Stab\fcr}
\newcommand\StabIdfcr{\StabId\fcr}
\newcommand\Stabf{\Stab(\mrsfunc)}
\newcommand\StabIdf{\Stab_{\id}(\mrsfunc)}

\newcommand\pr{p}

\newcommand\crpt[1]{c_{#1}}
\newcommand\crcnt{n}

\newcommand\Ahom{H}

\newcommand\tdifM{\widetilde{\difM}}
\newcommand\Zid{Z}
\newcommand\DiffMn{\Diff(\manif,\crcnt)}
\newcommand\DiffIdMn{\Diff_{\id}(\manif,\crcnt)}
\newcommand\Kleinb{K}

\newcommand\Dpartf{\Diff(\partitf)}
\newcommand\DpartfId{\Diff_{\id}(\partitf)}
\newcommand\Dflow{\Diff(\flow)}
\newcommand\DflowId{\Diff_{\id}(\flow)}
\newcommand\Cflow{\End(\flow)}
\newcommand\CflowId{\End_{\id}(\flow)}

\newcommand\tStabIdf{\widetilde{\Stab}_{\id}(\mrsfunc)}

\newcommand\apoint{{\rm point}}

\newcommand\Dinvflow{\Gamma}
\newcommand\Dinvplus{\Dinvflow^{+}}

\newcommand\invol{\xi}

\newcommand\excrlev{L}
\newcommand\ashift{\alpha}

\newcommand\Cyl{B}
\newcommand\restM{L}
\newcommand\adifM{g}
\newcommand\cdifM{\widehat{\difM}}
\newcommand\anbh{V}
\newcommand\bdif{q}
\newcommand\Seps{S_{\eps}}
\newcommand\Stwoeps{S_{2\eps}}
\newcommand\Interv{I} 
\newcommand\larc{l}

\newcommand\omanif{M}
\newcommand\nmanif{N}
\newcommand\omrsfunc{f}
\newcommand\nmrsfunc{g}

\newcommand\difoM{u}
\newcommand\difnM{v}

\newcommand\DiffId{\Diff_{\id}}

\newcommand\partitof{\Delta_{\omrsfunc}}
\newcommand\partitnf{\Delta_{\nmrsfunc}}

\newcommand\DiffnM{\Diff(\nmanif)}

\newcommand\DpartIdnf{\DiffId(\partitnf)}
\newcommand\DpartIdf{\DiffId(\partitf)}

\newcommand\Dpartfplus{\Diff^{+}(\partitf)}

\newcommand\Dflowplus{\Diff^{+}(\flow)}

\newcommand\oDflow{\widetilde{\Diff}(\flow)}
\newcommand\oDflowId{\widetilde{\Diff}_{\id}(\flow)}
\newcommand\oCflow{\widetilde{\End}(\flow)}
\newcommand\oCflowId{\widetilde{\End}_{\id}(\flow)}

\newcommand\oDpartofplus{\widetilde{\Diff}^{+}(\partitof)}
\newcommand\Dpartnfplus{\Diff^{+}(\partitnf)}
\newcommand\oDflowplus{\widetilde{\Diff}^{+}(\flow)}

\newcommand\StabIdnf{\StabId(\nmrsfunc)}

\newcommand\stmr{\Cinf(\omanif,\RRR)}

\newcommand\somr{\Cinf(\omanif,\RRR)}

\newcommand\comptstmr[1]{E_{#1}} 
\newcommand\perfunc{\theta}
\newcommand\Dinv[1]{\Gamma_{#1}}
\newcommand\isolift{\rho}

\newcommand\nfoliat{\Delta}
\newcommand\Dfol{\Diff(\nfoliat)}
\newcommand\DfolId{\Diff_{\id}(\nfoliat)}

\newcommand\Cfol{\End(\nfoliat)}
\newcommand\CfolId{\End_{\id}(\nfoliat)}

\newcommand\twgr{\mathcal{J}}
\newcommand\twgrId{\twgr_{0}}

\newcommand\Reebf{\Gamma(\mrsfunc)}
\newcommand\aReebf{\Gamma_{H}(\mrsfunc)}
\newcommand\prReebf{\hat{\mrsfunc}}
\newcommand\prReeb{\pr_{\mrsfunc}}

\newcommand\Reebfb{\Gamma(\bar\mrsfunc)}

\newcommand\AutfR{\Aut(\Reebf)}
\newcommand\AutfRHId{\Aut_{\id}(\Reebf)}
\newcommand\actSRall{\lambda}
\newcommand\actSR{\widetilde{\actSRall}}
\newcommand\edgeR{e}
\newcommand\intecnt{l}

\newcommand\totsupp{T}

\newcommand\crv{\gamma}
\newcommand\tcrv{\crv} 
\newcommand\Dtw{\tau}
\newcommand\tDtw{\Dtw} 

\newcommand\izerocr{j_{0}}
\newcommand\izero{i_{0}}
\newcommand\ione{i_{1}}

\newcommand\Dvert{$\partial$\,}
\newcommand\Evert{$e$\,}
\newcommand\Cvert{$c$\,}
\newcommand\generic{generic}
\newcommand\simple{simple}

\newcommand\typeA{{\rm($A$)}}
\newcommand\typeB{{\rm($B$)}}
\newcommand\typeC{{\rm($C$)}}
\newcommand\typeDk{{\rm($D$)}}
\newcommand\typeEk{{\rm($E$)}}
\newcommand\CyclesR{C}

\newcommand\atom{A}
\newcommand\datom{A'}

\newcommand\Homol{$H_1$}

\newcommand\Graph{\Gamma}

\newcommand\sphere{S^2}
\newcommand\torus{T^2}
\newcommand\parthom{\partial}
\newcommand\kerAut{A}

\newcommand\prF{p}

\newcommand\MobiusBand{M\text{\"o}}
\newcommand\prjplane{\RRR P^2}

\newcommand\tPsp{\widetilde{\Psp}}
\newcommand\Rdefect{m}
\newcommand\RdefectC{\Rdefect_{c}}
\newcommand\RdefectE{\Rdefect_{e}}
\newcommand\minGr{\Graph_{\min}}

\newcommand\rankpiO{k}
\newcommand\keractSRall{\ker\actSRall}
\newcommand\keractSRallId{(\ker\actSRall)_{\id}}
\newcommand\vectOm{v}
\newcommand\crcomp{K}

\newcommand\acycle{\gamma}

\newcommand\tomega{\widetilde{\omega}}

\newcommand\LLL{{\rm(LL)}}

\newcommand\condHHisID{{\rm(A)}}
\newcommand\condHkeespEdges{{\rm(B)}}
\newcommand\condHkeepsCycles{{\rm(C)}}

\newcommand\Kends[1]{[#1]_{\crcomp}}

\newcommand\tgfunc{\widetilde{g}}
\newcommand\thfunc{\widetilde{h}}

\newcommand\VMd{\Gamma_{\partial}(\manif)}
\newcommand\VMdcr{\Gamma_{\partial}(\manif,\singf)}

\newcommand\codimf{\mu(\mrsfunc)}
\newcommand\codimfcr{\mu(\mrsfunc,\singf)}

\newcommand\Mman{M}
\newcommand\Nman{N}
\newcommand\smmn{\Cinf(\Mman,\Nman)}

\newcommand\ssmftn{\Sinf(\Mman;\mrsfunc^{*}T\Nman)}

\newcommand\ttorus{\widetilde{T}}
\newcommand\Reebfn{\Gamma(\mrsfunc_{n})}
\newcommand\aReebfn{\Gamma_{H}(\mrsfunc_{n})}

\newcommand\KR{KR}

\newcommand\ZNbh{\mathcal{Z}}
\newcommand\feNbh{\mathcal{N}}

\newcommand\ael{a}

\newcommand\del{d}
\newcommand\fel{f}
\newcommand\gel{g}
\newcommand\hel{h}

\newcommand\muel{\mu}

\newcommand\rrel{\lambda}
\newcommand\xrel{\xi}

\newcommand\xel{x}
\newcommand\yel{y}

\newcommand\mhat[1]{\hat{#1}}

\newcommand\tgel{\mhat{\gel}}
\newcommand\thel{\mhat{\hel}}

\newcommand\trrel{\mhat{\rrel}}
\newcommand\txrel{\mhat{\xrel}}

\newcommand\txel{\mhat{\xel}}

\newcommand\FFr{F}
\newcommand\GFr{G}
\newcommand\HFr{H}

\newcommand\UFr{U}

\newcommand\XFr{X}

\newcommand\Bmp{B}

\newcommand\Dm[1]{D#1}
\newcommand\Dmp[2]{\Dm{#1}(#2)}
\newcommand\Dmpv[3]{\Dm{#1}(#2\,;\,#3)}

\newcommand\Invm[1]{V\!#1}

\newcommand\Grp{G} 
\newcommand\Hgrp{R} 
\newcommand\actGX{\alpha}

\newcommand\un{e}
\newcommand\xpnt{\fel}

\newcommand\codm[2]{\mu_{#1}(#2)}

\newcommand\ImDx[1]{\mathrm{im}\,\Dm{\Gact}(#1)}
\newcommand\ImD{\mathrm{im}D}

\newcommand\grgr{*}
\newcommand\grtgr{\,\overrightarrow{{}_{\grgr}}\,}
\newcommand\tgrgr{\,\overleftarrow{{}_{\grgr}}\,}
\newcommand\grx{{\cdot}}
\newcommand\grtx{\,\overrightarrow{{}_{\grx}}\,}

\newcommand\Gmult{\mu}
\newcommand\tGmult{\overrightarrow{\Gmult}}
\newcommand\trGmult{\overleftarrow{\Gmult}} 

\newcommand\Ginv{\nu}

\newcommand\Gact{\alpha}
\newcommand\tGact{\overrightarrow{\Gact}}

\newcommand\Hact{\beta}

\newcommand\Hmult{\nu}

\newcommand\GHact{\chi}

\newcommand\UUU{U}

\newcommand\UG{\UUU_{G}}
\newcommand\UX{\UUU_{\XFr}}
\newcommand\UR{\UUU_{R}}

\newcommand\FX{F_{\XFr}}
\newcommand\FG{F_{G}}

\newcommand\fcod{n}
\newcommand\Asur{\chi}

\newcommand\Rf{\RRR^{\fcod}}

\newcommand\phiel{\phi}
\newcommand\tdop[1]{\phiel_{#1}}
\newcommand\sprod[2]{\langle\,#1\,,\,#2\,\rangle}

\newcommand\sectGname[1]{#1^{g}}
\newcommand\sectRname[1]{#1^{r}}

\newcommand\sectGHX{S}
\newcommand\sectGHXG{\sectGname{\sectGHX}}
\newcommand\sectGHXR{\sectRname{\sectGHX}}

\newcommand\Lsect{L}
\newcommand\LsectG{\sectGname{\Lsect}}
\newcommand\LsectR{\sectRname{\Lsect}}

\newcommand\Msect{M}
\newcommand\MsectG{\sectGname{\Msect}}
\newcommand\MsectR{\sectRname{\Msect}}

\newcommand\Bsect{B}
\newcommand\BsectG{\sectGname{\Bsect}}
\newcommand\BsectR{\sectRname{\Bsect}}

\newcommand\Asect{\sectGHX}
\newcommand\AsectG{\sectGname{\Asect}}
\newcommand\AsectR{\sectRname{\Asect}}

\newcommand\UNbh{U}

\newcommand\WNbh{W}

\newcommand\TfX{T_{\fel}\XFr}
\newcommand\TeG{T_{\un}\Grp}
\newcommand\isot{\beta}

\newcommand\LsectRg[1]{\Gamma_{#1}}

\newcommand\hrrel{\widetilde\rrel}
\newcommand\hgel{\widetilde\gel}

\newcommand\Ract[2]{#2\oplus#1}

\newcommand\emb{W}
\newcommand\invemb{V}

\newcommand\LsectM{\Lsect^{m}}

\newcommand\VMPD{\Gamma_{\partial}(\Mman,\Sigma)}
\newcommand\DMS{\Diff(\Mman,\Sigma)}

\newcommand\sgradf{G}

\title[Stabilizers and orbits of Morse functions]{Homotopy types of stabilizers and orbits \\ of Morse functions on surfaces}

\author{Sergey Maksymenko}
\email{maks@imath.kiev.ua}
\date{}
\address{Topology Department, Institute of Mathematics, Ukrainian National Academy of Science, Te\-re\-shchen\-kiv\-ska str. 3, 01601 Kyiv, Ukraine}

\begin{document}

\begin{abstract}
Let $\manif$ be a smooth compact surface, orientable or not, with boundary or without it, $\Psp$ either the real line $\RRR^1$ or the circle $\aCircle$, and $\DiffM$ the group of diffeomorphisms of $\manif$ acting on $C^{\infty}(\manif,\Psp)$ by the rule $\difM\cdot\mrsfunc=\mrsfunc\circ\difM^{-1}$ for $\difM\in\DiffM$ and $\mrsfunc\in C^{\infty}(\manif,\Psp)$.
Let $\mrsfunc:\manif\to\Psp$ be an arbitrary Morse mapping, $\singf$ the set of critical points of $\mrsfunc$, $\DiffMcr$ the subgroup of $\DiffM$ preserving $\singf$, and $\Stabf$, $\Stabfcr$, $\Orbf$, and $\Orbfcr$ the stabilizers and the orbits of $\mrsfunc$ with respect to $\DiffM$ and $\DiffMcr$.
In fact $\Stabf=\Stabfcr$.

In this paper we calculate the homotopy type of $\Stabf$, $\Orbf$ and $\Orbfcr$.
It is proved that except for few cases the connected components of $\Stabf$ and $\Orbfcr$ are contractible,
$\pi_k\Orbf=\pi_k\manif$ for $k\geq 3$, $\pi_2\Orbf=0$, and $\pi_1 O(f)$ is an extension of $\pi_1\DiffM\oplus\ZZZ^{k}$ (for some $k\geq0$)  with a (finite) subgroup of the group of automorphisms of the Kronrod-Reeb graph of $\mrsfunc$.

We also generalize the methods of F.~Sergeraert to prove that a finite codimension orbit of a tame smooth action of a tame Lie group on a tame \Frechet\ manifold is a tame \Frechet\ manifold itself.
In particular, this implies that $\Orbf$ and $\Orbfcr$ are tame \Frechet\ manifolds.
\end{abstract}

\maketitle


\section{Introduction}\label{sect:intro}
The study of homotopy properties of Morse functions space on surfaces was stimulated
in recent years by the applications in symplectic topology and Hamiltonian
systems~\cite{BolsFom, Kudr, Sharko, Kulinich, Maks:PhD, Maks:PathComp, SaekiIkegami, SharkoUMZ}.
Connected components of this space were described in~\cite{Kudr, Sharko,Maks:PhD,Maks:PathComp}.
The cobordism group of Morse functions was calculated in~\cite{SaekiIkegami}.

In this paper we describe the homotopy types of stabilizers and orbits of Morse functions on
compact surfaces under the action of diffeomorphism groups of these surfaces.
Applications of these results will be published elsewhere.

Let $\manif$ be a smooth ($C^{\infty}$) compact connected manifold, orientable or not, with boundary or without it.
Let also $\Psp$ be either $\RRR$ or $\aCircle$.
The group $\DiffM$ of diffeomorphisms of $\manif$ naturally acts on $\smone$
by the following rule: if $\cmap\in\DiffM$ and $\mrsfunc\in\smone$, then
\chEq\begin{equation}\label{equ:action-DM}
\cmap\cdot\mrsfunc = \mrsfunc\circ\cmap^{-1}.
\end{equation}
For $\mrsfunc\in\smone$ let
$\Stabf = \{\difM\in\DiffM \, | \, \mrsfunc\circ\difM = \mrsfunc \}$
be the {\em stabilizer} and
$\Orbf=\{\mrsfunc\circ\difM \, | \, \difM\in\DiffM\}$
the {\em orbit} of $\mrsfunc$ under this action.

Let $\singf$ be the set of critical points of $\mrsfunc$ and $\DiffMcr$ the group of diffeomorphisms $\difM$ of $\manif$ such that $\difM(\singf)=\singf$.
Then the stabilizer $\Stabfcr$ and the orbit $\Orbfcr$ of $\mrsfunc$ under the restriction of the above action to $\DiffMcr$ are well defined.
Since $\Stabf\subset\DiffMcr$, we get $\Stabfcr = \Stabf$.

We endow $\DiffM$ and $\smone$ with $C^{\infty}$ topologies, and their subspaces $\Stabf$,   $\Orbf$, and $\Orbfcr$ with the induces ones.

Let $\mrsfunc:\manif\to\Psp$ be a smooth mapping.
By a {\em level-set\/} or {\em point-inverse\/} of $\mrsfunc$ we mean a set $\mrsfunc^{-1}(c)$, where $c\in\Psp$.
This set is {\em critical\/} if it contains a critical point of $\mrsfunc$, otherwise, $\mrsfunc^{-1}(c)$ is {\em regular}.
Similarly, a connected component $\omega$ of $\mrsfunc^{-1}(c)$ is called {\em critical\/} if it contains a critical point of $\mrsfunc$; otherwise, $\omega$ is {\em regular}.
\begin{defn}\label{defn:cond-i-ii-for-f}
Let $\mrsfunc:\manif\to\Psp$ be a smooth mapping satisfying the following two conditions:

{\rm(i)}
$\mrsfunc$ is constant on every connected component of $\partial\manif$ and

{\rm(ii)}
critical points of $\mrsfunc$ are isolated and belong to the interior of $\manif$.

Then $\mrsfunc$ will be called {\em generic} if every level-set of $\mrsfunc$
contains at most one critical point;
$\mrsfunc$ is {\em \simple} if every critical component
of a level-set of $\mrsfunc$ contains precisely one critical point;
$\mrsfunc$ is {\em Morse} provided all critical points of $\mrsfunc$ are non-degenerate.
\end{defn}
Evidently, every generic mapping is \simple.

Let us fix once and for all some orientation of $\Psp$.
Then for every Morse mapping $\mrsfunc:\manif\to\Psp$ the indices of its (non-degenerate) critical points are well-defined.

Throughout the paper we will use the following notations:
$\aCircle=\RRR/\ZZZ$; $\Interv = [0,1]$;
$D^2$ is the closed unit disk in $\RRR^2$;
$\sphere$ the unit sphere in $\RRR^3$;
$\prjplane$ the real projective plane;
$\MobiusBand$ the M\"obius band;
$\torus=\aCircle\times\aCircle$ the $2$-torus;
and $\Kleinb$ the Klein bottle.

The main results of this paper are Theorems~\ref{th:StabIdf-hom-type},~\ref{th:hom-gr-orbits-c1}, and~\ref{th:hom-gr-orbits-c10} below.
For each of them we will give at first a brief idea of proof.
This will help the reader interested only in some separate result of the paper extract the necessary proof.

\begin{thm}\label{th:StabIdf-hom-type}
Let $\manif$ be a smooth compact connected surface and $\mrsfunc:\manif\to\Psp$ a Morse mapping.
Denote by $\StabIdf$ the subset of $\Stabf$ consisting of diffeomorphisms that are isotopic  in $\Stabf$ to $\id_{\manif}$, and endow $\StabIdf$ with the induced $C^{\infty}$ topology.
Suppose that either $\mrsfunc$ has at least one critical point of index $1$ or $\manif$ is non-orientable, then $\StabIdf$ is contractible.
Otherwise, $\StabIdf$ is homotopy equivalent to $\aCircle$.
\end{thm}
\begin{rem}
\rm
Denote by $\Stabf^{\infty}$ and $\Stabf^{0}$ the space $\Stabf$ with $C^{\infty}$ and $C^{0}$ topologies respectively.
By the definition $\StabIdf$ consists of diffeomorphisms isotopic in $\Stabf$ to $\id_{\manif}$. Hence $\StabIdf$ is the identity path-component of $\Stabf^{0}$.
Let $\StabIdf'$ be the identity path-component of $\Stabf^{\infty}$.
Since the identity mapping $\id:\Stabf^{\infty}\to\Stabf^{0}$ is continuous, we have that 
$\StabIdf'\subset\StabIdf$.

On the other hand, by Theorem~\ref{th:StabIdf-hom-type}, $\StabIdf$ is connected when endowed with $C^{\infty}$ topology, whence $\StabIdf\subset\StabIdf'$. Thus $\StabIdf$ is also the connected component of $\Stabf^{\infty}$.
\end{rem}

Theorem~\ref{th:StabIdf-hom-type} is essentially new only for non-orientable surfaces.
For orientable ones it is a simple corollary of~\cite{Maks:Shifts}.
Notice that if $\manif$ is orientable then there is a flow on $\manif$ such that the level-sets of $\mrsfunc$ consist of full trajectories of $\flow$.
Let $\Dflow$ be the group of diffeomorphisms preserving every trajectory of $\flow$ and $\DflowId$ be its identity path-component.
Then it is easy to show that $\StabIdf=\DflowId$, see Lemma~\ref{lm:Dpart-Stabf} and (1) of Lemma~\ref{lm:exist-Fld}.
In~\cite{Maks:Shifts} the author gave the condition on a flow $\flow$ on a manifold $\manif$ for $\DflowId$ to be either contractible of homotopy equivalent to $\aCircle$: the flow $\flow$ must be linear (in some local coordinates) near its fixed points.
Since $\mrsfunc$ is a quadratic form near its critical points (by Morse lemma), we can choose $\flow$ to satisfy that condition.

If $\mrsfunc$ has at least one critical point of index $1$, then $\flow$ has non-closed trajectory, whence
by~\cite{Maks:Shifts} $\DflowId$ is contractible. Hence so is $\StabIdf$.

Otherwise $\mrsfunc$ has no critical points of index $1$, all non-constant trajectories of $\flow$ are closed, whence again by~\cite{Maks:Shifts} $\DflowId=\StabIdf$ is homotopy equivalent to $\aCircle$.

If $\manif$ is non-orientable, we reduce the situation to the orientable double covering of $\manif$.
Thus we can assume that $\manif$ is orientable but equipped with an orientation reversing involution $\xi$ without fixed points such that $\mrsfunc\circ\xi=\mrsfunc$.
We also have to study instead of $\StabIdf$ the group $\tStabIdf$ of $\xi$-equivariant diffeomorphisms.
Then we choose $\flow$ so that $\xi$ permutes trajectories of $\flow$ changing their orientation.
In this situation we show the contractibility of $\tStabIdf$, see Section~\ref{subsect:SkewSymmetrFlows}.

\begin{thm}\label{th:hom-gr-orbits-c1}
Let $\manif$ be a smooth compact connected surface and $\mrsfunc:\manif\to\Psp$ be a Morse mapping having at least one critical point of index $1$.
Let $\Orbff$ and $\Orbffcr$ be the corresponding path-components of $\Orbf$ and $\Orbfcr$ in $C^{\infty}$-topology containing $\mrsfunc$.
Then

{\rm(1)} $\Orbffcr$ is contractible;

{\rm(2)} $\pi_i\Orbff \approx \pi_i\manif$ for $i\geq3$ and
$\pi_2\Orbff=0$.
In particular, $\Orbff$ is aspherical provided $\manif$ is.
Moreover,
$\pi_1\Orbff$ is included in the following exact sequence:
\chEq\begin{equation}\protect\label{equ:pi1Orbf-ex-seq}
0 \to \pi_1\DiffM \oplus \ZZZ^{\rankpiO} \to \pi_1\Orbff \to \grp \to 1,
\end{equation}
where $\grp$ is a finite group and $\rankpiO\geq0$.

Let $\crpt{0}$, $\crpt{1}$, $\crpt{2}$ be the numbers of critical points of $\mrsfunc$ of the corresponding indices.
Then $\rankpiO\leq\bar\rankpiO$ where $\bar\rankpiO$ is shown in Table~\ref{tbl:contractions}.
Moreover, if $\manif$ is of type 1 or 2 of this table and $\mrsfunc$ is {\em \,simple\,}, then in fact $\rankpiO=\bar\rankpiO$.

\chTbl\begin{table}[ht]
\caption{}
\protect\label{tbl:contractions}
\begin{tabular}{|c|p{0.6\textwidth}|c|}\hline
& \centering $\manif$ & $\bar\rankpiO$  \\ \hline
1 & $\sphere$, $D^2$, $\aCircle\times\Interv$, $\torus$, $\prjplane$ with or without holes & $\crpt{1}-1$ \\ \hline
2 & $\manif$ is orientable and is not of type 1 & $\crpt{1}+\chi(\manif)=\crpt{0}+\crpt{2}$ \\ \cline{1-2}
3 & $\manif$ is non-orientable and is not of type 1 &  \\ \hline
\end{tabular}
\end{table}

{\rm(3)} Suppose that $\mrsfunc$ is {\em\,generic\,}.
Then the group $\grp$ in Eq.~\eqref{equ:pi1Orbf-ex-seq}
is trivial, whence
$\pi_1\Orbff\approx \pi_1\DiffM \oplus \ZZZ^{\rankpiO}$.
In particular, $\pi_1\Orbff$ is abelian.
The homotopy type of $\Orbff$ for this case is given in Table~\ref{tbl:hom_groups_of_orbits_c1}.
\chTbl\begin{table}[ht]
\caption{}
\protect\label{tbl:hom_groups_of_orbits_c1}
\begin{tabular}{|c|c|}\hline
$\manif$ & Homotopy type of $\Orbff$ \\ \hline\hline
$\sphere$, $\prjplane$  & $SO(3)\times(\aCircle)^{\crpt{1}-1}$ \\ \hline
$D^2$, $\aCircle\times\Interv$, $\MobiusBand$ & $(\aCircle)^{\crpt{1}}$ \\ \hline
$\torus$ & $(\aCircle)^{\crpt{1}+1}$ \\ \hline
$\Kleinb$ & $(\aCircle)^{\rankpiO+1}$ \\ \hline
other cases  & $(\aCircle)^{\rankpiO}$ \\ \hline
\end{tabular}
\end{table}
\end{thm}
Let us describe the plan of the proof of this theorem.
Following methods of F.~Sergeraert~\cite{Sergeraert} we will show in Appendix~\ref{app:Frech_struct} that the natural projection $p:\DiffM\to\Orbf$ is a locally-trivial principal $\Stabf$-fibration.
By Theorem~\ref{th:StabIdf-hom-type} we also have that $\StabIdf$ is contractible.
Then from the exact homotopy sequence of the fibration above we get isomorphisms between the higher homotopy groups of $\DiffM$ and $\Orbf$.
The similar statement holds for $\DiffMcr$ and $\Orbfcr$.
Then triviality of groups $\pi_{i}\Orbff$ and $\pi_{i}\Orbffcr$ for $i\geq2$ in (1) and the first sentence of (2) of Theorem~\ref{th:hom-gr-orbits-c1} are direct corollaries of the classification of the homotopy types of $\DiffIdM$ for compact surfaces $\manif$, see Table~\ref{tbl:hom-types}. Thus it remains to calculate $\pi_1\Orbff$ and prove that $\pi_1\Orbffcr=0$.

\smallskip
First consider the group $\pi_1\Orbff$.
The exact sequence~\eqref{equ:pi1Orbf-ex-seq} shows that $\pi_1\Orbff$ depends on three terms.
Let us briefly explain them.

Let $\mrsfunc_t:\manif\to\Psp$, $t\in[0,1]$, be a loop in $\Orbff$ based at $\mrsfunc=\mrsfunc_0=\mrsfunc_1$.
Using the fibration property we lift $\mrsfunc_t$ to a path $\difM_t$ in $\DiffIdM$ such that $\mrsfunc_t=\mrsfunc\circ\difM_t$ and $\difM_0=\id_{\manif}$.

1) Suppose that $\difM_0=\difM_1=\id_{\manif}$, i.e. $\{\difM_t\}$ is a loop in $\DiffIdM$.
Let $p_{*}:\pi_1\DiffIdM \to \pi_1\Orbff$ be a natural homomorphism, then $\{\mrsfunc_t\}=p_{*}\{\difM_t\}$.
From aspherity of $\StabIdf$ we get that $p_{*}$ is injective.
This gives the term $\pi_1\DiffM$ in~\eqref{equ:pi1Orbf-ex-seq}.

2) Suppose that $\difM_1\not=\id_{\manif}$.
Since $\mrsfunc=\mrsfunc_1=\mrsfunc\circ\difM_1$, we see that $\difM_1\in\Stabf$.
Thus another type of generators of $\pi_1\Orbff$ arises from diffeomorphisms of $\Stabf$ that are isotopic to $\id_{\manif}$.
Evidently, they constitute the kernel of the natural homomorphism $\izero:\pi_0\Stabf\to\pi_0\DiffM$.

Notice that $\difM_1\in\Stabf$ yields an automorphism $\theta$ of the Kronrod-Reeb (\KR-) graph $\Reebf$ of $\mrsfunc$ so that the correspondence $\{\mrsfunc_t\} \mapsto \theta$ induces a homomorphism $\pi_1\Orbff\to\AutfR$, see Section~\ref{sect:ReebGr}.
The group $\grp$ in~\eqref{equ:pi1Orbf-ex-seq} is just its image.

3) Suppose that $\theta$ is the identity automorphism.
Let $\gamma$ be a regular component of a level-set of $\mrsfunc$.
Then there is a Dehn twist $\Dtw$ along $\gamma$ such that $\mrsfunc\circ\Dtw=\mrsfunc$, i.e. $\Dtw\in\Stabf$, see Section~\ref{sect:pi-0Dpartf}.
Thus $\Dtw$, as well as $\difM_1$, yields the identity automorphism of $\Reebf$.
We prove that $\difM$ is in fact generated by such Dehn twists (Proposition~\ref{pr:actSRall-Dpartf}).

Moreover, the connected components of the group of the above Dehn twists constitute a free abelian group $\twgr$ (Theorem~\ref{th:twgr-descr}).
Since $\difM_1$ is also isotopic to $\id_{\manif}$, it belongs to the subgroup of $\twgr$ generated by the ``relations'' in $\DiffM$ between the Dehn twists of $\twgr$, i.e. to the kernel of the natural homomorphism $\twgr \to \pi_0\DiffM$.
This kernel is precisely the group $\ZZZ^{\rankpiO}$ of~\eqref{equ:pi1Orbf-ex-seq}.

Finally, the sum $\pi_1\DiffIdM\oplus\ZZZ^{\rankpiO}$ in~\eqref{equ:pi1Orbf-ex-seq} is direct since the image $p_{*}(\pi_1\DiffIdM)$ is included in the center of $\pi_1\Orbff$, see Lemma~\ref{lm:prop-pi-1Dm}.

\smallskip
Consider now the group $\pi_1\Orbffcr$.
By Theorem~\ref{th:StabIdf-hom-type} the projection $\DiffMcr\to\Orbfcr$ is also a locally trivial fibration.
Hence a loop $\{\mrsfunc_t\}\in\pi_1\Orbffcr$ based at $\mrsfunc$ can be represented in the form $\mrsfunc_t=\mrsfunc\circ\difM_t$, where $\difM_0=\id_{\manif}$, $\difM_1\in\Stabf$ and $\difM_t\in\DiffIdMcr$.
Let $\theta$ be the automorphism of the \KR-graph $\Reebf$ of $\mrsfunc$ induced by $\difM_1$.
Since $\mrsfunc_t$ is also a loop in $\Orbff$, it is generated by loops of the types 1)-3) as above.
We prove that in our case all such loops are trivial.

1) If $\difM_1=\id_{\manif}$, then the loop $\{\difM_t\}$ is trivial in $\DiffIdMcr$ due to the contractibility of $\DiffIdMcr$ (Lemma~\ref{lm:DiffIdMcr-is-contractible}).

2) Since every $\difM_t$ preserves critical points of $\mrsfunc$, it follows that $\theta$ is trivial (Proposition~\ref{pr:actSRall-Dpartf} and Lemma~\ref{lm:Dpart-Stabf}).

3) Hence $\difM_1$ is generated by the Dehn twists along regular level-sets of $\mrsfunc$.
We prove that $\difM_1$ preserves every connected component of every level-set of $\mrsfunc$ (Theorem~\ref{th:ends} and Proposition~\ref{pr:actSRall-Dpartf}).
This is the most difficult part of Theorem~\ref{th:hom-gr-orbits-c1}.
We then deduce that $\mrsfunc_t$ is a trivial loop.
Thus $\pi_1\Orbffcr=0$.

All other statements of Theorem~\ref{th:hom-gr-orbits-c1} are based on the study \KR-graph of $\mrsfunc$, see Section~\ref{sect:proof-Th3}.

\begin{thm}\label{th:hom-gr-orbits-c10}
Let $\mrsfunc:\manif\to\Psp$ be a Morse mapping having no critical points of index $1$.
Then $\mrsfunc$ can be represented in the following form
$$
\mrsfunc = \pr\circ\tmrsfunc:
\manif
\stackrel{\tmrsfunc}{\longrightarrow}
\tPsp
\stackrel{\pr}{\longrightarrow}
\Psp,
$$
where $\tmrsfunc$ is one of the mappings ``of type \typeA-\typeEk''
shown in Table~\ref{tbl:hom_groups_of_orbits},
$\tPsp$ is either $\RRR$ or $\aCircle$,
and $\pr$ is either a covering map or an embedding.
The homotopy types of $\Orbff$ and $\Orbffcr$
depend only on $\tmrsfunc$ and are
given in Table~\ref{tbl:hom_groups_of_orbits}.

{
\footnotesize
\chTbl\begin{table}[ht]
\caption{}\protect\label{tbl:hom_groups_of_orbits}
\begin{tabular}{|c|l|l|c|c|c|c|c|}\hline
Type & \multicolumn{2}{|c|}{$\tmrsfunc:\manif\to\tPsp$} & $\crpt{0}$ & $\crpt{1}$ & $\crpt{2}$ &
$\Orbff$ & $\Orbffcr$ \\ \hline
\typeA & $\sphere\to\RRR$ &  \ $\tmrsfunc(x,y,z)=z$                   & $1$  & $0$ & $1$ &  $\sphere$ & $\apoint$ \\ \hline
\typeB & $D^2\to\RRR$ &  \ $\tmrsfunc(x,y)=x^2+y^2$               & $1$ & $0$ & $0$ & \multicolumn{2}{|c|}{$\apoint$} \\ \hline
\typeC & $\aCircle\times\Interv\to\RRR$ &  \ $\tmrsfunc(\phi,t)=t$     &  $0$ & $0$ & $0$ & \multicolumn{2}{|c|}{$\apoint$}  \\ \hline
\typeDk & $\torus\to\aCircle$ & \ $\tmrsfunc(x,y) = x$,   &  $0$ & $0$ & $0$  & \multicolumn{2}{|c|}{$\aCircle$}  \\ \hline
\typeEk & $\Kleinb\to\aCircle$ & \ $\tmrsfunc(\{x\},\{y\}) = \{2x\}$
 &  $0$ & $0$ & $0$  & \multicolumn{2}{|c|}{$\aCircle$}  \\ \hline
\end{tabular}
The Klein bottle $\Kleinb$ is regarded here as the factor space of $\torus$ by the involution
$\invol(x,y)=(x+1/2,-y)$.
\end{table}
}
\end{thm}

\begin{rem}
\rm We do not consider another natural action of the group $\DiffM\times\DiffP$ on $\smone$
(see e.g.~\cite{Sergeraert}).
The comparison between the stabilizers and orbits of functions under this action and under the action of the group $\DiffM$ will appear in another paper (see~\cite{Maks:StOrb}).
\end{rem}

The paper is arranged in the following way.
In Section~\ref{sect:orbits} we formulate Theorem~\ref{th:loc-triv-fibering} claiming that $\Orbf$ and $\Orbfcr$ are \Frechet\ manifolds and that projections $\DiffM\to\Orbf$ and $\DiffMcr\to\Orbfcr$ are locally trivial fibrations with the same fiber $\Stabf$.
Hence $\Orbf$ and $\Orbfcr$ have the homotopy types of CW-complexes, see~\cite{Palais}, and we also get exact sequences of homotopy groups of these fibrations.
These results seem to be more or less known, but the author has not found the proof in the literature.
For the sake of completeness we prove Theorem~\ref{th:loc-triv-fibering} in Section~\ref{sect:proof_th_loctriv}.
The method of proof is mostly an extension of results of F.~Sergeraert~\cite{Sergeraert}, see also~\cite{Poenaru, Hendriks}.

In Section~\ref{sect:orbits} we also recall a description of the homotopy types of groups $\DiffIdM$ for compact surfaces.

Section~\ref{sect:prelim} contains definitions of the Kronrod-Reeb graph $\Reebf$ and the foliation $\partitf$ on $\manif$ induced by a smooth function $\mrsfunc$ satisfying conditions (i) and (ii) of Definition~\ref{defn:cond-i-ii-for-f}.
We also establish some simple properties of them.

Section~\ref{sect:sm_shifts}.
We briefly formulate the results of~\cite{Maks:Shifts} concerning smooth shift along trajectories of flows and extend them to {\em skew-symmetric} flows on the oriented double covering of a non-orientable manifold (Theorem~\ref{th:CDfolId-contr}).

Section~\ref{sect:proof-th:StabIdf-hom-type-or}.
Using the results of Section~\ref{sect:sm_shifts} we prove Theorem~\ref{th:StabIdf-hom-type}.
The proof does not use the results of Section~\ref{sect:orbits} and can be read independently.
Theorem~\ref{th:StabIdf-hom-type} will allow to compute the groups $\pi_k\Orbff$ and $\pi_k\Orbffcr$ for $k\geq2$ (see Section~\ref{sect:proof-Th3}). Therefore in the next three sections we concentrate upon the fundamental groups of orbits.

Section~\ref{sect:pi-0Dpartf}.
It deals with the group $\Dpartf$ of diffeomorphisms preserving every leaf of the foliation $\partitf$ on $\manif$ defined by a Morse function $\mrsfunc$.
We prove that $\pi_0\Dpartf$ is an abelian group generated by the Dehn twists along regular components of level-sets of $\mrsfunc$ and that except for few cases this group is free abelian.

Section~\ref{sect:ends}. We prove Theorem~\ref{th:ends} which guarantees an exactness of the sequence Eq.~\eqref{equ:pi1Orbf-ex-seq}.

Section~\ref{sect:ijzero}.
The kernels of the natural homomorphisms of $\pi_0\Stabf$ to $\pi_0\DiffM$ and $\pi_0\DiffMcr$ are studied.
They are homomorphic images of $\pi_1\Orbff$ and $\pi_1\Orbffcr$ respectively.

Section~\ref{sect:proof-Th3}. We prove Theorem~\ref{th:hom-gr-orbits-c1}
and in particular show how to calculate the number $\rankpiO$ in
Eq.~\eqref{equ:pi1Orbf-ex-seq}.

Section~\ref{sect:proof-orb-gr} contain the proof of Theorem~\ref{th:hom-gr-orbits-c10}.

Finally in Appendix~\ref{app:Frech_struct} we prove Theorem~\ref{th:loc-triv-fibering}.


\section{Orbits $\Orbf$ and $\Orbfcr$}\label{sect:orbits}
The following theorem can be established by the methods similar to F.~Sergeraert~\cite{Sergeraert}, see also~\cite{Poenaru,Hendriks}.
These methods are quite far from our main tecnique, therefore we shift the proof to Section~\ref{sect:proof_th_loctriv}.

\begin{thm}\label{th:loc-triv-fibering}
Let $\Mman$ be a smooth compact connected manifold of dimension $\dimM$ having $b$ connected components of $\partial\Mman$,
$\Psp$ either $\RRR$ or $\aCircle$, 
$\smoned$ the space of smooth mappings $\manif\to\Psp$ that are constant on the connected components of $\partial\manif$.
Let $\mrsfunc\in\smoned$ be a Morse mapping having $c$ critical points.
Endow  $\DiffM$, $\DiffMcr$, $\Stabf$, $\Orbf$, and $\Orbfcr$ with the corresponding $C^{\infty}$-topologies.
Then

{\rm (1)} 
$\Orbf$ and $\Orbfcr$ are \Frechet\ submanifolds of $\smoned$ of codimensions
$\codimf = c+b$  and \ $\codimfcr = c\dimM+c+b$ resp.
In particular, they have the homotopy types of CW-complexes, see e.g.{\rm~\cite{Palais}}.

{\rm (2)}
The projections $\DiffM\to\Orbf$ and $\DiffMcr\to\Orbfcr$
defined by $\difM \mapsto \mrsfunc\circ\difM$ are locally trivial principal $\Stabf$-fibrations.
\end{thm}

Consider some properties of the fibrations of (2) of Theorem~\ref{th:loc-triv-fibering}.
Let $\prF:\DiffM\to\Orbf$ be the projection defined by $p(\difM)=\mrsfunc\circ\difM$.
By Theorem~\ref{th:loc-triv-fibering} $\prF$ is a locally trivial fibration.
Consider the following part of its exact homotopy sequence:
$$
\pi_1(\DiffM,\id_{\manif}) \stackrel{\prF}{\longrightarrow}
\pi_1(\Orbf,\mrsfunc)  \stackrel{\parthom_1}{\longrightarrow}
\pi_0(\Stabf,\id_{\manif}),
$$
where $\partial_1$ is a boundary map.
Since $\Stabf$ is a topological group, so is 
$\pi_0(\Stabf,\id_{\manif}) \approx \Stabf/\StabIdf$.

\begin{lem}\label{lm:prop-pi-1Dm}
{\rm(1)}~$\parthom_1$ is a homomorphism of groups;

{\rm(2)}~$\prF(\pi_1\DiffM)$ is included in the center of $\pi_1\Orbf$.
\end{lem}
\begin{proof}
Let us briefly recall the construction of $\parthom_1$.
Let $\mrsfunc_t:\manif\to\Psp$ be a loop in $\Orbff$ based at $\mrsfunc$, i.e. $\mrsfunc_0=\mrsfunc_1=\mrsfunc$.
Then there exists a lifting of $\mrsfunc_t$ to $\DiffM$, i.e.$\!$ a path in $\adifM_t:\manif\to\manif$ such that $\mrsfunc_t=\mrsfunc\circ\adifM_t$, $\adifM_0=\id_{\manif}$, and $\mrsfunc \circ \adifM_1 = \mrsfunc$.
Then $\parthom_1(g) = [\adifM_1] \in \pi_0\Stabf$.

(1)
Let $g,h$ be two loops in $\Orbff$ based at $\mrsfunc$
and $\alpha,\beta\in\pi_1\DiffM$ their liftings.
Then the lifting of the loop $g\cdot h$ is $\alpha_1\circ\beta_1$
as it is illustrated in Figure~\ref{fig:composition}~a).
Thus $\parthom_1(g\cdot h) = \alpha_1\circ\beta_1 = \parthom_1(g)\  \parthom_1(h)$.
Similarly, it can be shown that
$\parthom_1(g^{-1}) = \parthom_1(g)^{-1}$.
So $\parthom_1$ is a homomorphism of groups.

(2)
Let $\adifM_t:\manif\to\manif$ be a loop in $\DiffM$ at $\id_{\manif}$.
Then $\mrsfunc\circ\adifM_t$ is a loop in $\Orbff$ at $\mrsfunc$ corresponding to $\prF(\adifM)$.

Let $\beta\in\pi_1\Orbf$ and $\difM_t:\manif\times\Interv\to\manif$ be a lifting of $\beta$ to $\DiffM$ such that $\difM_0=\id_{\manif}$, and $\mrsfunc=\mrsfunc\circ\difM_1$.
Thus $\beta(t,x) = \mrsfunc\circ\difM_t(x)$.

It follows that the commutator $[\prF(\adifM),\beta]$ has the form
$\mrsfunc\circ u_t$, where
\ $u = \adifM_{[0,1]} \ \difM_{[0,1]} \ (\difM_1\cdot\adifM_{[1,0]})\ \difM_{[1,0]}$ \
is a loop
in $\DiffM$ (see Figure~\ref{fig:composition} b)).

Then the following family of loops
$\adifM_{[0,1]}\ \difM_{[0,s]}\ (\difM_s \cdot \adifM_{[1,0]})\ \difM_{[s,0]}$
deforms $u$ to $\adifM_{[0,1]}\adifM_{[1,0]}$.
Thus $u$ is null-homotopic and so is $[\prF(\adifM),\beta]$.
\end{proof}

\SHOWFIG{
\chFig\begin{figure}[ht]
\begin{tabular}{ccc}
\includegraphics[height=2.8cm]{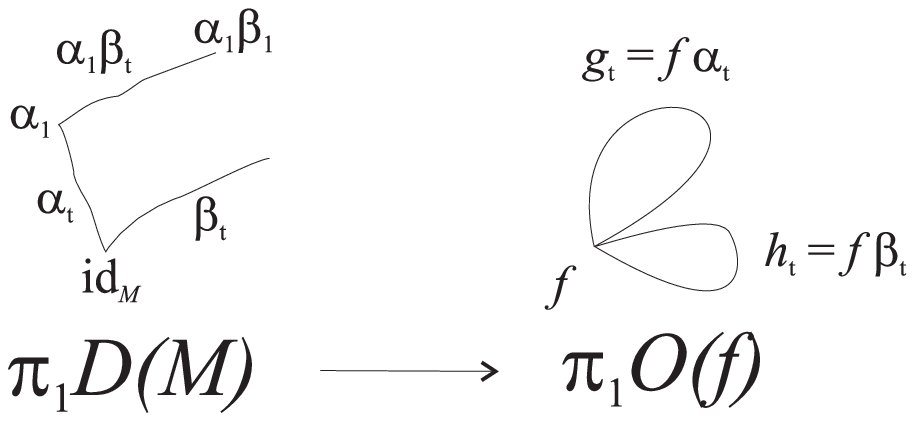}
& \qquad &
\includegraphics[height=2.8cm]{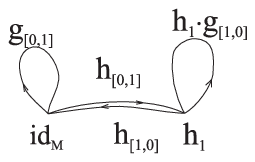} \\
a) & & b)
\end{tabular}
\caption{}\protect\label{fig:composition}
\end{figure}
}

\subsection{Homotopy type of $\DiffIdMn$ for compact surfaces} \label{subsect:HomGr}
Let $\manif$ be a compact surface and $x_1,\ldots,x_{\crcnt} \subset \interior\manif$ mutually distinct points.
Denote by $\DiffMn$ the group of diffeomorphisms of $\manif$ preserving the set of these points.
We endow this group with $C^{\infty}$-topology. 
Thus, if $\crcnt$ is the number of critical points of a function $\mrsfunc$, then $\DiffMn=\DiffMcr$ and $\DiffM = \Diff(\manif,0)$.

The homotopy types of the groups $\DiffIdMn$ for the case $\manif$ is oriented and closed were described in C.~J.~Earle and J.~Eells~\cite{EarleEells}.
For arbitrary compact surface $\manif$ the homotopy type of $\DiffIdM$ was studied by~C.~J.~Earle, A.~Schatz~\cite{EarleSchatz} and A.~Gramain~\cite{Gramain}.
These results can easily be extended to the groups $\DiffIdMn$ in sense that a puncture can be replaced by hole.

\begin{thm}\label{th:hom_type_DM}
Let $\manif$ be a connected compact surface and $\manif_{\bndcnt}$ a surface obtained by shrinking $\bndcnt$ connected components of $\partial\manif$ into points.
Then $\Diff_{\id}(\manif_{\bndcnt},\crcnt+\bndcnt)\sim\Diff_{\id}(\manif,\crcnt)$ (homotopy equivalent), so 
a puncture can be replaced by a hole, e.g.
$$
\Diff_{\id}(\sphere,2) \ \sim \ \Diff_{\id}(D^2,1) \ \sim \ \Diff_{\id}(\aCircle\times\Interv).
$$
Thus it suffices to know the homotopy types of $\DiffIdM$.
This information is collected in Table~\ref{tbl:hom-types}.
\chTbl\begin{table}[ht]
\caption {Homotopy type of $\DiffIdM$}
\protect\label{tbl:hom-types}
\begin{tabular}{|p{3.8cm}|c|c|c|c|c|} \hline
\centering{$\manif$} & 
$\sphere$, $\prjplane$ & 
$D^2$, $\aCircle\times\Interv$, $\MobiusBand$, $\Kleinb$ & 
$\torus$ & 
other cases \\ \hline
\centering{Homotopy type of $\DiffIdM$} & $SO(3)$ & $\aCircle$ & $\torus$ & point \\ \hline
\end{tabular}
\end{table}
\end{thm}
\begin{proof}
If $\manif$ has no punctures, i.e. $n=0$, then the homotopy type of $\DiffIdM$ was calculated in~\cite{EarleSchatz,Gramain} and is represented in Table~\ref{tbl:hom-types}.

Let us shrink one of the connected component of $\partial\manif$ into a point $x_1$ and denote the obtained surface by $\manif_1$.
We have to show that $\DiffId(\manif_1,1)\sim\DiffId(\manif)$.

Define an {\em evaluation map} $e:\DiffId(\manif_1)\to\manif_1$ at $x_1$ by $e(\difM)=\difM(x_1)$.
It is well known~\cite{FadellNeuwirth}, that $e$ is a locally trivial fibration with fiber $\DiffId(\manif_1,1)$.
Since we know the homotopy types of $\DiffId(\manif_1)$ and $\manif_1$, we will be able to establish the homotopy type of the fiber $\DiffId(\manif_1,1)$ via exact homotopy sequence of this fibration.

If $\manif_1=\sphere$ or $\prjplane$, so $\manif=D^2$ or $\MobiusBand$, then $\DiffId(\manif_1)\sim SO(3)$ and the corresponding exact sequence coincides with the exact sequence of the fibration of $SO(3)$ over $\sphere$ or $\prjplane$ by circles.
Hence $\DiffId(\manif_1,1)\sim\aCircle \sim \DiffId(\manif)$ by Table~\ref{tbl:hom-types}.  

Otherwise, $\manif_1$ and $\DiffId(\manif_1)$ are aspherical, whence so is $\DiffId(\manif_1,1)$ and therefore its homotopy type is determined only by the fundamental group.
The calculation of $\pi_1\DiffId(\manif_1,1)$ is based on the fact that the image of $e_{*}(\pi_1\DiffId(\manif_1))$ is included in the center of $\pi_1\manif_1$, e.g.~\cite[Lemma~1]{Gramain}.
If $\manif_1$ is neither $\torus$ nor $\aCircle\times\Interv$ nor $\MobiusBand$, then 
it is well known that $\pi_1\manif_1$ is centerless, whence we get an exact sequence:
$$
 \pi_2\manif_1 (=0) \ \stackrel{\partial_1}{\to} \
 \pi_1\DiffId(\manif_1,1) \ \to \
 \pi_1\DiffId(\manif_1)\  \stackrel{e_{*}}{\to} \  
 0 = (\text{center of $\pi_1\manif_1$}),
$$
implying $\DiffId(\manif_1,1)\sim\DiffId(\manif_1)$.

If $\manif_1$ is either $\torus$ or $\aCircle\times\Interv$ or $\MobiusBand$, then 
we have homotopy equivalences $\DiffId(\manif_1)\sim\manif_1$, implying $\pi_1\DiffId(\manif_1,1)=0$, whence $\DiffId(\manif_1,1)$ is contractible. But in this case $\manif$ belongs to the last column of Table~\ref{tbl:hom-types}, whence $\DiffId(\manif)$ is also contractible.

Similar arguments with the evaluation map by induction on $n$ show that $\DiffId(\manif_n,n)\sim\DiffId(\manif_{n-1}, n-1)$.
The details are left to the reader.
\end{proof}
\begin{rem}\label{rem:pi2D=0}
\rm Evidently, $\pi_2\Diff_{\id}(\manif,\crcnt)=0$
and $\pi_k\Diff_{\id}(\manif,\crcnt)\approx \pi_k\manif$ for $k\geq 3$.
Also $SO(3)=\RRR{}P^3$, whence $\pi_1 SO(3) = \ZZZ_2$,
$\pi_2 SO(3) = 0$, and $\pi_k SO(3) = \pi_k \sphere  = \pi_k \prjplane$ for $k\geq 3$ (cf. (2) of Theorem~\ref{th:hom-gr-orbits-c1}).
\end{rem}
\begin{lem}\label{lm:DiffIdMcr-is-contractible}
Let $\mrsfunc:\manif\to\Psp$ be a Morse mapping and $\crpt{i}$ $(i=0,1,2)$ the number of critical points of $\mrsfunc$ of index $i$.
If $\crpt{1}\geq1$, then $\DiffIdMcr$ is contractible.
\end{lem}
\begin{proof}
Suppose that $\manif$ is orientable.
Then $\chi\manif+\bndcnt=2-2\genus$ and from Morse equality
$\crpt{0}-\crpt{1}+\crpt{2} =\chi\manif$
we get:
$$
\crcnt+\bndcnt =
\crpt{0} + \crpt{1} + \crpt{2} + \bndcnt =
(\crpt{0} - \crpt{1} + \crpt{2})  + \bndcnt  + 2\crpt{1}  \geq
4-2\genus.
$$
A non-orientable surface $\manif$ is a connected sum of $\genus$ projective planes and $2$-sphere with $\bndcnt$ holes.
Moreover, $\chi\manif+\bndcnt=2-\genus$.
Hence similarly to the oriented case we obtain: $\crcnt+\bndcnt \geq 4-\genus$.
Then the contractibility of $\DiffIdMcr$ in both cases of $\manif$ follows from Table~\ref{tbl:hom-types}.
\end{proof}


\section{Two constructions related to a smooth mapping}\label{sect:prelim}
Let $\manif$ be a compact connected surface and $\mrsfunc:\manif\to\Psp$ a smooth mapping satisfying the conditions (i) and (ii) of Definition~\ref{defn:cond-i-ii-for-f}.

\subsection{Kronrod-Reeb graph of $\mrsfunc$}\label{sect:ReebGr}
Consider the partition of $\manif$ by the connected components of level-sets of $\mrsfunc$.
The corresponding factor-space (further denoted by $\Reebf$) has the structure of a one-dimensional CW-complex and is called the {\em Kronrod-Reeb (KR-) graph} of $\mrsfunc$.
Notice that the vertices of $\Reebf$ are of the following three types:
(a)~connected components of $\partial\manif$;
(b)~local extremes of $\mrsfunc$, i.e.$\!$ critical points of indices either $0$ or $2$;
(c)~critical components of level-sets of $\mrsfunc$.
We will call them \Dvert-, \Evert-, and \Cvert-vertices respectively.

Let $\prReeb:\manif\to\Reebf$ be the factor-map.
Then $\mrsfunc$ yields a unique (\KR-) mapping $\prReebf:\Reebf\to\Psp$ such that $\mrsfunc=\prReebf\circ\prReeb$.

Let $\edgeR$ be an edge of $\Reebf$, i.e. an open one-dimensional cell of $\Reebf$,
$\overline{\edgeR}$ the closure of $\edgeR$, and $\partial\edgeR=\overline{\edgeR}\setminus\edgeR$ the {\em boundary\/} of $\edgeR$.
Evidently, $\partial\edgeR$ consists of at most two points. 
Moreover, the case $|\partial\edgeR|\leq 1$ is possible only for $\Psp=\aCircle$.

An edge $\edgeR$ of $\Reebf$ will be called {\em external\/} if $\partial\edgeR$ contains a vertex of degree $1$, otherwise, $\edgeR$ is {\em internal.}

A homeomorphism $\theta:\Reebf\to\Reebf$ will be called an {\em automorphism\/} of the \KR-graph of $\mrsfunc$ provided $\prReebf=\prReebf\circ\theta$ and $\theta$ preserves each of the sets of \Dvert-, \Evert-, and \Cvert-vertices. 
Let $\AutfR$ be the group of all automorphisms of $\Reebf$. 
Evidently, each $\difM\in\Stabf$ yields a unique automorphism $\gdifM$ of $\Reebf$ so that the correspondence $\difM\mapsto\gdifM$ is a homomorphism
\chEq\begin{equation}\protect\label{equ:actSRall-hom}
\actSRall:\Stabf \to \AutfR
\end{equation}
which is not necessarily onto.

The following lemma is evident and can be easily deduced e.g. from A.~V.~Bolsinov and A.~T.~Fomenko~\cite{BolsFom} or E.~V.~Kulinich~\cite{Kulinich}. 
See also A.~Hatcher and W.~Thurston~\cite{HatcherThurston}.
\begin{lem}\label{lm:actSRall_Dpartf}
If $\mrsfunc:\manif\to\Psp$ is a simple Morse mapping, then $\actSRall$ is onto. 
\end{lem}

\subsection{Foliation of $\mrsfunc$}\label{subsect:foliation}
It is proved in~\cite{Prishlyak} that if $0\in\RRR^{2}=\CCC$ is an isolated non-extremal critical point of a smooth function $\gfunc:\RRR^{2}\to\RRR$, then there is a {\em homeomorphism\/} $\difM$ of $\CCC$ such that $\difM(0)=0$ and $\gfunc\circ \difM(z)={\rm Re}(z^k)$ for some $k\geq 1$.

It follows that the critical level-sets of a smooth map $\mrsfunc:\manif\to\Psp$ having only isolated critical points are embedded {\em graphs}, i.e. $1$-dimensional CW-complexes.

Hence $\mrsfunc$ yields on $\manif$ a one-dimensional foliation $\partitf$ with singularities: a subset $\omega \subset \manif$ is a leaf of this foliation iff $\omega$ is either a critical point of $\mrsfunc$ or a path-component of a set $\mrsfunc^{-1}(c)\setminus\singf$ for some $c\in\Psp$.

Let $\Dpartf$ be the group of diffeomorphisms of $\manif$ preserving every leaf of $\partitf$ and $\Dpartfplus$ the subgroup of $\Dpartf$ preserving orientations of these leaves.
Evidently, $\Dpartfplus$ and $\Dpartf$ are normal in $\Stabf$.

Let also $\DpartfId\subset\Dpartf$ and $\keractSRallId \subset \keractSRall$ be the subsets consisting of diffeomorphisms isotopic to $\id_{\manif}$ in $\Dpartf$, resp. in $\keractSRall$.
\begin{lem}\label{lm:Dpart-Stabf}
$\DpartfId = \keractSRallId = \StabIdf$.
\end{lem}
\begin{proof}
As $\Dpartf\subset \keractSRall\subset\Stabf$, we have
$\DpartIdf \subset \keractSRallId \subset \StabIdf.$
Thus it suffices to show that $\StabIdf \subset \DpartfId$.

Let $\cmap \in \StabIdf$.
Then there exists an isotopy $\Ahom:\manif\times I \to \manif$
such that $\Ahom_0=\id_{\manif}$, $\Ahom_1=\cmap$, and $\Ahom_t\in\Stabf$ so $\mrsfunc \circ \Ahom_t = \mrsfunc$, $t\in I$.

It follows that $\Ahom_t$ preserves $\singf$ and each level-set $\mrsfunc^{-1}(c)$.
Moreover, since $\singf$ is discrete and $\Ahom_0=\id_{\manif}$, it follows that $\Ahom_t(\pnt)=\pnt$ for each $\pnt\in\singf$ and that $\Ahom_t$ also preserves path-components of $\mrsfunc^{-1}(c)\setminus\singf$.
Hence $\Ahom_t\in\DpartfId$.
In particular, $\cmap=\Ahom_1\in\DpartfId$.
\end{proof}
\begin{lem}\label{lm:keep_orient_of_levels}
Let $\difM\in\Stabf$ be a diffeomorphism such that
$\difM(\omega)=\omega$ for some one-dimensional leaf $\omega$ of $\partitf$.

{\rm(1)}~Suppose that $\manif$ is oriented.
Then $\difM$ preserves the orientation of $\manif$ if and only if
it preserves the orientation of $\omega$.

{\rm(2)}~If $\difM\in\StabIdf$, then $\difM$ preserves the orientation of $\omega$.
\end{lem}
\begin{proof}
(1)
Not loosing generality we can assume that $\difM$ has a fixed point $\pnt\in\omega$.
Then $\difM$ preserves orientation of $\manif$ iff it preserves orientation of $T_{\pnt}\manif$.

Let $\vectOm$ be a non-zero tangent vector to $\omega$ at $\pnt$ and $\grad\mrsfunc(\pnt)$ the gradient vector of $\mrsfunc$ at $\pnt$
in some Riemannian metric on $\manif$.
Then the pair $(\grad\mrsfunc(\pnt), \vectOm)$ forms a basis of $T_{\pnt}\manif$.
Since $\mrsfunc\circ\difM=\mrsfunc$, we obtain that $T\difM(\grad\mrsfunc(\pnt))=\grad\mrsfunc(\pnt)$.
Moreover, from $\difM(\omega)=\omega$, we get $T\difM(\vectOm)=\alpha \vectOm$ for some $\alpha\not=0$.
Hence $\difM$ preserves orientation of $T_{\pnt}\manif$ iff $\alpha>0$, i.e. $\difM$ preserves the orientation of $\omega$.

(2)
Suppose that $\difM\in\StabIdf$.
If $\manif$ is oriented, then $\difM$ preserves the orientation of $\manif$,
and by (1) preserves the orientation of $\omega$.

Suppose that $\manif$ is non-orientable.
Let $\tmanif$ be an oriented double covering of $\manif$
and $\pr:\tmanif\to\manif$ the corresponding projection.
Then $\pr^{-1}(\omega)$ consists of two components $\tomega_1$ and $\tomega_2$.
Indeed, this is obvious for the case when $\omega$ is an open interval.
Otherwise, $\omega$ is a regular component of some level-set of $\mrsfunc$.
Then $\omega$ is {\em two-sided}, whence
$\pr^{-1}(\omega)$ consists of two components.

It follows that $\difM$ yields a diffeomorphism $\tdifM$ of $\tmanif$ which is isotopic to $\id_{\tmanif}$ and $\tdifM(\tomega_i)=\tomega_i$, $(i=1,2)$.
By (1) we obtain that $\tdifM$ preserves orientations of $\tomega_i$.
Hence, $\difM$ preserves the orientation of $\omega$.
\end{proof}


\section{Smooth shifts along trajectories of flows}\label{sect:sm_shifts}
First we briefly recall the results obtained in~\cite{Maks:Shifts}.
Let $\manif$ be a smooth compact $m$-dimensional manifold, $\Fld$ a vector field on $\manif$ tangent to $\partial\manif$, $\flow$ a flow generated by $\Fld$, and  $\Fix\flow$ the set of fixed points of $\flow$.

Define the following mapping $\Shift:\smr\to\smm$ by the formula:
$\Shift(\afunc)(x) = \flow(x,\afunc(x))$ for $\afunc\in\smr$ and $x\in\manif$.
We will call it a {\em shift-map along trajectories of $\flow$.}

Let $\difM:\manif\to\manif$ be a mapping and $X\subset\manif$.
We will say that a function $\afunc:X\to\RRR$ defined on a subset $X$ of $\manif$ is a {\em partial shift-function} for $\difM$ on $X$ provided $\difM(x)=\flow(x,\afunc(x))$ for all $x\in X$.
When $X=\manif$, such a function $\afunc$ is {\em global}.
In turn, $\difM$ will be called a {\em shift} along trajectories of $\flow$ by the function $\afunc$.

The set $\Zid=\Shift^{-1}(\id_{\manif})\subset\smr$ will be called the {\em kernel\/} of $\Shift$.
Thus for every $\alpha\in\Zid$ we have $\flow(x,\alpha(x))=x$.
\begin{lem}\label{lm:Shift_kernel}{\rm\cite[Corollary~6]{Maks:Shifts}} \ \
$\Shift(\afunc_1)=\Shift(\afunc_2)$ iff $\afunc_1 - \afunc_2 \in \Zid$.
Thus the image $\Shift(\smr) \subset C^{\infty}(\manif,\manif)$ is in the one-to-one correspondence with the factor group $\smr/\Zid$. 
\end{lem}

\begin{lem}\label{lem:Shift_kernel_class}{\rm\cite[Theorem~12]{Maks:Shifts}}
{\rm(a)}~Suppose that $\interior(\Fix\flow)\not=\varnothing$. Then $\Zid$ consists of functions vanishing on $\manif\setminus\interior(\Fix\flow)$.

Otherwise, we have another two possibilities: 

{\rm(b)}~$\Zid=\{0\}$. In this case $\Shift$ is injective. In particular, this is true when $\flow$ has at least one non-closed trajectory;

{\rm(c)}~$\Zid=\{n\perfunc\}_{n\in\ZZZ}\approx \ZZZ$, where $\perfunc>0$ is some smooth strictly positive function on $\manif$. In this case all non-fixed trajectories are closed, $\perfunc$ is constant on them, i.e. it is an integral of $\Fld$, and for every closed trajectory $\omega$ the value $\perfunc(\omega)$ is equal to the period of $\omega$.
\end{lem}
Denote by $\Cflow$ the subset of $\smmap$ consisting of smooth mappings
$\difM:\manif\to\manif$ such that 
(i) $\difM(\omega)\subset\omega$ for every trajectory $\omega$ of $\flow$ and  
(ii) for each $\pnt\in\Fix\flow$ the tangent linear map $T_{\pnt}\difM:T_{\pnt}\manif\to T_{\pnt}\manif$ at $\pnt$ is an isomorphism.

Let also $\Dflow=\Cflow\cap\DiffM$ and $\Dflowplus$ be the subset of $\Dflow$
consisting of diffeomorphisms preserving orientation of trajectories $\flow$.
Finally, let $\DflowId\subset\Dflow$ and $\CflowId\subset\Cflow$ be the subsets consisting of mappings that are homotopic in $\Dflow$, resp. in $\Cflow$, to $\id_{\manif}$.

Evidently, that the following set
$$ 
 \Dinvplus = \bigl\{ \afunc\in\smr \ | \ d\afunc(\Fld)(x) > -1 , \ x\in\manif \bigr\}.
$$ 
is open and convex in $\smr$.
\begin{lem}\label{lm:Shift_smr__CflowId}
{\rm\cite[Lemma~23]{Maks:Shifts}}
$$\Shift(\smr) \subset \CflowId \quad \text{and} \quad  \Shift(\Dinvplus)\subset\DflowId\subset\Dflowplus.$$
\end{lem}
\begin{defn}
\rm 
Let $\Fld$ be a vector field tangent to $\partial\manif$ and generating a flow $\flow$.
We say that $\Fld$ is \ \LLL\ {\em (locally linear)}\  if
for each $\pnt\in\Fix\flow$ there are local coordinates $(x_1,\ldots,x_m)$ in which $\pnt=0\in\RRR^{m}$ and $\Fld$ is a {\em linear\/} vector field, i.e. $\Fld(x)=Vx$, where $V$ is a constant $(m\times m)$-matrix.
A flow generated by \LLL\ vector field will also be called \LLL.
\end{defn}
\begin{thm}\label{th:Maks}{\rm\cite[Theorem~1]{Maks:Shifts}}~
If $\Fld$ is \LLL\ vector field, then 

{\rm(i)}~$\interior(\Fix\flow)=\varnothing$, thus the kernel $\Zid$ of shift-mapping $\Shift$ is of type either {\rm(b)} or {\rm(c)} of Lemma~\ref{lem:Shift_kernel_class}\,{\rm;}

{\rm(ii)}~$\Shift(\smr)=\CflowId$ \ and \ $\Shift(\Dinvplus) = \DflowId${\rm;}

{\rm(iii)}~$\Shift:\smr\to\CflowId$ and the restriction 
$\Shift:\Dinvplus\to\DflowId$ are (mutually) either homeomorphisms or a covering maps in the $C^{\infty}$ topologies of these spaces{\rm;}

{\rm(iv)}~the embedding $\DflowId\subset\CflowId$ is a homotopy equivalence.
Moreover, if $\Shift$ is injective $(\Zid=\{0\})$, then these spaces are contractible.
Otherwise, they are homotopy equivalent to $\aCircle$.
\end{thm}

We will also need the following partial variant of Theorem~\ref{th:Maks}(ii).

\begin{lem}\label{lm:local_th_Maks}
Let $\Fld$ be a vector field defined on a neighborhood $\nbh$ of $(0,0)\in\RRR^2$, $\flow$ a local flow generated by $\Fld$, and $\difM\in\Dflow$.
Then $\difM$ admits a smooth partial shift-function $\afunc$ defined on some neighborhood $V\subset\nbh$ of $(0,0)$, i.e. $\difM(x)=\flow(x,\afunc(x))$ for all $x\in V$ provided one of the following condition holds true:

{\rm1)} $\Fld(0,0)\not=0$, i.e. $(0,0)$ is regular,{\rm~\cite[Eq.~(10)]{Maks:Shifts}};

{\rm2)} $\Fld(x,y)=(y,x)$, i.e. $(0,0)$ is a ``saddle'' point,{\rm~\cite[Theorem~27]{Maks:Shifts}};

{\rm3)} $\Fld(x,y)=(-y,x)$, i.e. $(0,0)$ is a ``focus'',{\rm~\cite[Eq.~(25)]{Maks:Shifts}}.

In the cases {\rm1)} and {\rm2)} $\afunc$ is unique, and in the case {\rm3)} $\afunc$ is determined up to constant summand $2\pi k$, $k\in\ZZZ$.
\end{lem}

\subsection{Skew-symmetric flows}\label{subsect:SkewSymmetrFlows}
Let $\nmanif$ be a smooth non-orientable manifold, $\pr:\omanif\to\nmanif$ the oriented double covering of $\nmanif$, and $\invol:\omanif\to\omanif$ the smooth involution without fixed points generating the group $\ZZZ_2$ of covering slices of $\omanif$. 
Thus $\invol^2=\id_{\omanif}$ and $\pr\circ\invol=\pr$.

A continuous mapping $\difoM:\omanif\to\omanif$ will be called {\em symmetric} provided
$\difoM\circ\invol=\invol\circ\difoM$, and {\em skew-symmetric} if  $\difoM\circ\invol=\invol\circ\difoM^{-1}$.

Let $T\invol:T\manif\to T\manif$ be the corresponding tangent mapping.
We say that a vector field $\Fld$ on $\manif$ is {\em skew-symmetric} if $\Fld\circ\invol = -T\invol \circ\Fld$. 
This is equivalent to the condition that the flow $\flow$ generated by $\Fld$ is {\em skew-symmetric\/} in sense that $\flow_t\circ\invol=\invol\circ\flow_{-t}$ for all $t\in\RRR$.

Let $\flow$ be a skew-symmetric flow on $\manif$.
Then $\invol$ preserves foliation of $\flow$ but permutes its trajectories.
Therefore $\flow$ yields on $\nmanif$ some one-dimensional foliation $\nfoliat$ with singularities.
Let $\Cfol$ be the subset of $C^{\infty}(\nmanif,\nmanif)$
consisting of mappings $\difM:\nmanif\to\nmanif$
such that (i)~$\difM(\omega)\subset\omega$ for every leaf $\omega$ of $\nfoliat$ and (ii)~$\difM$ is a local diffeomorphism at each fixed point of $\flow$.
Put $\Dfol=\Cfol\cap\DiffnM$.

Let also $\CfolId\subset\Cfol$ and $\DfolId\subset\Dfol$ be the subsets consisting of mappings that are homotopic to $\id_{\nmanif}$ in $\Cfol$, resp. in $\Dfol$.

\begin{thm}\label{th:CDfolId-contr}
Let $\flow$ be a skew-symmetric flow on $\manif$.
If $\flow$ is \LLL, then $\CfolId$ and $\DfolId$ are contractible.
\end{thm}
First we establish two lemmas.
%
\begin{lem}\label{lm:lift_iso}
Let $\oCflowId\subset\CflowId$ and $\oDflowId\subset\DflowId$ be the subsets consisting of symmetric mappings.
Then we have the following homeomorphisms in the corresponding $C^{\infty}$-topologies:
$\CfolId \approx \oCflowId$ and $\DfolId \approx \oDflowId$.
\end{lem}
\begin{proof}
Every symmetric map $\difoM:\omanif\to\omanif$ projects to a unique map $\difnM:\nmanif\to\nmanif$, whence we have a natural projection $\isolift:\oCflowId\to\CfolId$ defined by $\isolift(\difoM)=\difnM$.

Let $\difnM\in\CfolId$, and $\difnM_t\in\CfolId$ be the homotopy between $\difnM_0=\id_{\nmanif}$ and 
$\difnM_1=\difnM$.
By the covering homotopy property there is a unique symmetric homotopy $\difoM_t\in\oCflowId$ such that  $\difoM_0=\id_{\omanif}$. Then $\isolift(\difoM_1)=\difnM$, whence $\isolift$ is onto.
Notice that another lifting of $\difnM$ is a map $\invol\circ\difoM_1$ having no invariant trajectories of $\flow$. Hence $\invol\circ\difoM_1\not\in\oCflowId$.
Therefore for each $\difnM\in\CfolId$ there is a unique lifting $\difoM\in\oCflowId$, whence $\isolift$ is bijective.
Since the projection $\pr:\omanif\to\nmanif$ is a smooth local diffeomorphism, it follows that $\isolift$ a homeomorphism in $C^{\infty}$-topologies.
It remains to note that $\isolift(\oDflowId)=\DfolId$.
\end{proof}
Thus for the proof of Theorem~\ref{th:CDfolId-contr}, we have to show the contractibility of $\oCflowId$ and $\oDflowId$.

Let $\Shift:C^{\infty}(\manif,\RRR) \to C^{\infty}(\manif,\manif)$ be the shift-mapping along trajectories of $\flow$, $\Zid=\Shift^{-1}(\id_{\manif})$ its kernel, and $\Dinvplus=\Shift^{-1}(\Dflowplus)$.
Denote 
$$\comptstmr{} = \{\ashift\in\stmr \ | \ \afunc\circ\invol +\afunc\in\Zid\},$$
and for every $\mu\in\Zid$ let 
$$\comptstmr{\mu}=\{\ashift\in\stmr \ | \ \afunc\circ\invol +\afunc = \mu\}$$ and $\Dinv{\mu}=\comptstmr{\mu}\cap\Dinvplus.$
Evidently, $\comptstmr{\mu}$ and $\Dinv{\mu}$ are convex.

Since $\flow$ is \LLL, we have by Theorem~\ref{th:Maks}(i) that either $\Zid=\{0\}$, then $\comptstmr{}=\comptstmr{0}$, or $\Zid=\{n\perfunc\}_{n\in\ZZZ}$, $(\perfunc>0)$.
In the latter case we will denote $\comptstmr{n\perfunc}$ and $\Dinv{n\perfunc}$ simply by $\comptstmr{n}$ and $\Dinv{n}$ respectively.
Then each $\comptstmr{n}$ is a connected component of $\comptstmr{}$.
\begin{lem}\label{lm:prop-comptstmr}
{\rm (1)} $\comptstmr{0}\cap\Zid=\{0\}$, whence the restriction of $\Shift$ to $\comptstmr{0}$ is injective.

{\rm (2)}
Suppose that $\Zid=\{n\perfunc\}_{n\in\ZZZ}$, $\perfunc>0$.
Then $\perfunc=\perfunc\circ\invol$.
Moreover, $\comptstmr{n} = \comptstmr{0}+\frac{n\perfunc}{2}$,
whence $\Shift(\comptstmr{a})=\Shift(\comptstmr{b})$ iff \ $a\equiv b\,\mathrm{mod}\,2$.
Thus $\Shift(\comptstmr{})$ consists of two connected components $\Shift(\comptstmr{0})$ and $\Shift(\comptstmr{1})$.

{\rm (3)} Let $\ashift\in\smr$. 
Then $\Shift(\ashift)$ is symmetric iff $\afunc\in\comptstmr{}$. 
This means that $\comptstmr{}=\Shift^{-1}(\oCflow)$. 
Moreover, $\comptstmr{}\cap\Dinvplus=\Shift^{-1}(\oDflow)$.
In particular we have $\Shift(\comptstmr{})\subset\oCflow$, whence by {\rm(2)\/} $\Shift(\comptstmr{0})= \Shift(\comptstmr{})_{\id}\subset\oCflowId$.
Similarly, $\Shift(\Dinv{0})\subset\oDflowId$.
\end{lem}
\begin{proof}
(1)
If $\Zid=\{0\}$, then $\comptstmr{0}\cap\Zid=\{0\}$ holds trivially.
Suppose that $\Zid=\{n\perfunc\}_{n\in\ZZZ}$ and $m\perfunc\in\comptstmr{0}$ for some $m\in\ZZZ$. 
Then $m\perfunc\circ\invol= -m\perfunc$. 
Since $\perfunc>0$, we get $m=0$, whence $\comptstmr{0}\cap\Zid=\{0\}$.

By Lemma~\ref{lm:Shift_kernel} $\Shift$ yields a bijection between $\somr/\Zid$ and the image $\IM\Shift$, whence the restriction of $\Shift$ to $\comptstmr{0}$ is injective.

(2)
First we prove that $\perfunc=\perfunc\circ\invol$.
Notice that for every point $x\in\omanif\setminus\Fix\flow$ the value $\perfunc(x)$ is equal to the period of $x$ with respect to $\flow$.
Therefore $\perfunc$ is constant along trajectories of $\flow$.
Moreover, since $\flow$ is skew-symmetric, it follows that for every trajectory $\omega$ of $\flow$ its period coincides with the period of the trajectory $\invol(\omega)$, whence $\perfunc=\perfunc\circ\invol$.

Suppose that $\ashift\in\comptstmr{n}$, i.e. $\ashift+\ashift\circ\invol = n\perfunc$.
It suffices to show that $\ashift_1=\ashift+\perfunc/2 \in\comptstmr{n+1}$.
Indeed,
$$
\ashift_1+\ashift_1\circ\invol =
\ashift+\perfunc/2 + \ashift\circ\invol+
(\perfunc\circ\invol)/2 =
n\perfunc + \perfunc = (n+1)\perfunc.
$$
Since $\Shift(\ashift+\perfunc) =  \Shift(\ashift)$ for all $\ashift\in\stmr$, it follows that
$$
  \Shift(\comptstmr{n}) = \Shift(\comptstmr{n}+2\perfunc/2)= \Shift(\comptstmr{n+2}).
$$

(3)
Suppose that $\difoM=\Shift(\ashift)$, where $\ashift\in\smfunc$.
Then
$$
\invol\circ\difoM\circ\invol(x) = \invol\circ\flow(\invol(x),\ashift\circ\invol(x))=\flow(x,-\ashift\circ\invol(x)).
$$
Hence this mapping coincides with $\difoM$ iff $\ashift+\ashift\circ\invol\in\Zid$.
In other words, $\difoM\in\oCflow$ iff $\ashift\in\comptstmr{}$.
\end{proof}

\begin{proof}{Proof of Theorem~\ref{th:CDfolId-contr}} 
Suppose that $\flow$ is \LLL. 
Then by Theorem~\ref{th:Maks} $\Shift(\smr)=\CflowId$, $\Shift(\Dinvplus)=\DflowId$, and $\Shift$ is either a homeomorphism or a covering map.
Hence from (3) of Lemma~\ref{lm:prop-comptstmr} it follows that $\Shift(\comptstmr{0})=\oCflowId$ and $\Shift(\Dinv{0})=\oDflowId$.

From (1) of Lemma~\ref{lm:prop-comptstmr} we get that $\Shift$ homeomorphically maps convex sets $\comptstmr{0}$ and $\Dinv{0}$ onto $\oCflowId$ and $\oDflowId$ respectively.
Hence $\oCflowId$ and $\oDflowId$ are contractible, and by Lemma~\ref{lm:lift_iso} so are $\CfolId$ and $\DfolId$.
\end{proof}

\subsection{Flows on $\aCircle\times[0,1]$}
Let $\flow$ be a flow on the cylinder $\Cyl = \aCircle\times\Interv$ such that the trajectories of $\flow$ are of the form $\aCircle\times\{t\}$, $t\in\Interv$.
Thus $\Fix\flow=\varnothing$, whence $\flow$ is \LLL.
For every $\eps\in(0,\frac{1}{2})$ denote $\Seps=\aCircle\times[0,\eps] \cup [1-\eps,1].$ 

The following lemma is easy and we left its proof to the reader.
\begin{lem}\label{lm:ext-shift-func}
{\rm (1)}
$\Dflowplus=\DflowId$. 
Hence for every $\difM\in\Dflowplus$ there is a shift-function with respect to $\flow$.

{\rm (2)}
Suppose that a diffeomorphism $\difM\in\Dflowplus$ has a partial shift-function $\ashift:\aCircle\times[0,\eps]\to\RRR$, where $\eps\in(0,1)$.
Then $\ashift$ uniquely extends to a global shift-function for $\difM$.

{\rm (3)}
Let $\Dtw$ be a Dehn twist about the curve $\aCircle\times\{1/2\}$.
Suppose that $\difM$ is the identity on $\Seps$.
Then $\difM$ is isotopic in $\DflowId$ to some degree $\Dtw^{k}$ for $k\in\ZZZ$.
Moreover, if in addition, $\ashift=0$ on $\aCircle\times[0,\eps]$, then $\ashift=0$ also on $\aCircle\times[1-\eps,1]$ if and only if $k=0$, i.e. $\difM$ is isotopic to the identity with respect to $\Seps$.

{\rm (4)}
Let $\difM:\Stwoeps\to\Stwoeps$ be a diffeomorphism that preserves sets of the form $\aCircle\times\{t\}$, changes their orientation, and such that $\difM^2=\id_{\Stwoeps}$.
Then there exists a diffeomorphism $\adifM$ of $\Cyl$ such that $\adifM=\difM$ on $\Seps$, $\adifM\in\Dflow$, and $\adifM^2=\id_{\Cyl}$.
\end{lem}

\subsection{Extension of a shift-function}
Let $\flow$ be a flow on a manifold $\manif$, $\restM$ closed and invariant under $\flow$ subset of $\manif$, and $\nbh$ an invariant neighborhood of $\restM$.

Suppose that $\difM\in\Dflowplus$ has a shift-function $\ashift$ in $\nbh$, i.e. $\difM(x)=\flow(x,\ashift(x))$ for all $x\in\nbh$.
Then, in general, $\ashift$ does not extends to a global shift-function for $\difM$.
A possible obstruction is that $\difM$ can be non-isotopic to $\id_{\manif}$.
Nevertheless, the existence of a partial shift-function allows us to simplify $\difM$, provided $\flow$ admits an integral on $\nbh$.

\begin{lem}\label{lm:sympl-h}
Suppose that there exist an invariant neighborhood $\anbh$ of $\restM$ and a smooth function $\mu:\manif\to\Interv$ such that $\overline{\anbh}\subset\nbh$ and 
{\rm (i)}~$\mu=1$ on $\anbh$; 
{\rm (ii)}~$\mu=0$ in a neighborhood of the set $\overline{\manif\setminus\nbh}$; 
{\rm (iii)}~$\mu$ is constant on every trajectory of $\flow$.
Then $\difM$ is isotopic in $\Dflowplus$ to a diffeomorphism that is the identity on $\overline{\anbh}$ and coincides with $\difM$ on $\manif\setminus\nbh$.
\end{lem}
\begin{proof}
Consider the following function $\nu = \ashift \cdot \mu$.
It follows from (ii) that $\nu$ is well defined on all of $\manif$.

\begin{claim}
For every $t\in\Interv$ the mapping $\bdif_t:\manif\to\manif$ defined by
$\bdif_t(x) = \flow(x,t\,\nu(x))$ is a diffeomorphism of $\manif$.
\end{claim}
\proof 
By (ii) of Theorem~\ref{th:Maks} $\Shift(\Dinvplus) = \DflowId$, i.e. $\bdif_t=\Shift(t\nu)$ is a diffeomorphism iff $t\nu\in\Dinvplus$, i.e.
\chEq\begin{equation}\protect\label{equ:dmu-g-m1}
d(t \, \nu(\Fld))(x) = t \, d\nu(\Fld)(x)>-1
\end{equation}
for all $x\in\manif$, where $\Fld$ is a vector field generating $\flow$, see also~\cite[Theorem~19]{Maks:Shifts}.

Since $\difM$ is a diffeomorphism, we have that $d\ashift(\Fld)>-1$.
Moreover, the condition (iii) of this lemma implies $d\mu(\Fld)=0$.
Therefore,
$$
t \, d\nu(\Fld) = t \, d(\mu\,\ashift)(\Fld) =
t \, \ashift \, d\mu(\Fld) + t \, \mu \, d\ashift(\Fld) =
t \, \mu \, d\ashift(\Fld) > -1.
\text{\qed}
$$

Now the following isotopy $\difM_t = \difM \circ \bdif_t$ deforms $\difM$ in $\Dflowplus$ to a diffeomorphism $\difM_1$ that is the identity on $\overline{\anbh}$ and coincides with $\difM$ on $\manif\setminus\nbh$.
Lemma~\ref{lm:sympl-h} is proved.
\end{proof}


\section{Proof of Theorem~\ref{th:StabIdf-hom-type}}
\label{sect:proof-th:StabIdf-hom-type-or}
First suppose that $\manif$ is oriented.
The idea is to identify $\StabIdf$ with $\DflowId$ for some flow on $\manif$ and then apply Theorem~\ref{th:Maks}.

\begin{lem}\label{lm:exist-Fld}
Let $\mrsfunc:\manif\to\Psp$ be a Morse mapping on a smooth orientable surface $\manif$.

{\rm (1)}
There exists an \LLL\ vector field $\Fld$ on $\manif$ whose trajectories are precisely the leaves of $\partitf$.

{\rm (2)}
Suppose that there is a smooth changing orientation involution $\invol$ without fixed points such that $\mrsfunc\circ \invol=\mrsfunc$.
Then $\Fld$ can be chosen skew-symmetric with respect to $\invol$.
\end{lem}
\begin{proof}
(1) Let $\omega$ be a $2$-form determining some symplectic structure on $\manif$. 
Then $\omega$ yields an isomorphism $\psi:T^{*}\manif\to T\manif$ between the cotangent and tangent bundles of $\manif$.
Consider the {\em skew-gradient\/} vector field (also called Hamiltonian vector field for $\mrsfunc$) $\sgradf = \psi (d\mrsfunc)$.
Evidently, the trajectories of $\sgradf$ are precisely the leaves of $\partitf$.
Moreover, since $\mrsfunc$ is constant on the connected components of $\partial\manif$, we obtain that $\sgradf$ is tangent to $\partial\manif$.

Let $\pnt$ be a critical point of $\mrsfunc$ and $(x,y)$ local coordinates at $\pnt$ such that $\pnt=(0,0)$ and $\mrsfunc(x,y)= \eps_{x} x^2 + \eps_y y^2$, where $\eps_x,\eps_y=\pm1$.
Then we define following vector field near $\pnt$ by: 
$\Fld_{\pnt}(x,y)=(\eps_y y,-\eps_x x)$.
Evidently, $\Fld_{\pnt}$ is linear and its trajectories are subsets of leaves of $\partitf$.
Moreover, $\Fld_{\pnt}$ is collinear with $\sgradf$.

Now using the partition of unity technique we can glue $\sgradf$ with all $\Fld_{\pnt}$ $(\pnt\in\singf)$ so that the resulting vector field $\Fld$ will be \LLL.

(2) Let $\invol\in\Stabf$ be a smooth changing orientation involution without fixed points
and $\Fld$ is \LLL\ vector field constructed as just above.
We can make $\Fld$ skew-symmetric by replacing it by the vector field
$\frac{1}{2}\bigl( \Fld - T\invol\circ\Fld\circ\invol \bigr)$, but this vector field may loose \LLL\ property.
To solve this problem, we should improve the construction of $\Fld$.

First divide the set of singular points $\singf$ into two disjoint sets $\singf^1$ and $\singf^2$ so that $\invol(\singf^1)=\singf^2$.
Then for $\pnt\in\singf^1$ we define $\Fld_{\pnt}$ as in (1) and for $\pnt\in\singf^2$ by the formula $\Fld_{\pnt} = -T\invol\circ\Fld_{\invol(\pnt)}\circ\invol$.

Gluing $\sgradf$ with $\Fld_{\pnt}$ $(\pnt\in\singf)$ as above we get a vector field $\Fld$ which is \LLL\ and skew-symmetric with respect to $\invol$ near its singular points.
Then the vector field $\frac{1}{2}\bigl( \Fld - T\invol\circ\Fld\circ\invol \bigr)$ is \LLL\ and skew-symmetric.
\end{proof}

Let $\flow$ be a flow generated by the vector field $\Fld$ of this lemma.
Then $\Dflow = \Dpartf$ and
by Lemma~\ref{lm:Dpart-Stabf} we get $\DflowId = \DpartfId = \StabIdf.$

Since $\flow$ is \LLL, it follows from (iv) of Theorem~\ref{th:Maks} that $\StabIdf$ is either contractible or has the homotopy type of $\aCircle$.
Suppose that $\mrsfunc$ has at least one critical point of index $1$.
Then from the structure of level-sets of $\mrsfunc$ near such a point we obtain that $\flow$ has non-closed trajectories.
Hence by (iv) of Theorem~\ref{th:Maks} $\StabIdf$ is contractible.

If $\mrsfunc$ has no critical points of index $1$, then $\mrsfunc$ belongs to one of the  types \typeA-\typeDk\ of Table~\ref{tbl:hom_groups_of_orbits}.
In each of these cases it is easy to construct a flow on $\manif$ satisfying the statement of Lemma~\ref{lm:exist-Fld}, and such that the kernel of the corresponding shift-mapping coincides with the set constant functions taking integral values.
Hence $\StabIdf$ is homotopy equivalent to $\aCircle$.
Theorem~\ref{th:StabIdf-hom-type} is proved for oriented surfaces.

Let $\nmanif$ be a non-orientable surface and $\nmrsfunc:\nmanif\to\Psp$ a Morse mapping.
We have to prove that $\StabIdnf$ is contractible.
Let $\pr:\omanif\to\nmanif$ be an oriented double covering, $\omrsfunc=\nmrsfunc\circ\pr:\omanif\to\Psp$, $\invol$ the corresponding involution on $\omanif$, and $\partitof$ and $\partitnf$ the corresponding foliations on $\omanif$ and $\nmanif$ defined by $\omrsfunc$ and $\nmrsfunc$ respectively.

Let $\Fld$ be a skew-symmetric \LLL\ vector field on $\manif$ constructed in (2) of Lemma~\ref{lm:exist-Fld} for $\omrsfunc$.
Then $\partitnf$ is induced by $\Fld$, whence by Theorem~\ref{th:CDfolId-contr} $\DpartIdnf$ is contractible and by Lemma~\ref{lm:Dpart-Stabf} coincides with $\StabIdnf$.
\qed


\section{The group $\pi_0\Dpartfplus$}\label{sect:pi-0Dpartf}
In the sequel, $\mrsfunc:\manif\to\Psp$ will be a Morse mapping, $\Reebf$ the Kronrod-Reeb graph of $\mrsfunc$, and $\prReeb:\manif\to\Reebf$ the factor-map,

Let $\edgeR$ be an edge of $\Reebf$. 
Then $\prReeb^{-1}(\edgeR)$ contains no critical points of $\mrsfunc$, whence there exist a diffeomorphism $\psi:\aCircle\times(0,1) \to \prReeb^{-1}(\edgeR)$ and a homeomorphism $\delta:(0,1)\to\edgeR$ such that the following diagram is commutative:
$$
\begin{CD}
\aCircle\times(0,1) @>{\psi}>> \prReeb^{-1}(\edgeR) \\
@V{p_2}VV @VV{\prReeb}V \\
(0,1) @>{\delta}>> \edgeR @>{\prReebf}>> \Psp.
\end{CD}
$$

Here $p_2$ is a standard projection.
Regarding $\aCircle$ as the unit circle in the complex plane, choose a smooth function $\mu:(0,1)\to[0,1]$ such that $\mu(0,1/4]=0$ and $\mu[3/4,1)=1$, and define the following diffeomorphism $\Dtw$ of $\aCircle\times(0,1)$ by the formula: 
$\Dtw(z,t) = \bigl(e^{2\pi i \mu(t)} z, t\bigr).$
Thus $\Dtw$ is a Dehn twist along $\aCircle\times\{1/2\}$.
Then the following diffeomorphism
$$\Dtw_{\edgeR}(x) = \left\{
\begin{array}{cl}
x, & x \in \manif\setminus \prReeb^{-1}(\edgeR), \\
\psi \circ \Dtw \circ \psi^{-1}, & x \in \prReeb^{-1}(\edgeR)
\end{array}
\right.
$$
is a Dehn twist along the curve $\crv_{\edgeR} = \psi(\aCircle\times\{1/2\})$.
Notice that $\Dtw_{\edgeR}$ is the identity on a neighborhood of $\manif\setminus\prReeb^{-1}(\edgeR)$ and $\Dtw_{\edgeR}$ reserves each leaf of the foliation $\partitf$, i.e.$\!$ $\Dtw_{\edgeR} \in \Dpartfplus \subset \Stabf$.
We will call $\Dtw_{\edgeR}$ a Dehn twist about the edge $\edgeR$.

The Dehn twist $\Dtw_{\edgeR}$ and the corresponding curve $\crv_{\edgeR}$
will be called either {\em external\/} or {\em internal} with respect to $\edgeR$.
Evidently, $\edgeR$ is an external edge if and only if
$\Dtw_{\edgeR}$ is isotopic to $\id_{\manif}$ relative to $\singf$.
Denote by $\intecnt$ the number of internal edges of $\mrsfunc$.
\begin{rem}
\rm Dehn twists along regular components of level-sets of Morse functions play a crucial role in~\cite{Maks:PathComp} for the description of connected components of the space of Morse mappings on surfaces.
\end{rem}
\begin{thm}\label{th:twgr-descr}
The group $\pi_0\Dpartfplus$ is generated by the internal Dehn twists.
If $\mrsfunc$ is of type \typeEk\ of Table~\ref{tbl:hom_groups_of_orbits}, then $\pi_0\Dpartfplus \approx \ZZZ_{2}$.
Otherwise, internal Dehn twists form a basis of $\pi_0\Dpartfplus$, and thus 
$\pi_0\Dpartfplus \approx \ZZZ^{\intecnt}$, where $l$ is the number of internal Dehn twists.
\end{thm}
\proofstyle{Proof. Oriented case.}
Suppose that $\manif$ is oriented.
Let $\edgeR_{1}, \ldots,\edgeR_{\intecnt}$ be internal edges of $\Reebf$ and $\twgr$ the subgroup of $\pi_0\Dpartf$ generated by the Dehn twists about them.
Since the internal edges are disjoint, we see that $\twgr$ is abelian.
For simplicity denote $\crv_{\edgeR_{i}}$ and $\Dtw_{\edgeR_{i}}$ by $\crv_{i}$ and
$\Dtw_{i}$ respectively. 
Let also $\totsupp=\cup_i\supp\Dtw_i$ be the union of supports of $\Dtw_i$.

Let $\flow$ be a flow on $\manif$ constructed in (1) of Lemma~\ref{lm:exist-Fld}.
Then we can identify $\Dpartfplus$ with $\Dflowplus$.
Now the proof for oriented case consists of the following two lemmas.

\begin{lem}\label{lm:twgr_Zintecnt}
The isotopy classes of internal Dehn twists are independent in $\pi_0\DiffMcr$.
Therefore they are independent in $\pi_0\Dpartfplus$, whence $\twgr\approx\ZZZ^{\intecnt}$.
\end{lem}
\begin{proof}
Notice that the group $\twgr$ acts on $H_1(\manif\setminus\singf,\ZZZ)$ by the formula:
$\Dtw_i\cdot x = x + \langle x, \crv_i \rangle \crv_i,$ for $x\in  H_1(\manif\setminus\singf,\ZZZ)$, where $\langle \cdot, \cdot\rangle$ is the intersection form.
Hence,
\begin{equation}\protect\label{equ:prod-Dehn-twist}
\prod_{i=1}^{\intecnt} \Dtw^{m_i}_i \cdot x =
x + \sum_{i=1}^{\intecnt} m_i \langle x, \crv_i \rangle \crv_i.
\end{equation}
Since the internal curves represent linearly independent $1$-cycles in $H_1(\manif\setminus\singf,\ZZZ)$, it follows from~\eqref{equ:prod-Dehn-twist} that $\Dtw_i$ are independent in $\twgr$.
\end{proof}
\begin{lem}\label{lm:ex-shift-for-difM}
For every $\difM\in\Dflowplus$ there exists a partial shift-function defined on $\manif\setminus\totsupp$.
\end{lem}
It follows from Lemma~\ref{lm:ex-shift-for-difM} that each $\difM\in\Dflowplus$ is isotopic in $\Dflowplus$ to a diffeomorphism that is the identity on $\manif\setminus\totsupp$ and therefore belongs to $\twgr$.
Then from Lemma~\ref{lm:twgr_Zintecnt} we get $\pi_0\Dflowplus=\twgr\approx\ZZZ^{\intecnt}$ which proves our theorem.

\proofstyle{Proof of Lemma~\ref{lm:ex-shift-for-difM}.}
First we define a certain subset $\excrlev$ of $\manif$ and construct a shift-function $\ashift$ for $\difM$ near $\excrlev$.
Consider two cases.

\chFig\begin{figure}[ht]
\begin{tabular}{ccccc}
\includegraphics[height=3.5cm]{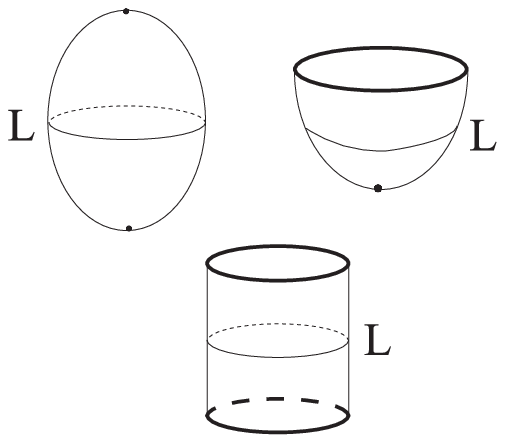}
& \qquad &
\includegraphics[height=3.5cm]{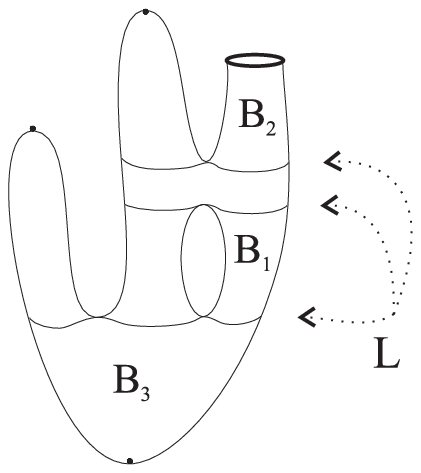}
& \qquad &
\includegraphics[height=3.5cm]{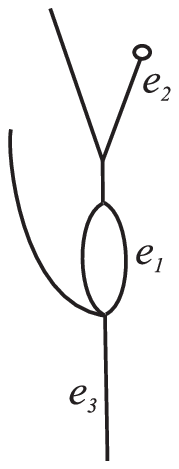} \\
a) & \qquad & b) & \qquad & c)
\end{tabular}
\caption{}
\protect\label{fig:constr-L}
\end{figure}

1) Suppose that $\crpt{1}=0$.
Then $\flow$ has no non-closed trajectories and the kernel of shift-mapping $\Shift$ along trajectories of $\flow$ is $\Zid=\{n\perfunc\}_{n\in\ZZZ}$.
Let $\excrlev$ be an arbitrary regular component of a level-set of $\mrsfunc$, see Figure~\ref{fig:constr-L}~a).
Then $\excrlev$ is a closed trajectory of $\flow$ and
a shift-function $\ashift$ for $\difM$ can be defined in a neighborhood of $\excrlev$
up to the summand $n\perfunc$.

2) If $\crpt{1}>0$, then we let $\excrlev$ to be the union of those critical components of level-sets of $\mrsfunc$ that contain critical points of index $1$, see Figure~\ref{fig:constr-L}~b).
Let $\pnt$ be a critical point of $\mrsfunc$ of index $1$.
Since in some local coordinates at $\pnt=0$, $\mrsfunc(x,y)=x^2-y^2$, we may choose $\Fld$ so that $\Fld(x,y)=(y,x)$. 
Then by 2) of Lemma~\ref{lm:local_th_Maks}, in some neighborhood $\nbh$ of $\pnt$ there exists a unique smooth function $\ashift$ such that $\difM(x)=\flow(x,\ashift(x))$ for all $x\in\nbh$.

Let now $\omega$ be a non-closed trajectory of $\flow$.
Since $\difM(\omega)=\omega$, it follows that for every $x\in\omega$ there exists a unique number $\ashift(x)$ such that $\difM(x)=\flow(x,\ashift(x))$.
Moreover, by 1) of Lemma~\ref{lm:local_th_Maks} $\ashift$ uniquely and smoothly extends onto some neighborhood of $\omega$. 

Thus $\difM$ admits a shift-function near $\excrlev$.
It remains to extend $\ashift$ to $\manif\setminus\totsupp$.
Let $\Cyl$ be a connected component of $\manif\setminus\excrlev$, it corresponds to some edge $\edgeR$ of $\Reebf$.
Consider three cases, see Figure~\ref{fig:constr-L}~b) and c).

{\em Case 1:}~$\edgeR$ is an internal edge of $\Reebf$.
Then $\Cyl$ includes $\supp\Dtw_i$ for some $i=1,\ldots,\intecnt$ and we should extend $\ashift$ to $\Cyl\setminus\supp\Dtw_i$.
Notice that $\Cyl\setminus\supp\Dtw_i$ is a union of two open cylinders and $\ashift$ is defined at ``one side'' of each of these cylinders.
Then by (2) of Lemma~\ref{lm:ext-shift-func} $\ashift$ extends to a shift function for $\difM$ on  $\Cyl\setminus\supp\Dtw_i$.

{\em Case 2:}~$\edgeR$ is external and has a \Dvert-vertex, i.e. $\Cyl\cap\partial\manif\not=\varnothing$.
Then $\Cyl\cap\totsupp=\varnothing$ and $\ashift$ extends to $\Cyl$ by (2) of Lemma~\ref{lm:ext-shift-func}.

{\em Case 3:}~$\edgeR$ is external and has an \Evert-vertex.
Then $\Cyl$ is an open $2$-disk containing a unique critical point $\pnt$ of $\mrsfunc$.
Again by (2) of Lemma~\ref{lm:ext-shift-func} and $\ashift$ extends to a shift-function for $\difM$ on $\Cyl\setminus\pnt$.
Moreover, since $\mrsfunc(x,y)=\pm(x^2+y^2)$ in some local coordinates $(x,y)$ at $\pnt=0$, we can choose $\Fld(x,y)=(-y,x)$, whence by 3) of Lemma~\ref{lm:local_th_Maks} $\ashift$ can be defined smoothly at $\pnt$.
\qed

\proofstyle{Proof. Non-orientable case.}
Let $\nmanif$ be a non-orientable surface, $\pr:\omanif\to\nmanif$ its oriented double covering, $\nmrsfunc:\nmanif\to\Psp$ a Morse mapping, and $\omrsfunc=\nmrsfunc\circ\pr$.
We will show that $\pi_0\Dpartnfplus$ is generated by internal Dehn twists.

1) Suppose that $\nmrsfunc$ is not of type \typeEk.
Let $\flow$ be a skew-symmetric \LLL\ flow on $\omanif$ generated by the vector field $\Fld$ of (2) Lemma~\ref{lm:exist-Fld}.
Then $\Dpartnfplus$ is naturally identified with $\oDpartofplus = \oDflowplus$.

Notice that the internal curves $\crv_i$ ($i=1,\ldots,\intecnt$) are two-sided, whence
$\pr^{-1}(\crv_i) \subset\omanif$ consists of two connected components $\tcrv_{i1}$ and $\tcrv_{i2}$ which are internal with respect to $\omrsfunc$.
Then there are two internal Dehn twists $\tDtw_{i1}$ and $\tDtw_{i2}$ about $\tcrv_{i1}$ and $\tcrv_{i2}$ (resp.) such that $\tDtw_{i2}= \invol\circ\tDtw_{i1}\circ\invol$ and $\tDtw_{i2}\circ\tDtw_{i1}$ is a lifting of $\Dtw_i$ belonging to $\oDflowplus$. 

By the oriented case of this theorem $\Dtw_{ij}$ ($i=1,\ldots,\intecnt$, $j=1,2$) form a basis of $\pi_0\Dflowplus$.
Hence $\tDtw_{i2} \circ\tDtw_{i1}$ are independent in $\pi_0\oDflowplus$ and generate the subgroup isomorphic with $\ZZZ^{\intecnt}$. 
Then it remains to prove the following statement.
\begin{claim}
The isotopy classes of $\tDtw_{i2} \circ\tDtw_{i1}$ generate $\pi_0\oDflowplus$.
\end{claim}
\begin{proof}
Denote by $\totsupp$ the union of supports of $\Dtw_{ij}$.
Suppose that $\nmrsfunc$ (and therefore $\omrsfunc$) has at least one critical point of index $1$.
Then the kernel $\Zid$ of the shift-mapping $\Shift$ along trajectories of $\flow$ is trivial: $\Zid=\{0\}$.

Let $\difoM\in\oDflowplus$.
Then by Lemma~\ref{lm:ex-shift-for-difM} we have $\difoM(x)=\flow(x,\ashift(x))$ for some smooth function $\ashift$.
Since $\difoM$ is symmetric, it follows from Lemma~\ref{lm:prop-comptstmr} that $\ashift+\ashift\circ\invol\in\Zid =\{0\}$, whence $\ashift\circ\invol = -\ashift$ on $\manif\setminus\totsupp$.

Let $\mu$ be a function of the statement of Lemma~\ref{lm:sympl-h}.
Then the function $\mu_1=(\mu+\mu\circ\invol)/2$ is symmetric and also satisfies the statement of Lemma~\ref{lm:sympl-h}.
Therefore $\mu_1 \cdot \ashift$ is a smooth skew-symmetric function on all of $\manif$.
Whence $\difoM_t(x)=\difoM\circ\flow(x,-t\cdot \mu_1(x) \cdot \ashift(x))$ is a symmetric isotopy of $\difoM$ in $\oDflowplus$ to a symmetric diffeomorphism that is the identity on $\manif\setminus\totsupp$.
So we can assume that $\difoM$ is a product of some $\Dtw_{ij}$.
Since $\difoM$ is symmetric, we obtain that $\difoM$ is in fact a product of some $\tDtw_{i2}\circ\tDtw_{i1}$.
\end{proof}

2) Suppose that $\nmrsfunc$ is of type \typeEk.
Then $\nmanif=\Kleinb$ is a Klein bottle, $\nmrsfunc:\Kleinb\to\aCircle$, $\omanif=\torus$,
the Kronrod-Reeb graph of $\nmrsfunc$ is a circle, and there is a unique internal Dehn twist $\Dtw_1$.

Let $\Dtw_{11}$ and $\Dtw_{12}$ be the liftings of $\Dtw_1$.
Then it is easy to see that $\Dtw_{11}\circ\Dtw_{12}$ is isotopic to $\id_{\manif}$ in $\Dflow$, i.e. $\Dtw_{11}\circ\Dtw_{12}\in\DflowId$.

\begin{claim}\label{clm:oDpartofplus=oDflow-cap-DflowId}
$\oDpartofplus = \oDflow\cap\DflowId$, whence $\pi_0\oDpartofplus=\ZZZ_2$
by Lemma~\ref{lm:prop-comptstmr}.
\end{claim}
\begin{proof}
Evidently, $\oDflow\cap\DflowId \subset\oDpartofplus$.
Conversely, Let $\difoM\in\oDflowplus \subset \oDflow$.
We have to show that $\difoM\in\DflowId$.

Notice that $\omanif\setminus\totsupp$ consists of two connected
components $\Cyl_1$ and $\Cyl_2$ that are diffeomoprhic to $\aCircle\times(0,1)$,
and such that $\invol(\Cyl_1)=\Cyl_2$.
Then by Lemma~\ref{lm:ext-shift-func}, there is a partial
shift-function $\ashift$ for $\difM$ on $\Cyl_1$.
We extend $\ashift$ on $\Cyl_2$ by
$\ashift(x) = -\ashift\circ\invol(x)$, $x\in\Cyl_2$.
Then $\ashift$ is skew-symmetric
and by the arguments used in the previous case 1) $\difM$ is isotopic
in $\oDflowplus$ to some degree of $\Dtw_{12} \circ\Dtw_{11} \in \DflowId$.

Claim~\ref{clm:oDpartofplus=oDflow-cap-DflowId} and Theorem~\ref{th:twgr-descr} are proved.
\end{proof}


\section{Ends at a subgraph}\label{sect:ends}
Let $\manif$ be a compact surface.
By a subgraph $\crcomp \subset \interior\manif$ we mean a one-dimensional CW-subcomplex of some cellular division of $\manif$.

Let $\crcomp \subset \interior\manif$ be a finite connected subgraph, $\nbh_1,\nbh_2$ two regular neighborhoods of $\crcomp$, and $V_i~(i=1,2)$ a connected component of $\nbh_i\setminus\crcomp$.
We will say that $V_1$ and $V_2$ represent the same {\em end} (of $\manif\setminus\crcomp$) at $\crcomp$ if there exists a regular neighborhood $\nbh \subset \interior(\nbh_1 \cap \nbh_2)$ of $\crcomp$ and a connected component $V$ of $\nbh\setminus\crcomp$ such that $V \subset V_1\cap V_2$.
The end determined by $V_i$ will be denoted by $\Kends{V}$.

Evidently, if $\Kends{V_1}=\Kends{V_2}$, then $\crcomp\cap\overline{V_1}=\crcomp\cap\overline{V_2}$.
We will say that an edge $\edgeR$ of $\crcomp$ {\em belongs} to $\Kends{V_i}$ if $\edgeR\subset\crcomp\cap\overline{V_1}$.

Let $\nbh$ be a regular neighborhood of $\crcomp$, $V$ a connected component of $\nbh\setminus\crcomp$, and $\difM:\manif\to\manif$ a diffeomorphism such that $\difM(\crcomp)=\crcomp$.
We will say that $\difM$ {\em preserves} $\Kends{V}$ if $\Kends{\difM(V)}=\Kends{V}$.

Also notice that $V\approx \aCircle\times [0,1)$, where $\aCircle\times 0$ corresponds to a connected component of $\partial\nbh$. So we can choose some orientation on $V$.
Let $\nbh_1 \subset \nbh\cap\difM(\nbh)$ be a regular neighborhood of $\crcomp$
and $V_1$ be the connected component of $\nbh_1\setminus\crcomp$ such that $\Kends{V_1}=\Kends{V}$.
Then $\difM$ {\em preserves the orientation} of $\Kends{V}$ if $\difM|_{V_1}:V_1\to V$ preserves orientation.

The following theorem is crucial for the proof of the exactness of Eq.~\eqref{equ:pi1Orbf-ex-seq} and will be applied to the case when $\crcomp$ is a critical component of a level-set of $\mrsfunc$.

\begin{figure}[ht]
\begin{tabular}{ccc}
\includegraphics[height=2cm]{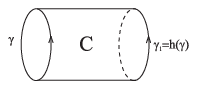} &
\includegraphics[height=2cm]{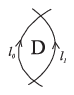} &
\includegraphics[height=2cm]{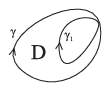} \\
a) & b) & c)
\end{tabular}
\protect\label{fig:cyl-bigon}
\end{figure}

\begin{thm}\label{th:ends}
Let $\crcomp\subset\interior\manif$ be a finite connected subgraph having no vertices of degrees $1$ and $2$ and such that every edge $\edgeR$ of $\crcomp$ belongs to precisely two ends of $\manif\setminus\crcomp$ at $\crcomp$.
Let $\difM:\manif\to\manif$ be a diffeomorphism such that $\difM(\crcomp)=\crcomp$ and $\difM$ preserves the ends at $\crcomp$ with their orientation.
Let $\cdifM$ be the combinatorial automorphism of $\crcomp$ induced by $\difM$.
Then each of the following conditions {\rm (1)} and {\rm (2)} implies that $\cdifM=\id_{\crcomp}$:

{\rm(1)} a regular neighborhood of $\crcomp$ is {\em flat}
(can be embedded in $\RRR^2$);

{\rm(2)} $\difM$ is isotopic to $\id_{\manif}$.
\end{thm}
\begin{proof}
First we prove the following two claims.

\begin{claim}\label{clm:one_edge}
The following statements are equivalent:

\condHHisID~$\cdifM=\id_{\crcomp}$;

\condHkeespEdges~$\difM$ preserves at least one edge of $\crcomp$ with its orientation;

\condHkeepsCycles~$\difM$ preserves every simple cycle $\acycle$  in $\crcomp$ with its orientation;
\end{claim}
\begin{proof}
Evidently, \condHHisID\ implies \condHkeespEdges\ and \condHkeepsCycles.

\condHkeespEdges$\Rightarrow$\condHHisID~
Let $\edgeR$ be an oriented edge of $\crcomp$ preserved by $\difM$ and $\pnt$ a vertex of $\edgeR$. Then $\difM(\pnt)=\pnt$.
Since $\difM$ preserves the ends at $\crcomp$ with their orientation, it follows that $\difM$ preserves the cyclic order of edges at $\pnt$, and thus preserves the edges incident to $\pnt$ with their orientation.
By the same arguments applied to another vertex of $\edgeR$ (if it exists), we obtain that $\cdifM$ is fixed on a neighborhood of $\edgeR$ in $\crcomp$.
It follows that the fixed-point set of $\cdifM$ on $\crcomp$ is open-closed.
From the connectedness of $\crcomp$ we get $\cdifM = \id_{\crcomp}$.

\condHkeepsCycles$\Rightarrow$\condHkeespEdges~
Suppose that $\difM$ preserves every simple cycle in $\crcomp$ with its orientation.
Since every edge of $\crcomp$ belongs to precisely two different ends at $\crcomp$ and $\crcomp$ has no vertices of degrees $1$ and $2$, it follows that there are two different simple cycles $\acycle_1$ and $\acycle_2$ such that either $\acycle_1\cap\acycle_2$ consists of a unique edge $\edgeR$ of $\crcomp$ or $\acycle_1\cap\acycle_2=\varnothing$ and there exists a unique simple path $\larc$ connecting these cycles.
Evidently, in both cases $\difM$ fixes some edge of $\crcomp$ with its orientation.
\end{proof}

\begin{claim}\label{clm:cycle_bnd_disk}
Let $\acycle$ be a simple cycle in $\crcomp$.
Suppose that $\acycle$ bounds a $2$-disk $D$ in $\manif$.
Then $\difM$ preserves $\acycle$ with its orientation.
\end{claim}
\proof
If $\interior D \cap \crcomp=\varnothing$, then $D$ determines a unique end at $\crcomp$, which is preserved with its orientation by $\difM$.
Otherwise, $D$ is a union of several $2$-disks of the previous type.
They are invariant under $\difM$, whence so is $D$.
In particular, $\difM$ preserves the boundary $\partial D=\acycle$ with its orientation.
\qed

Now we can complete Theorem~\ref{th:ends}.

(1)
Let $\atom$ be a regular neighborhood of $\crcomp$ in $\manif$, $i:\atom \subset \RRR^2$ an embedding, and $\datom\subset\interior\atom$ another regular neighborhood of $\crcomp$ such that $\difM(\datom)\subset\interior\atom$.
Then $\difM|_{\atom'}:\atom'\to\atom\subset\RRR^2$ extends to a diffeomorphism of $\RRR^2\to\RRR^2$ keeping the ends at $\crcomp$ with their orientation.
So we will assume that $\manif=\RRR^2$.
Then every simple cycle in $\crcomp$ bounds a $2$-disk in $\manif=\RRR^2$,
whence by Claim~\ref{clm:cycle_bnd_disk} $\difM(\acycle)=\acycle$.

(2) Suppose that $\difM$ is isotopic to $\id_{\manif}$.
Let $\acycle$ be an oriented simple cycle in $\crcomp$. 
Then $\acycle_1=\difM(\acycle)$ is also oriented cycle isotopic to $\acycle$.
In particular, homology classes of these cycles are equal: $[\acycle]=[\acycle_1]\in H_1(\manif)$.
By \condHkeepsCycles\ of Claim~\ref{clm:one_edge} it suffices to prove that $\acycle_1=\acycle$ and that $\difM$ preserves orientation of $\acycle$. 

By Claim~\ref{clm:cycle_bnd_disk} we can assume that $\acycle$ does not bound a $2$-disk in $\manif$.
Consider four cases.

(2.1) Suppose that $\acycle=\crcomp$, thus $\crcomp$ is a simple closed curve.
Since $\crcomp$ has no vertices of degree $2$, we see that $\crcomp$ has no vertices at all.
Moreover, $\manif\setminus\crcomp$ has two ends at $\crcomp$, whence $\crcomp$ is two-sided.
Since $\difM$ preserves the orientation of these ends, it follows that $\cdifM=\id_{\crcomp}$.

(2.2) Suppose that $\acycle\cap\acycle_1=\varnothing$. 
This is possible only if $\acycle$ is two-sided.
Then~\cite{Epstein} $\acycle \cup \acycle_1$ bounds a $2$-cylinder $C$ see in Figure~\ref{fig:cyl-bigon}a).

Notice that $\interior C\cap\crcomp\not=\varnothing$, otherwise the ends of $C$ determined by $\acycle$ and $\acycle_1$ will be not invariant under $\difM$.
Since $\crcomp$ is connected, it follows that $C\setminus\crcomp$ is a union of several open $2$-disks.
Each of them determines some end of $\manif\setminus\crcomp$ at $\crcomp$, and thus is invariant under $\difM$ with its orientation.
Therefore $\difM$ yields a preserving orientation homeomorphism of a subsurface $C\subset\manif$ that exchange boundary components $\acycle$ and $\acycle_1$ of $C$.
This contradicts to the assumption that $\difM$ is isotopic to $\id_{\manif}$.

(2.3) Suppose that $\acycle \cap \acycle_1 \not=\varnothing$ but $\acycle\not=\acycle_1$.
Then $\acycle_1$ and $\acycle$ must bound a {\em bigon}, i.e. a $2$-disk $D$ whose boundary consists either of two arcs $l_0\subset\acycle$ and $l_1\subset\acycle_1$ as in Figure~\ref{fig:cyl-bigon}b)
or of a union $\acycle\cup\acycle$ as in Figure~\ref{fig:cyl-bigon}c).
Then $D$ determines several ends at $\crcomp$ and therefore $\difM$ preserves $D$ with its orientation. 
Hence $\difM$ maps $l_0$ onto $l_1$ preserving their orientations.
But similarly to the previous case, $\difM$ preserves orientation of $D$ iff it maps $l_0$ onto $l_1$ or in case c) $\acycle$ onto $\acycle_1$ with opposite orientation.
We get a contradiction.

(2.4) Suppose that $\acycle=\acycle_1$ but $\difM$ reverses orientation of $\acycle$.
It follows that $[\acycle]=-[\acycle_1]$.
On the other hand $[\acycle]=[\acycle_1]$, whence $2[\acycle]=0$. 
Then we have two possibility.

Suppose that $\acycle$ bounds a subsurface $P$ in $\manif$.
Since $\acycle$ is isotopic to itself with opposite orientation, it follows that $P$ is a $2$-disk, whence by Claim~\ref{clm:cycle_bnd_disk} $\difM$ preserves orientation of $\acycle$.

Otherwise, $\acycle$ bounds no subsurfaces in $\manif$.
Since this cycle is isotopic to itself with opposite orientation, it follows that $\manif$ is a Klein bottle and $\acycle$ represents a unique element of order $2$ of $H_1\manif$.
Notice that $\manif\setminus\acycle$ is an open cylinder, whence the connected components of $\manif\setminus\crcomp$ are open $2$-disks.
Each of them determines some end of $\manif\setminus\crcomp$ at $\crcomp$.
Since $\difM$ preserves orientation of these ends, is also preserves the orientation of $\manif\setminus\crcomp$, and hence the orientation of $\acycle$.
\end{proof}
\chFig
\section{Homomorphisms $\izero$ and $\izerocr$.}\label{sect:ijzero}
Let $\izero:\pi_0\Stabf\to\pi_0\DiffM$ and $\izerocr:\pi_0\Stabf\to\pi_0\DiffMcr$ be the natural homomorphisms.
Then we have the following exact sequences:
\chEq\begin{equation}\label{equ:exact-Orbf}
\pi_1\DiffM \to \pi_1\Orbf \to \ker\izero \to 0
\end{equation}
\chEq\begin{equation}\label{equ:exact_Orbfcr}
\pi_1\DiffMcr \to \pi_1\Orbfcr \to \ker\izerocr \to 0.
\end{equation}
Thus in order to obtain estimates on the fundamental groups of orbits we should study the kernels of $\izero$ and $\izerocr$.
Notice that
$$
\ker\izero   \approx \pi_0\left(\Stabf\cap\DiffIdM   \right),
 \qquad
\ker\izerocr \approx \pi_0\left(\Stabf\cap\DiffIdMcr \right).
$$
Recall that we have a homomorphism $\actSRall:\Stabf \to \AutfR$, which by Lemma~\ref{lm:Dpart-Stabf}, reduces to a homomorphism $\actSR:\pi_0\Stabf \to \AutfR$.
So we obtain the following exact sequence
$$
   0                    \, \longrightarrow\,
   \ker\actSRall        \, \longrightarrow\,
   \Stabf               \, \stackrel{\actSRall}{\,\longrightarrow\,}
   \AutfR.
$$
Intersecting $\ker\actSRall$ and $\Stabf$ with $\DiffIdM$ and with $\DiffIdMcr$ and then taking $\pi_0$-groups we get the following two exact sequences:
\chEq\begin{equation}\protect\label{equ:izero-exact-seq}
   0                                 \, \longrightarrow\,
   \pi_0(\ker\actSRall\cap \DiffIdM) \, \longrightarrow\,
   \ker\izero                        \, \stackrel{\actSR}{\,\longrightarrow\,}
   \AutfR,
\end{equation}
\chEq\begin{equation}\protect\label{equ:izerocr-exact-seq}
   0                                 \, \longrightarrow\,
   \pi_0(\ker\actSRall\cap \DiffIdMcr) \, \longrightarrow\,
   \ker\izerocr                        \, \stackrel{\actSR}{\,\longrightarrow\,}
   \AutfR.
\end{equation}
\begin{prop}\label{pr:actSRall-Dpartf}
\chEq\begin{equation}\protect\label{equ:ker-cap-DidM}
\keractSRall \cap \DiffIdM = \Dpartfplus \cap \DiffIdM,
\end{equation}
\chEq\begin{equation}\protect\label{equ:ker-cap-DidMcr}
\keractSRall \cap \DiffIdMcr = \DpartfId,
\end{equation}
\chEq\begin{equation}\protect\label{equ:j0-inject}
 \Stabf\cap\DiffIdMcr = \StabIdf, \text{if either $\partial\manif\not=\varnothing$ or $\singf\not=\varnothing$.}
\end{equation}
\end{prop}
\proofstyle{Proof. Eq.~\eqref{equ:ker-cap-DidM}.}~Since $\keractSRall\supset\Dpartfplus$, we have to show that
$$\keractSRall \cap \DiffIdM \, \subset \, \Dpartfplus \cap \DiffIdM.$$
Suppose that $\difM\in\keractSRall \cap \DiffIdM$, so $\difM$ yields the identity automorphism of $\Reebf$ and is isotopic to $\id_{\manif}$.
We should prove that $\difM$ preserves leaves of the foliation $\partitf$ with their orientations.

Since $\difM$ trivially acts on $\Reebf$, it follows that $\difM$ preserves regular components of level-sets of $\mrsfunc$ and local extremes of $\mrsfunc$.
Moreover, as $\difM$ is isotopic to $\id_{\manif}$, it follows from (2) of Lemma~\ref{lm:keep_orient_of_levels} that $\difM$ preserves orientation of regular components of level-sets of $\mrsfunc$.

Let $\crcomp$ be a critical component of a level-set of $\mrsfunc$ containing a critical point of index $1$ and $\pnt\in\Reebf$ be a \Cvert-vertex corresponding to $\crcomp$.
Notice that the ends of $\manif\setminus\crcomp$ at $\crcomp$ corresponds to the edges of $\Reebf$ that are incident to $\pnt$.
Therefore $\difM$ preserves the ends of $\manif\setminus\crcomp$ at $\crcomp$ with their orientation and so the condition (2) of Theorem~\ref{th:ends} holds true.
Whence $\difM$ yields the identity automorphism of $\crcomp$, i.e. preserves the foliation $\partitf$ with its orientation.

\proofstyle{Eq.~\eqref{equ:ker-cap-DidMcr}.}
Evidently, $\keractSRall \cap \DiffIdMcr \supset \DpartfId$.
Suppose that $\difM\in\keractSRall \cap \DiffIdMcr$.
Then by Eq.~\eqref{equ:ker-cap-DidM}
$\difM\in \keractSRall \cap \DiffIdM \subset \Dpartfplus$.
Hence by Theorem~\ref{lm:twgr_Zintecnt}
$\difM$ is isotopic in $\Dpartfplus$ to a product of the internal Dehn twists
$\adifM=\Dtw_1^{m_1}\circ\cdots\circ\Dtw_{\intecnt}^{m_{\intecnt}}$
and therefore yields some automorphism $\adifM_{*}$ of $H_1(\manif\setminus\singf)$.
Since $\difM\in\DiffIdMcr$,
we see that $\adifM_{*}=\id$.
Whence $m_i=0$ for all $i$ and $\adifM=\id_{\manif}\in\DpartfId$.
Thus $\difM\in\DpartfId$.
\qed

\proofstyle{Eq.~\eqref{equ:j0-inject}.}
Evidently $\StabIdf \subset \Stabf\cap\DiffIdMcr$.
Conversely, suppose that $\difM\in\Stabf\cap\DiffIdMcr$.
Then $\actSRall(\difM)$ preserves the vertices of $\Reebf$ and yields
the identity isomorphism of $H_1\Reebf$, whence
$\actSRall(\difM)=\id_{\Reebf}$, i.e.$\!$ $\difM\in\keractSRall$.
Then from Eq.~\eqref{equ:ker-cap-DidMcr} we get:
$$\Stabf\cap\DiffIdMcr \, \subset \, \keractSRall  \cap \DiffIdMcr
\, = \,  \DpartfId \, = \, \StabIdf. \text{\qed}$$


\section{Proof of Theorem~\ref{th:hom-gr-orbits-c1}}\label{sect:proof-Th3}
Let $\mrsfunc:\manif\to\Psp$ be a Morse mapping and $\crpt{i}$ be the number of critical points of index $i$. Suppose that $\crpt{1}\geq 1$.

\proofstyle{(1).}~We have to show that $\Orbffcr$ is contractible.
By (3) of Theorem~\ref{th:loc-triv-fibering} and Whitehead theorem it suffices to prove that $\pi_i\Orbffcr=0$ for all $i\geq 1$.
By Theorem~\ref{th:StabIdf-hom-type} and Lemma~\ref{lm:DiffIdMcr-is-contractible} $\StabIdfcr$ and $\DiffMcr$ are contractible.
Moreover, by Eq.~\eqref{equ:j0-inject} $\ker\izerocr=0$.
Then from exact sequence of the $\Stabf$-fibration $\DiffMcr\to\Orbfcr$ we get $\pi_i\Orbffcr=0$ for $i\geq1$.
\qed

\proofstyle{(2).}~As $\crpt{1}\geq 1$, we have by Theorem~\ref{th:StabIdf-hom-type} that $\StabIdf$ is contractible.
By exact sequence of $\Stabf$-fibraiton $\DiffM\to\Orbf$ we get $\pi_k\DiffM\approx\pi_k\Orbf$ for $k\geq 2$.
Moreover, by Remark~\ref{rem:pi2D=0} (for $n=0$) we get $\pi_k\Orbf\approx\pi_k\DiffM\approx\pi_k\manif$ for $k\geq3$ and $\pi_2\Orbff=\pi_2\DiffM=0$.
\qed

\proofstyle{Eq.~\eqref{equ:pi1Orbf-ex-seq}.}~Denote $\twgrId=\pi_0(\Dpartfplus\cap\DiffIdM)$.
Thus $\twgrId$ consists of the ``relations'' between the internal Dehn twists in $\DiffIdM$.
We may regard $\twgrId$ as a subgroup of $\pi_0\Dpartfplus\approx\ZZZ^{\intecnt}$, whence $\twgrId\approx\ZZZ^{\rankpiO}$ for some $\rankpiO \leq \intecnt$.
Notice that for the case $\manif$ is a Klein bottle, it is possible that $\twgrId$ is not a direct summand of $\pi_0\Dpartfplus$.

Let $\kerAut = \parthom_1^{-1}(\twgrId) \subset \pi_1\Orbff$.
Then we have the following commutative diagram in which horizontal and vertical sequences are exact:

{\footnotesize
$$
\begin{CD}
 &    &   & &      0 & &  0 \\
 &    &   &  &        @V{}VV  @V{}VV  \\
0 @>{}>> \pi_1\DiffM / \IM(\ione)  @>{\prF}>> \kerAut  @>{\parthom_1}>> \twgrId\approx\ZZZ^{\rankpiO} @>{}>> 0 \\
&   &  @|   @V{}VV   @V{}VV \\
0 @>{}>> \pi_1\DiffM / \IM(\ione) @>{\prF}>> \pi_1\Orbff  @>{\parthom_1}>> \ker\izero @>{}>> 0 \\
  &    &   &  &    @V{\actSR\circ\parthom_1}VV     @V{\actSR}VV \\
  &    &   &   &               \grp    @=  \grp \\
  &    &   &   &   @V{}VV     @V{}VV \\
  &    &   &   &    0   & & 0
\end{CD}
$$
}
Here $\ione:\pi_1\Stabf \to \pi_1\DiffM$ is a natural homomorphism, and $\grp\subset\AutfR$ is the image of $\ker\izero$ under $\actSR$.
Then $\grp$ is a finite group.

Since $\prF(\pi_1\DiffM)$ is in the center of $\pi_1\Orbf$ (see (2) of Lemma~\ref{lm:prop-pi-1Dm}) and $\twgrId$ is free abelian, we see that $\kerAut$ is abelian
and that the upper sequence splits.
Since $\StabIdf$ is contractible, we have $\IM(\ione)=0$, whence $\kerAut\approx\pi_1\DiffM \oplus \twgrId$.
Thus the left vertical exact sequence coincides with Eq.~\eqref{equ:pi1Orbf-ex-seq}.
\qed

\proofstyle{Table~\ref{tbl:contractions}. Calculation the rank $\rankpiO$ of $\twgrId$.}
Suppose that $\mrsfunc$ is {\em not simple\/}.
Then by small perturbation of $\mrsfunc$ in a neighborhood of its critical level-sets we can find a {\em simple} Morse mapping $\bar\mrsfunc:\manif\to\Psp$ having same critical points as $\mrsfunc$.
By Theorem~\ref{th:twgr-descr} the group $\bar\twgrId=\pi_0(\Diff^{+}(\Delta_{\bar\mrsfunc})\cap\DiffIdM)$ is a free abelian group of some rank $\bar\rankpiO$.
\begin{claim}\label{clm:bark_k}
$\rankpiO \leq \bar\rankpiO$.
\end{claim}
\begin{proof}
Notice that \KR-graph $\Reebfb$ of $\bar\mrsfunc$ can be obtained by ``blowing up'' some internal vertices of $\Reebf$, i.e. replacing them with certain graphs, and we have a natural factorization  $q:\Reebfb\to\Reebf$ that shrinks these graphs into the corresponding points.

Let $\bar\edgeR$ be an edge of $\Reebfb$ such that $q(\bar\edgeR)$ is an edge of $\Reebf$.
Then both $\bar\edgeR$ and $q(\bar\edgeR)$ are internal or external simultaneously.
Moreover, since $\bar\mrsfunc$ differs from $\mrsfunc$ only near critical level-sets there is a Dehn twist $\Dtw:\manif\to\manif$ about $\bar\edgeR$ which is also a Dehn twist about $q(\bar\edgeR)$. In particular, $\Dtw$ preserves both mappings $\bar\mrsfunc$ and $\mrsfunc$.

Since $\pi_0\Dpartfplus$ is freely generated by internal Dehn twists, it follows that $\pi_0\Dpartfplus$ can be regarded as a subgroup of $\pi_0\Diff^{+}(\Delta_{\bar\mrsfunc})$.
Intersecting these groups with $\pi_0\DiffIdM$ we obtain that $\twgrId\subset\bar\twgrId$, whence $\rankpiO\leq\bar\rankpiO$.
\end{proof}

Suppose now that $\mrsfunc$ is simple.
Then there is a bijection between the critical points of $\mrsfunc$ and vertices of $\Reebf$.
\begin{defn}
An edge $\edgeR$ of $\Reebf$ is {\em contractible} if the following two conditions hold true:

{\rm(a)}
$\partial\edgeR$ consists of one \Evert-vertex and one \Cvert-vertex of degree $3$;

{\rm(b)}
let $\edgeR_1$ and $\edgeR_2$ be the other edges that are incident to the
vertex of degree $3$. Then $\edgeR_1 \not= \edgeR_2$ and at least one of
them is internal.
\end{defn}
Suppose that $\edgeR$ is a contractible edge.
Let us delete $\edgeR$ from $\Reebf$ and replace $\edgeR_1\cup\edgeR_2$ with one edge.
Denote the obtained graph by $\Graph_1$ and preserve the notations \Cvert-, \Dvert-, and \Evert- of its remaining vertices.
If $\Graph_1$ has a contractible edge, then we can repeat contractions as far as possible and obtain a {\em minimal\/} graph $\minGr$ having no contractible edges.

\begin{figure}[ht]
\begin{tabular}{ccccc}
\includegraphics[height=2.5cm]{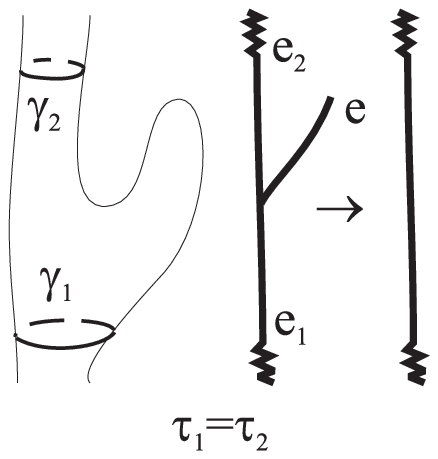} & \qquad &
\includegraphics[height=2.5cm]{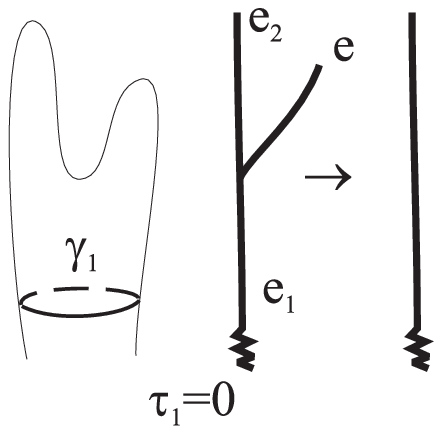} & \qquad &
\includegraphics[height=2.5cm]{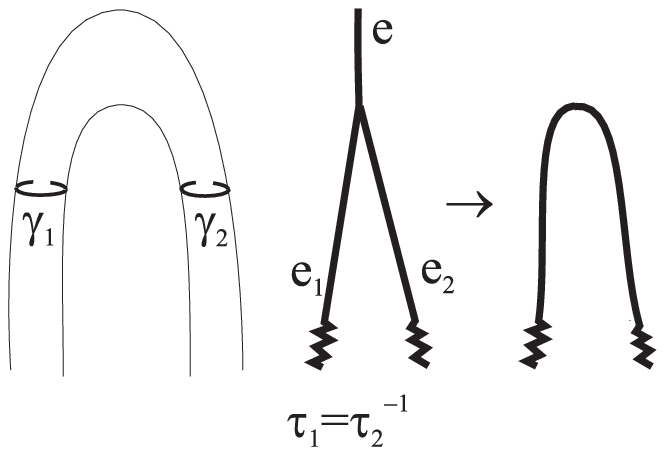} \\
 a) & \qquad & b)  & \qquad & c)
\end{tabular}
\caption{}
\protect\label{fig:contracting}
\end{figure}

As it is shown in Figure~\ref{fig:contracting}, every contractible edge $\edgeR$ corresponds to some relation in $\pi_0\DiffM$ between the Dehn twists about $\edgeR_1$ and $\edgeR_2$.
Moreover, the contraction of an edge yields a cancellation of the corresponding pair of singular points of $\partitf$.
It also decreases by $1$ the number of relations between the remaining internal Dehn twists, and the number of \Evert- and \Cvert-vertices.

\begin{figure}[ht]
\includegraphics[height=1.5cm]{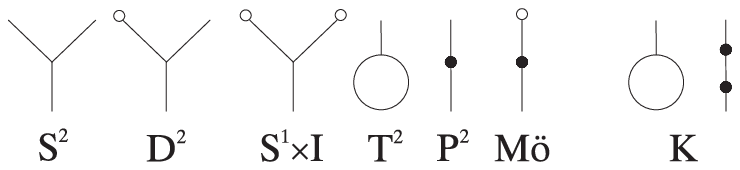}
\caption{}
\protect\label{fig:J0-generators}
\end{figure}

Let $\RdefectC$ and $\RdefectE$ be the numbers of \Cvert- and \Evert-vertices of $\minGr$ resp.
Since the numbers of \Cvert- and \Evert-vertices of $\Reebf$ are respectively $\crpt{1}$ and $\crpt{0}+\crpt{2}$, we see that $\minGr$ is obtained from $\Reebf$ by $\crpt{1}-\RdefectC=\crpt{0}+\crpt{2}-\RdefectE$ contractions.
Hence $\rankpiO \geq \crpt{1}-\RdefectC=\crpt{0}+\crpt{2}-\RdefectE$.
\chFig
Suppose that $\manif$ is of type 1 of Table~\ref{tbl:contractions}, i.e one of the surfaces $\sphere$, $D^2$, $\aCircle\times\Interv$, $\torus$, $\prjplane$ with or without holes.
Since $\crpt{1}>0$, $\minGr$ has a unique \Cvert-vertex and coincides with the corresponding graph in Figure~\ref{fig:J0-generators}.
Hence $\rankpiO \geq \crpt{1}-1$.
If $\manif=\torus$, then there is a unique internal Dehn twist which is non-isotopic to $\id_{\torus}$.
For other surfaces there are no internal edges at all.
Therefore $\rankpiO = \crpt{1}-1$ in all the cases.
\chFig
\begin{claim}
If $\manif$ is of type 2 or 3 of Table~\ref{tbl:contractions}, then 
$\rankpiO=\crpt{1}-\RdefectC=\crpt{0}+\crpt{2}-\RdefectE$.
\end{claim}
\begin{proof}
Suppose that $\manif$ is of types 2 or 3, i.e. it differs from the surfaces above. 
Then it is easy to see that the internal Dehn twists in $\minGr$ are mutually independent in $\DiffM$ since the corresponding simple closed curves are mutually disjoint and non-isotopic each to other.
Hence $\rankpiO=\crpt{1}-\RdefectC=\crpt{0}+\crpt{2}-\RdefectE$.
\end{proof}
If $\manif$ is of type 2, i.e. orientable, but not one of the surfaces above, then {\em all\/} \Evert-vertices can be removed by contractions.
Hence $\RdefectE=0$ and $\rankpiO=\crpt{0}+\crpt{2}$.

If $\manif$ is of type 3, i.e. non-orientable, but is neither $\prjplane$ with or without holes, then it is possible that we could not remove all \Evert-vertices.
The obstruction is that such vertices may be ``locked'' by vertices of degree $2$, see Figure~\ref{fig:J0-generators} for the Klein bottle $\Kleinb$.
Hence $\rankpiO=\crpt{0}+\crpt{2}-\RdefectE\leq\crpt{0}+\crpt{2}$.
The calculation of Table~\ref{tbl:contractions} is completed.
\qed

\subsection{\Homol-subgraph of $\Reebf$.}
For the proof of statement (3) of Theorem~\ref{th:hom-gr-orbits-c1} we have to study the group $\grp = \actSR(\ker\izero) \subset \AutfR$.

Suppose that $\difM\in\Stabf\cap\DiffIdM$, i.e. $[\difM]\in\ker\izero$ and $\actSRall(\difM)\in\grp$.
Then $\difM$ preserves every connected component of $\partial\manif$ and yields the identity automorphism of $H_1(\manif,\ZZZ)$.

Let $\AutfRHId$ be the subgroup of $\AutfR$ consisting of automorphisms trivially acting on $H_1(\Reebf,\ZZZ)$ and fixing every \Dvert-vertex of $\Reebf$.
Then $\actSRall(\difM) \in \AutfRHId$.
Thus $\grp\subset\AutfRHId.$

\begin{figure}[ht]
\includegraphics[width=0.3\textwidth]{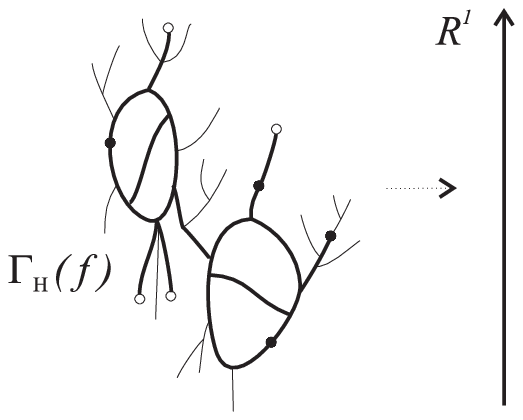}
\caption{}
\protect\label{fig:h1graph}
\end{figure}

It is easy to see that there exists a unique minimal connected subgraph $\aReebf$ of $\Reebf$ containing all simple cycles and all \Dvert-vertices of $\Reebf$, see Figure~\ref{fig:h1graph} in which $\aReebf$ is a subgraph in {\em bold}.
Then the inclusion $\aReebf\subset\Reebf$ yields an isomorphism $H_1\aReebf \approx H_1\Reebf$, therefore we will call $\aReebf$ an {\em \Homol-subgraph\/} of $\aReebf$.
Evidently that $\overline{\Reebf\setminus\aReebf}$ is a disjoint union of trees.

\chFig

Notice that every automorphism $\gdifM$ of $\AutfR$
that fixes $\aReebf$ point-wise belongs to $\AutfRHId$.
The following example shows that
the converse statement is not true.

\begin{exmp}
\em Let $\alpha,\beta$ be the standard generators of $\pi_1\torus=\ZZZ\oplus\ZZZ$.
Let $\mrsfunc:\torus\to\aCircle$ be a Morse function such that $\singf\not=\varnothing$, $\mrsfunc_{*}(\alpha)=1$, and
$\mrsfunc_{*}(\beta)=0\in\pi_1\aCircle=\ZZZ$.
For $n=2,3,\ldots$ let $\ttorus_{n}$ be a finite covering of $\torus$ corresponding to the subgroup of $\pi_1\torus$ generated by $\alpha^n$.
Let $\pr_n:\ttorus_{n}\to\torus$ be the covering projection,
$\mrsfunc_{n}=\mrsfunc\circ\pr_{n}$ the Morse function on $\ttorus_{n}$, and $\Reebfn$ the \KR-graph of $\mrsfunc_{n}$.
Evidently, $\Reebfn$ has a unique simple cycle which coincides with $\aReebfn$.

Let $\gdifM_{n}$ be the generator of the group of covering slices of $\ttorus_{n}$.
Then $\gdifM_{n}$ preserves $\mrsfunc_{n}$ and is isotopic to $\id_{\ttorus_{n}}$ but yields an automorphism of $\Reebfn$  which does not fix $\aReebfn$ point-wise.
\end{exmp}

\begin{prop}\label{pr:onto-AutfRHId}
{\rm (1)}~Suppose that $\mrsfunc$ is simple and
every automorphism $\gdifM\in\AutfRHId$ fixes $\aReebf$ point-wise.
Then $\actSR$ is onto, i.e.$\!$ $\grp = \AutfRHId$.

{\rm (2)}~If $\mrsfunc$ is generic, then $\AutfRHId$ is trivial, whence so is $\grp$.
\end{prop}

For the proof we need the following lemma.

\begin{lem}\label{lm:AutfRHId-prop}
If $\gdifM\in\AutfRHId$, then $\gdifM(\aReebf)=\aReebf$,
$\gdifM$ is fixed on $\overline{\aReebf\setminus\CyclesR}$,
where $\CyclesR$ is a union of simple cycles in $\Reebf$, and 
each of the following conditions implies that $\gdifM$ is fixed on $\aReebf$:

{\rm(1)} $\gdifM$ has a fixed point $\pnt\in\CyclesR$;

{\rm(2)} $\Reebf$ is a tree;

{\rm(3)} $\partial\manif\not=\varnothing$;

{\rm(4)} $\CyclesR$ has at least two connected components;

{\rm(5)} $\rank H_1\Reebf\geq 2$.

Thus, if $\gdifM$ is not fixed on $\aReebf$,
then $H_1\Reebf=\ZZZ$ and $\partial\manif=\varnothing$.
\end{lem}
\proof
Evidently, every simple cycle of $\Reebf$ is invariant under $\gdifM$.
Then so are the connected components of $\CyclesR$.
Moreover, since $\gdifM$ fixes \Dvert-vertices of $\Reebf$,
it follows that $\gdifM$ is fixed on simple paths
between the connected components of $\CyclesR$ and
between \Dvert-vertices of $\Reebf$ and $\CyclesR$.
Thus $\gdifM(\aReebf)=\aReebf$ and $\gdifM$ is fixed on
$\overline{\aReebf\setminus\CyclesR}$.

(1)~If $\gdifM(\pnt)=\pnt\in\CyclesR$, then $\gdifM$ is the identity on all edges of $\aReebf$ incident to $\pnt$.
Hence the fixed-points set of $\gdifM$ is open in $\aReebf$.
Since this set is also closed, it coincides with $\aReebf$.

(2)~If $\aReebf$ is a tree, then $\CyclesR=\varnothing$, whence
$\gdifM$ is fixed on $\aReebf$. 

(3) and (4)~In these cases $\gdifM$ is fixed on all simple paths
connecting \Dvert-vertices and the components of $\CyclesR$.
Then by (1) $\gdifM$ is fixed on $\aReebf$.

(5)~Suppose that $\rank H_1\Reebf\geq2$.
By (4) we can assume that $\CyclesR$ is connected.
Then there are two different simple cycles in $\Reebf$
whose intersection is either a point or a simple path.
In both cases $\gdifM$ is fixed on this intersection
and therefore on $\aReebf$.
\qed

\begin{figure}[ht]
\protect\label{fig:simple-atoms}
\begin{tabular}{ccc}
\includegraphics[height=2cm]{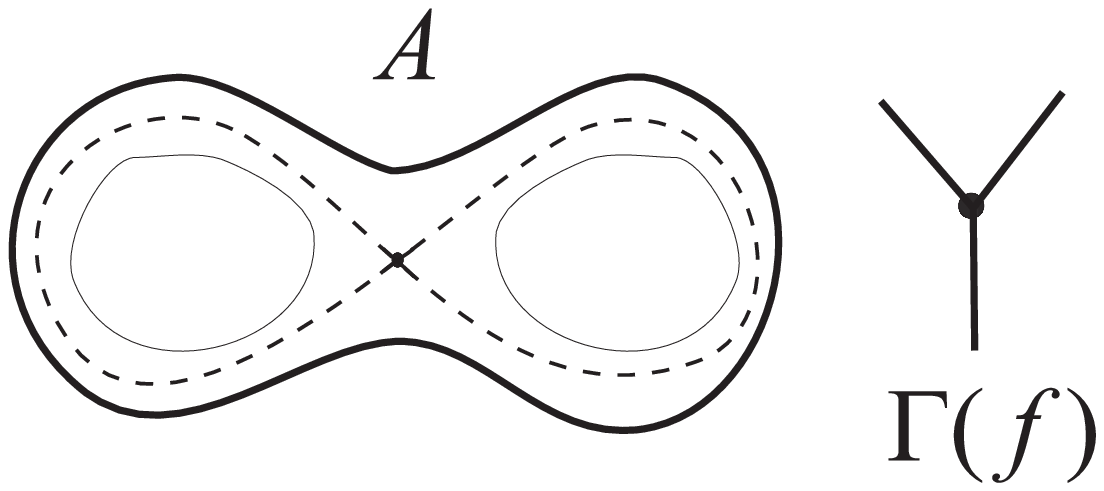}
& \qquad\qquad &
\includegraphics[height=2cm]{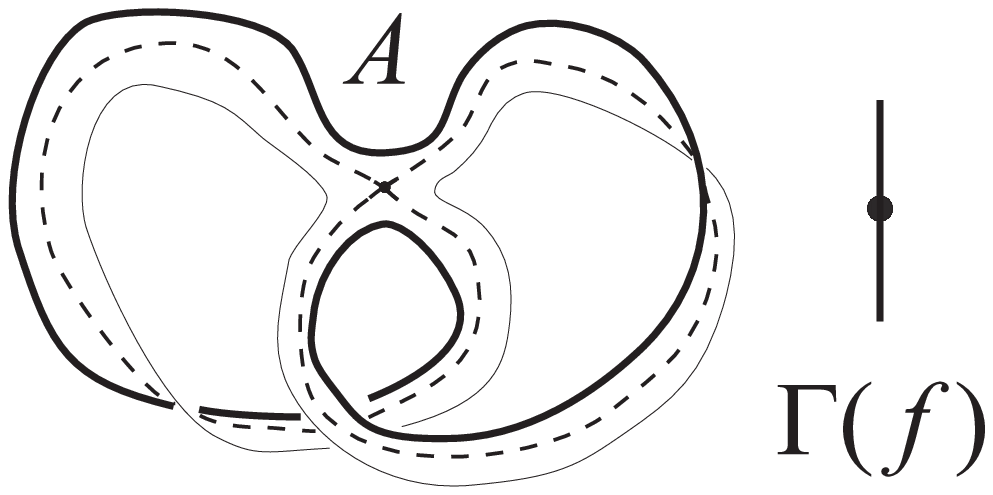} \\
a) & \qquad\qquad &  b)
\end{tabular}
\caption{}
\end{figure}

\proofstyle{Proof of Proposition~\ref{pr:onto-AutfRHId}}
(1)~We have to show that for every $\gdifM\in\AutfRHId$ there is a diffeomorphism $\difM\in\Stabf\cap\DiffIdM$ such that $\actSRall(\difM)=\gdifM$.
First we prove the following statement.
\begin{claim}\label{clm:dif-change-or}
There is a diffeomorphism $\bdifM_1$ of $\manif$ changing orientation of $\partitf$.
\end{claim}
\begin{proof}
Let $\crcomp$ be a union of all critical components of all critical level-set of $\mrsfunc$.
Since $\mrsfunc$ is simple, we see that every component of $\crcomp$ contains a unique critical point.
Let $\atom$ be a regular neighborhood of $\crcomp$.
Then every connected component $\datom$ of $\atom$ is either $2$-disk,
or has one of the forms shown in Figure~\ref{fig:simple-atoms}~a) and~b).
Evidently, $\datom$ admits an automorphism $\bdifM_{\datom}$
changing the orientation of the foliation $\partitf$ on $\atom$.

Since $\manif\setminus\atom$ is a disjoint union of
cylinders, it follows from (4) of Lemma~\ref{lm:ext-shift-func}, that
these automorphisms extends to a diffeomorphism $\bdifM_1$ of $\manif$
changing orientation of $\partitf$.
\end{proof}
\chFig
Since $\mrsfunc$ is simple, it follows from Lemma~\ref{lm:actSRall_Dpartf} that there exists a diffeomorphism $\adifM\in\Stabf$ such that $\actSRall(\adifM)=\gdifM$.
If $\adifM$ changes orientation of some leaf in $\prReeb^{-1}(\aReebf)$, then we replace $\adifM$ with $\adifM\circ\bdifM_1$, where $\bdifM_1$ is constructed in Claim~\ref{clm:dif-change-or}.
Thus we may assume that $\adifM$ preserves orientation of $\partitf$ in $\prReeb^{-1}(\aReebf)$.
Then similarly to the proof of Lemma~\ref{lm:ex-shift-for-difM}, we can find an isotopy of $\adifM$ in $\Stabf$ to a diffeomorphism $\adifM_1$ whose restriction to $\prReebf^{-1}(\aReebf)$ is a product of some internal Dehn twists $\adifM_1=\Dtw_1^{m_1}\circ \cdots \circ\Dtw_{\intecnt}^{m_{\intecnt}}$.
Notice that $\actSRall(\adifM) = \actSRall(\adifM\circ\adifM_1^{-1})$, so we set $\difM=\adifM\circ\adifM_1^{-1}$ and may assume that $\difM$ is the identity on $\prReeb^{-1}(\aReebf)$.

Since $\mrsfunc$ is \simple, and $\aReebf$ includes all simple cycles and \Dvert-vertices of $\Reebf$, we see that $\manif\setminus\prReeb^{-1}(\aReebf)$ is a disjoin union of $2$-disks.
Hence $\difM$ is isotopic to $\id_{\manif}$. 
This proves (1).
\qed

(2) Suppose that $\mrsfunc$ is \generic.
Let $\gdifM\in\AutfRHId$.
We have to show $\gdifM=\id_{\Reebf}$.
Since every critical level-set of $\mrsfunc$ contains a unique critical point, we see that $\gdifM$ fixes every vertex of $\Reebf$.
Then by (1) of Lemma~\ref{lm:AutfRHId-prop} $\gdifM$ is fixed on $\aReebf$.
Since $\overline{\Reebf\setminus\aReebf}$ is a disjoint union of trees and $\gdifM$ is fixed on the vertices of these trees, we obtain that $\gdifM$ is also fixed on their edges.
Thus $\gdifM=\id_{\Reebf}$.
\qed

\proofstyle{Proof of (3) of Theorem~\ref{th:hom-gr-orbits-c1}}.
Suppose that $\mrsfunc$ is generic.
Then $\grp=\id_{\Reebf}$ by (2) of Proposition~\ref{pr:onto-AutfRHId}, whence $\pi_1\Orbff\approx\pi_1\DiffM\oplus\ZZZ^{\rankpiO}$.
Since $\mrsfunc$ is also simple, it follows from statement (2) of this theorem and Table~\ref{tbl:hom-types} that the homotopy groups of $\Orbff$ for the surfaces of the left column of Table~\ref{tbl:hom_groups_of_orbits_c1} coincide with the homotopy groups of the corresponding spaces in the right column.
The construction of homotopy equivalences between these spaces is direct and we left it to the reader.
\qed


\section{Proof of Theorem~\ref{th:hom-gr-orbits-c10}}\label{sect:proof-orb-gr}
Suppose that $\mrsfunc$ has no critical points of index $1$.
Then by Theorem~\ref{th:StabIdf-hom-type}, $\StabIdf$ has the homotopy
type of $\aCircle$. We should describe the homotopy types of
$\Orbff$ and $\Orbffcr$. We will consider only the most non-trivial cases
\typeA\ and \typeEk.

\proofstyle{Type \typeA.}
Let $\mrsfunc:\sphere\to\Psp$ be a Morse mapping without critical points of index $1$.
We claim that $\Orbff$ is homotopy equivalent to $\sphere$.

Notice that the fibration of $\Diff(\sphere)$ over $\Orbf$ with fiber $\Stabf$ includes the fibration of $SO(3)$ over $\sphere$ with fiber $SO(2)$.
So we have the following commutative diagram:

\chEq\begin{equation}\protect\label{equ:fib-typeA}
\begin{array}{cccccccc}
 O(2)   & \to &  O(3)          & \to & O(3) / O(2)    & \approx & \sphere \\
      \cap   &     & \cap           &     & \cap   &     & \cap   &     \\
 \Stabf & \to & \Diff(\sphere) & \to & \Diff(\sphere)/\Stabf & \approx & \Orbf .
\end{array}
\end{equation}
It is easy to see that
the embedding $O(2)\subset\Stabf$ is a homotopy equivalence.
Moreover, by Smale~\cite{Smale} so is the embedding $O(3) \subset \Diff(\sphere)$.
Then from Eq.~\eqref{equ:fib-typeA} we obtain that the embedding
$\sphere\subset \Orbf$ yields isomorphisms of all homotopy groups.
Therefore it is also a homotopy equivalence.

\proofstyle{Type \typeEk.}
Let $\mrsfunc:\Kleinb\to\aCircle$ be a Morse map without critical points on the Klein bottle $\Kleinb$.
We represent $\Kleinb$ as the factor-space of $\torus = \RRR^2/\ZZZ^2$
by the relation $([x],[y]) \sim ([x+1/2], [-y])$ for $(x,y)\in\RRR^2$
and define $\mrsfunc$ by the formula
$\mrsfunc([x],[y]) = 2kx(\mathrm{mod} 1)$, for some $k=1,2,\ldots$

Since $\Kleinb$ is non-orientable,
$\StabIdf$ is contractible by Theorem~\ref{th:StabIdf-hom-type}.
Consider the embedding
$\difM:\aCircle\subset\DiffId(\Kleinb)$ defined by
$\difM(t)([x],[y]) = ([t+x],[y])$, where $t\in\RRR/\ZZZ=\aCircle$.
It is proved by C.~J.~Earle, A.~Schatz~\cite{EarleSchatz} and A.~Gramain~\cite{Gramain}
that $\difM$ is a homotopy equivalence.

Then from exact homotopy sequence we obtain that
$\pi_m\Orbff=0$ for $m\geq2$.
Consider the remaining part of this sequence:
$$ 0 \to \pi_1\Diff(\Kleinb) \to \pi_1\Orbff \to \pi_0\Stabf
 \to \pi_0\Diff(\Kleinb) \to \pi_0\Orbff \to 0.
$$
Recall, W.~B.~R.~Lickorish~\cite{Lickorish-non-or},
that $\pi_0\Diff(\Kleinb)\approx \ZZZ_2 \oplus\ZZZ_2$, where
one the former summand is generated
by the Dehn twist along some level-set of $\mrsfunc$
and the latter one is a ``$Y$-diffeomorphism'' in the terminology of~\cite{Lickorish-non-or}.

Moreover, it is easy to see that $\pi_0\Stabf \approx \ZZZ_k \oplus \ZZZ_2$,
where $\ZZZ_k$ is generated by the isotopy class of the diffeomorphism
$\difM_{1/k}([x],[y]) = ([1/k + x],[y])$,
and $\ZZZ_2$ is again generated by the Dehn twist along some level-set of $\mrsfunc$.
Then we obtain the following sequence:
$$ 0 \to \ZZZ
\to
\pi_1\Orbff \to \ZZZ_k \oplus \ZZZ_2
 \to \ZZZ_2 \oplus \ZZZ_2 \to \pi_0\Orbff \to 0.
$$
It follows that $\pi_1\Orbff\approx\ZZZ$, the homomorphism
$\pi_1\Diff(\Kleinb) \to \pi_1\Orbff$ coincides with the multiplication by $k$, and $\pi_0\Orbff\approx\ZZZ_2$.
Since $\Orbff$ is aspherical we obtain that
it is homotopy equivalent to $\aCircle$.
\qed


\section{Appendix. Orbits of tame actions}
\label{app:Frech_struct}
In this appendix we will prove Theorem~\ref{th:loc-triv-fibering}.
First we give one sufficient condition (Theorem~\ref{th:Orb_struct}) when a finite codimension orbit of a tame action of a tame Lie group $\Grp$ is a tame \Frechet\ manifold and the projection of $\Grp$ to this orbit is a locally trivial fibration.
Then we prove that $\DiffMcr$ and $\smoned$ are tame \Frechet\ manifolds (Section~\ref{sect:examples_tame_Frechet}).
Finally in Section~\ref{sect:proof_th_loctriv} we show that in the case of Theorem~\ref{th:loc-triv-fibering} that condition is satisfied.

We will assume that the reader is familiar with basic facts of \Frechet\ spaces. 
In particular, we will use differential calculus on \Frechet\ manifolds and the inverse function theorem, see~\cite{Hamilton}. 

\subsubsection{Derivatives}
Let $\FFr$ and $\GFr$ be \Frechet\ spaces, $\UFr \subset \FFr$ an open subset, and $\Bmp:\UFr\to\GFr$ a continuous mapping.
The {\em derivative of $\Bmp$ at a point $\fel\in\UFr$ in the direction $\thel\in\FFr$} is the limit:
$\Dmpv{\Bmp}{\fel}{\thel} = \lim\limits_{t\to0}\frac{1}{t}(\Bmp(\fel+t\thel)-\Bmp(\fel)).$
The mapping $\Bmp$ is {\em differentiable (of class $C^1$)} provided $\Dmpv{\Bmp}{\fel}{\thel}$ exists for all $\fel\in\UFr$, all $\thel\in\FFr$ and is continuous as a map $\Dm{\Bmp}:\UFr\times\FFr\to\GFr$. 
The mapping $\Bmp$ is {\em smooth ($C^{\infty}$)\/} if all its derivatives are differentiable.

Let $\FFr'$ another \Frechet\ space, $\UFr\subset\FFr\times\FFr'$ an open subset, and $\Bmp:\UFr\to\GFr$ a continuous mapping.
If $(\fel,\fel')\in\UFr$ and $(\thel,\thel')\in\FFr\times\FFr'$, then 
we shall write $\Dmpv{\Bmp}{\fel,\fel'}{\thel,\thel'}$ instead of $\Dmpv{\Bmp}{(\fel,\fel')}{(\thel,\thel')}$.

Moreover, when it is clear that $\thel\in\FFr$ and $\thel'\in\FFr'$, the partial derivatives $\Dmpv{\Bmp}{\fel,\fel'}{\thel,0}$ and $\Dmpv{\Bmp}{\fel,\fel'}{0,\thel'}$ of $\Bmp$ in the directions of $\FFr$ and $\FFr'$ will be denoted by $\Dmpv{\Bmp}{\fel,\fel'}{\thel}$ and $\Dmpv{\Bmp}{\fel,\fel'}{\thel'}$ respectively.
This would not lead to the confusion.
Notice that 
\chEq\begin{equation}\label{equ:D_d1_d2}
\Dmpv{\Bmp}{\fel_1,\fel_2}{\thel_1,\thel_2} =
\Dmpv{\Bmp}{\fel_1,\fel_2}{\thel_1} + \Dmpv{\Bmp}{\fel_1,\fel_2}{\thel_2}.
\end{equation}

\subsubsection{Tame group actions}
For the definition of tame linear mappings, tame smooth mapping, tame \Frechet\ spaces, and tame \Frechet\ manifolds we refer the reader to the parer of R.~Hamilton~\cite{Hamilton}.
In fact we will use the following statements about tame mappings: the composition of tame (linear or smooth) mappings is tame, a closed subspace of a tame \Frechet\ space is tame.
We will also use the inverse function theorem for tame mappings.
But we will never exploit the direct definition of a tame mapping.
Thus in a certain sense the proof belongs to categories theory.
 
Let $\Grp$ be a {\em tame Lie group}, i.e. a tame \Frechet\ manifold $\Grp$ which has a groups structure such that the multiplication map $\Gmult:\Grp\times\Grp\to\Grp$ and the inverse map $\Ginv:\Grp\to\Grp$ are smooth tame maps.

Let also $\XFr$ be a tame \Frechet\ manifold and $\Gact:\Grp\times\XFr\to\XFr$ a smooth tame left action of $\Grp$, thus $\Gact$ is a smooth tame mapping. 
Then the partial derivative of $\Gact$ with respect to $\XFr$ gives rise a smooth tame action
$\tGact:\Grp\times T\XFr\to T\XFr$ of $\Grp$ on the tangent bundle $T\XFr$ defined for $\gel\in\Grp$, $\xel\in\XFr$, and $\txel\in T_{\xel}\XFr$ by 
$\tGact(\gel,[\xel,\txel])=[\Gact(\gel,\xel), \Dmpv{\Gact}{\gel,\xel}{0,\txel}].$

In particular, since the left and right multiplications in $\Grp$ are also actions of $\Grp$ on itself, we have a {\em left\/} action $\tGmult:\Grp\times T\Grp\to T\Grp$ and a {\em right\/} action $\trGmult:T\Grp\times\Grp\to T\Grp$.
These multiplications commute in $\Grp$, therefore so do $\tGmult$ and $\trGmult$.

We will often use for such actions the following abbreviations:	
$$
\gel\grgr\gel_1 = \Gmult(\gel,\gel_1),  \qquad 
\gel\grtgr\tgel  =  \tGmult(\gel,\tgel),  \qquad 
\tgel\tgrgr\gel  =  \trGmult(\tgel,\gel),
$$
$$
\gel\grx\xel  =  \Gact(\gel,\xel),  \qquad 
\gel\grtx\txel  =  \tGact(\gel,\txel),
$$
where $\gel,\gel_1\in\Grp$, $\tgel\in T_{\gel}\Grp$, $\xel\in\XFr$, and $\txel\in T_{\xel}\XFr$.

Let $\UX\subset\XFr$ be an open set, and $\UG\subset\Grp$ an open neighborhood of $\un\in\Grp$.
A smooth tame map $\Gact:\UG \times \UX \to \XFr$ is a {\em smooth tame local action} if
\ $\Gact(\gel_1,\Gact(\gel_2,\xel))=\Gact(\gel_1\grgr\gel_2,\xel)$ \ 
provided $x,\Gact(\gel_2,\xel)\in\XFr$ and $\gel_1, \gel_2, \gel_1\grgr\gel_2\in\UG$.
If $\UG=\Grp$, then the action will be called {\em global}.

\subsubsection{Combined actions of two groups} \label{subsubsect:adj_actions}
Let $\Hgrp$ be another tame Lie group with multiplication $\Hmult:\Hgrp\times\Hgrp\to\Hgrp$ which will also be denoted by $\grgr$, and $\Hact:\Hgrp\times\XFr\to\XFr$ a left smooth tame action of $\Hgrp$ on $\XFr$ which we will write as $\Ract{}{}$.
Combining the actions of $\Hgrp$ and $\Grp$ we get the following smooth tame mapping
$\GHact:\Hgrp\times\Grp\times\XFr\to\XFr$ defined by 
$$\GHact(\hel,\gel,\xel)=\Hact(\hel,\Gact(\gel,\xel))=\Ract{(\gel\grx\xel)}{\hel}.$$

Let $\gel,\gel_1\in\Grp$, $\hel,\hel_1\in\Hgrp$, $\xel\in\XFr$,
$\tgel\in T_{\gel}\Grp$, and $\thel\in T_{\hel}\Grp$. 
Then the following identities can easily be verified:
$$
\begin{array}{lcl}
\Dmpv{\GHact}{\hel_1\grgr\hel, \gel,\xel}{\hel_1\grtgr\thel}
 & = & 
\hel_1\grtx\Dmpv{\GHact}{\hel,\gel,\xel}{\thel} \\ [1.5mm]
\Dmpv{\GHact}{\hel,\gel\grgr\gel_1,\xel}{\tgel\tgrgr\gel_1}
& = &
\Dmpv{\GHact}{\hel,\gel,\gel_1\grx\xel}{\tgel}
\end{array}
$$
They follow from the following ones:
$$
\GHact(\hel_1\grgr\hel,\gel,\xel)=\hel_1\grx\GHact(\hel,\gel,\xel), \qquad \GHact(\hel_1,\gel\grgr\gel_1,\xel)=\GHact(\hel,\gel,\gel_1\grx\xel)
$$
by differentiating the former of them in $\hel$ and the latter in $\gel$.
As a corollary we obtain the following relation:
\chEq\begin{equation}\label{equ:D_GHact}
\Dmpv{\GHact}{\hel,\gel,\xel}{\hel\grtgr\thel,\tgel\tgrgr\gel} =
\hel\grtx\Dmpv{\GHact}{\un,\un,\gel\grx\xel}{\thel,\tgel}.
\end{equation}

\subsubsection{Codimension of a point with respect to an action}
Let $\XFr$ be an open subset of a tame \Frechet\ space $\FFr$ and $\fel\in\XFr$.
Then we can identify the tangent space $\TfX$ with $\FFr$.
Let also $\Grp$ be a tame Lie group, $\UG$ an open neighborhood of $\un$ and $\Gact:\UG\times\XFr\to\FFr$ a smooth tame local action.
Since $\un\grx\fel=\fel$, we have the following linear mapping:
$$\Dmp{\Gact}{\un,\fel}:\TeG\to\TfX \equiv \FFr$$
obtained by differentiating $\Gact$ at $(\un,\fel)$ with respect to $\Grp$. Denote its image by $\ImD$. Thus $\ImD$ is a linear subspace of $\TfX$.

\begin{defn}
The number \ $\codm{\actGX}{\xpnt} \ = \ \dim_{\RRR} \, [ T_{\xpnt}\XFr \, /\, \ImD ] $ \ is called the {\em codimension\/} of $\xpnt\in\XFr$ with respect to the action $\actGX$.
\end{defn}

Suppose that $\codm{\actGX}{\fel} = \fcod < \infty$.
Thus $\ImD$ is a linear subspace of $\TfX$ of finite codimension $\fcod$.
Then we can find $\fcod$ linearly independent elements $\tdop{1},\ldots,\tdop{\fcod}\in\TfX$ which constitute a complementary basis to $\ImD$ in $\TfX$.
Denote $\tdop{}=(\tdop{1},\ldots,\tdop{\fcod})$.

Due to the identification $\TfX\equiv\FFr\supset\XFr$ the addition of linear combinations of $\tdop{i}$ to elements of $\XFr$ is well defined.
Since $\XFr$ is open in $\FFr$, there is a neighborhood $\UR\times\UX$ of $(0,\fel)$ in $\Rf\times\XFr$ such that the following smooth tame mapping $\isot:\UR\times\UX\to\XFr$ 
$$
\isot(\rrel,\xel)=\xel + \sum\limits_{j=1}^{\fcod}\rrel_{j}\tdop{j} =
\xel + \sprod{\rrel}{\tdop{}}.
$$
is well-defined for all $\xel\in\UX$, and $\rrel=(\rrel_1,\ldots,\rrel_{\fcod})\in\UR$.

Evidently, $\isot$ is a local action of $\Rf$ near $\fel\in\XFr$. 
Therefore it will be convenient sometimes to write down $\isot(\rrel,\xel)$ as $\Ract{\xel}{\rrel}$.
Also notice that 
\chEq\begin{equation}\label{equ:isot_prop}
\Dmpv{\isot}{\rrel,\xel}{\trrel,\txel}=\txel + \sprod{\trrel}{\tdop{}}. 
\end{equation}

Combining the local actions $\Gact$ and $\isot$ of $\Grp$ and $\Rf$ on $\UX$ we obtain the following smooth tame mapping
$\Asur:\UR\times\UG\times\XFr\to\XFr$ defined by:
\chEq\begin{equation}\label{equ:Asur_def}
\Asur(\rrel,\gel,\fel) = \isot(\rrel,\Gact(\gel,\fel)) = \Ract{(\gel\grx\fel)}{\rrel}
\end{equation}
for $\rrel\in\UR$ and $\gel\in\UG$.
This is an analogue of the mapping $\chi$ defined in~\cite[\S8.2]{Sergeraert}.
It is easy to see that 
\chEq\begin{equation}\label{equ:DAsur_DBsur_e_0}
 \Dmpv{\Asur}{\rrel,\gel,\fel}{\trrel,\tgel} = 
  \Dmpv{\Gact}{\gel,\fel}{\tgel} + \sprod{\trrel}{\tdop{}} 
\end{equation}

\begin{thm}\label{th:Orb_struct}
Let $\Gact:\Grp\times\XFr\to\XFr$ be a smooth tame action,
$\fel\in\XFr$ a point of finite codimension $\fcod$, $\tdop{1},\ldots,\tdop{\fcod}$ the complementary basis to $\ImD$, $\Orbit_{\fel}$ and $\Stab_{\fel}$ the $\Gact$-orbit and the $\Gact$-stabilizer of $\fel$.
 
Suppose that the tangent mapping $\Dmp{\Asur}{0,\un,\fel}:\Rf\times T_{\un}\Grp\to\TfX$
at $(0,\un,\fel)$ with respect to $\Rf\times\Grp$, see~Eq.\eqref{equ:DAsur_DBsur_e_0}, has a tame linear section
$$
\Lsect=(\LsectR, \LsectG) :\TfX\to \Rf\times T_{\un}\Grp,
$$
i.e. $\Dmp{\Asur}{0,\un,\fel}\circ\Lsect=\id(\TfX)$.
Thus for every $\txel\in\TfX$ we have
\chEq\begin{equation}\label{equ:repr_sect_L}
\txel = \Dmpv{\Gact}{\un,\fel}{\LsectG(\txel)} + \sprod{\LsectR(\txel)}{\tdop{}}.
\end{equation}

{\rm(1)} Then the natural projection $p:\Grp\to\Orbit_{\fel}$ is a locally trivial principal $\Stab_{\fel}$-fibration, and

{\rm(2)} the orbit $\Orbit_{\fel}$ is a smooth tame \Frechet\ manifold.
\end{thm}

For the proof of (1) it suffices to show that $p$ admits a local section at $\fel$, see Corollary~\ref{cor:p_loc_sect}.
To prove (2) we will show that there is a smooth tame embedding of a neighborhood $\UX$ of $\fel$ into some tame \Frechet\ space $\HFr$ such that the image of the intersection $\Orbit_{\fel}\cap\UX$ is an open subset of some closed linear subspace of $\HFr$, see Corollary~\ref{cor:Orf_Frechet}.

The crucial part of the proof is the following Lemma~\ref{lm:Asu_sect} below.
It follows the line of section~8.2 in~\cite{Sergeraert}. 
\begin{lem}\label{lm:Asu_sect}
The restriction 
\ $\Asur:\UR\times\UG\times\{\fel\} \to \XFr$ \
has a smooth tame local section at $\fel$, i.e. there exists a neighborhood $\UFr$ of $\fel$ in $\XFr$ and a smooth tame mapping $\sectGHX=(\sectGHXR,\sectGHXG):\UR\times\UG\to\UFr$ such that $\Asur\circ\Asect=\id_{\UFr}$, i.e. for every $\yel\in\UFr$ we have  
$$
\yel \,=\,
\Asur\circ\Asect(\yel) \,=\,
\Ract{\bigl(\sectGHXG(\yel)\grx\fel\bigr)}{\sectGHXR(\yel)} \,=\,
\sectGHXG(\yel)\grx\fel + \sprod{ \sectGHXR(\yel) }{\tdop{}}.
$$
\end{lem}
\begin{rem}\rm
This lemma is a combination of Theorem~4.2.5 and the first part of~\cite[Theorem~8.1.1]{Sergeraert}; 
the section $\sectGHX$ corresponds to the mapping $s_1$ of~\cite[Theorem 8.1.1]{Sergeraert}.
\end{rem}
\begin{proof}
We shall use the implicit function theorem for smooth tame mappings, see~\cite[Theorem~III.1.1.1]{Hamilton}.
 
We may assume that $\UG$ and $\UX$ are open subsets of some tame \Frechet\ spaces $\FG$ and $\FX$ respectively.
Then the tangent mapping of $\Asur:\UR\times\UG\to\XFr$ can be regarded as a {\em family of tame linear mappings}:
$$
\Dm{\Asur}:(\UR\times\UG)\times(\Rf\times\FG)\to\FX.
$$
Inverse function theorem claims that  $\Asur$ has a smooth tame local section provided $\Dm{\Asur}$ admits a smooth tame family of inverses:
$$
\Invm{\Asur}:(\UR\times\UG)\times\FX\to\Rf\times\FG.
$$
In other words we have to resolve smooth and tamely the following equation 
\chEq\begin{equation}\label{equ:DAsur_init_sect}
 \txel' =\Dmpv{\Asur}{\rrel,\gel,\fel}{\trrel',\tgel'}
\end{equation}
with respect to $\trrel'$ and $\tgel'$.
By~\eqref{equ:D_GHact} it can be rewritten as follows
$$
\rrel^{-1}\grtx\txel = \Dmpv{\Asur}{\,0,\,\un,\,\gel\grx\fel\,}{\,\rrel^{-1}\grtx\trrel',\,\tgel'\tgrgr\gel^{-1}\,}.
$$
Suppose that we can resolve the following equation smooth and tamely in $\tgel$ and $\trrel$:
\chEq\begin{equation}\label{equ:txel_need}
\txel = \Dmpv{\Asur}{0,\un,\gel\grx\fel}{\trrel,\tgel} = \Dmpv{\Gact}{\un,\gel\grx\fel}{\tgel} + \sprod{\trrel}{\tdop{}}.
\end{equation}
i.e. there is a smooth tame mapping $\Bsect=(\BsectR,\BsectG):\FX\times\UG\to\Rf\times\FG$ such that 
$$
\txel = \Dmpv{\Asur}{0,\un,\gel\grx\fel}{\BsectR(\txel,\gel),\BsectG(\txel,\gel)}.
$$
Then the solution of~\eqref{equ:DAsur_init_sect} can be given by the following formulas: 
$$
\trrel' = \rrel\grtx\BsectR(\rrel^{-1}\grtx\txel',\gel),
\qquad
\tgel' = \BsectG(\rrel^{-1}\grtx\txel',\gel)\tgrgr\gel.
$$
Thus we are reduced to resolve~\eqref{equ:txel_need}.
These arguments constitute Theorem~4.2.5. of~\cite{Sergeraert}.
\begin{claim}\label{clm:expr_for_tdop}
Denote $\gel^{-1}\grtx\tdop{}=(\gel^{-1}\grtx\tdop{1},\ldots,\gel^{-1}\grtx\tdop{\fcod}).$
There is a neighborhood $\UNbh_{\un}$ of $\un$ in $\UG$ such that for $\gel\in\UNbh_{\un}$ 
\chEq\begin{equation}\label{equ:tdop_g_1_tdop}
\tdop{} =
- \Dmpv{\Gact}{\un,\fel}{\LsectRg{\gel}^{-1}\cdot\LsectG(\gel^{-1}\grtx\tdop{})} +
\LsectRg{\gel}^{-1} \cdot (\gel^{-1}\grtx\tdop{}),
\end{equation}
where $\LsectRg{\gel}$ is a real non-singular $\fcod\times\fcod$ matrix which smooth and tamely depends on $\gel$.
\end{claim}
\begin{proof}
Applying~\eqref{equ:repr_sect_L} to $\gel^{-1}\grtx\tdop{i}$ for each
$i=1,\ldots\fcod$ we get the following system of equations:
$$
\left\|\begin{array}{c}
\gel^{-1}\grtx\tdop{1}\\ \cdots \\ \gel^{-1}\grtx\tdop{\fcod}
\end{array}\right\|
=
\Dmpv{\Gact}{\un,\fel}{
\left\|\begin{array}{c}
\LsectG(\gel^{-1}\grtx\tdop{1}) \\ \cdots \\ \LsectG(\gel^{-1}\grtx\tdop{\fcod})
\end{array}\right\|
} +
\left\|\begin{array}{c}
\LsectR(\gel^{-1}\grtx\tdop{1}) \\ \cdots \\ \LsectR(\gel^{-1}\grtx\tdop{\fcod})
\end{array}\right\|
\cdot
\left\|\begin{array}{c}
\tdop{1}\\ \cdots \\ \tdop{\fcod},
\end{array}\right\|
$$
which can be written in a simpler form:
\chEq\begin{equation}\label{equ:g_1_tdop}
\gel^{-1}\grtx\tdop{} =
\Dmpv{\Gact}{\un,\fel}{\LsectG(\gel^{-1}\grtx\tdop{})} +
\underbrace{ \LsectR(\gel^{-1}\grtx\tdop{}) }_{\LsectRg{\gel}} \,\tdop{}.
\end{equation}
Notice that $\LsectR(\gel^{-1}\grtx\tdop{i})$ is a real $\fcod$-vector, whence $\LsectR(\gel^{-1}\grtx\tdop{})$ is a real $\fcod\times\fcod$-matrix. 
Let us denote it by $\LsectRg{\gel}=\LsectR(\gel^{-1}\grtx\tdop{})$.

Evidently, $\LsectRg{\un}$ is the identity matrix.
Then there is a neighborhood $\UNbh_{\un}$ of $\un$ in $\UG$ on which $\LsectRg{\gel}$ is non-singular.
Hence, if $\gel\in\UNbh_{\un}$, then~\eqref{equ:g_1_tdop} is equivalent to~\eqref{equ:tdop_g_1_tdop}.
\end{proof}
Now we can complete Lemma~\ref{lm:Asu_sect}.
Let $\gel\in\UNbh_{\un}$, $\txel\in\FX$. Denote $\txel_1=\gel^{-1}\grtx\txel$.
Then 
\chEq\begin{equation}\label{equ:txel1}
\begin{array}{lcl}
\txel_1 & \stackrel{\eqref{equ:repr_sect_L}}{=\!=\!=} &
\Dmpv{\Gact}{\un,\fel}{\LsectG(\txel_1)}+\sprod{\LsectR(\txel_1)}{\tdop{}}
= \\ [1.5mm]
& \stackrel{\eqref{equ:tdop_g_1_tdop}}{=\!=\!=} &
 \Dmpv{\Gact}{\un,\fel}
 {\ \underbrace{ \LsectG(\txel_1) - \LsectRg{\gel}^{-1}\cdot\LsectG(\gel^{-1}\grtx\tdop{})}_{\MsectG(\txel,\gel)}\ } \ + \\ [2mm]
& & \hfill +\ \sprod{\ \underbrace{\LsectR(\txel_1)\cdot (\LsectRg{\gel}^{-1})^{t}}_{\MsectR(\txel,\gel)}\
}{\gel^{-1}\grtx\tdop{}}.  
\end{array}
\end{equation}
The expressions $\MsectG(\txel,\gel)$ and $\MsectR(\txel,\gel)$ in this formula are compositions of smooth tame mappings and therefore they are smooth tame themselves.

Then from~\eqref{equ:txel1} we get the following representation for $\txel$:
$$ 
\begin{array}{rl}
\txel=\!\gel\grtx\txel_1 \! \stackrel{\eqref{equ:txel1}}{=\!=\!=} \! &
\gel\grtx
\bigl[ \,
\Dmpv{\Gact}{\un,\fel}{\MsectG(\txel,\gel)}+\sprod{\MsectR(\txel,\gel)}{\gel^{-1}\grtx\tdop{}}
 \, \bigr] = \\ [1.5mm]
 \stackrel{\eqref{equ:D_GHact}}{=\!=\!=} &
\Dmpv{\Gact}{\un,\gel\grx\fel}{\gel\grtgr\MsectG(\txel,\gel)\tgrgr\gel^{-1}}+\sprod{\MsectR(\txel,\gel)}{\tdop{}},
\end{array}
$$ 
which has the form~\eqref{equ:txel_need} with smooth tame $\tgel$ and
$\trrel$.
This proves Lemma~\ref{lm:Asu_sect}.
Notice that we have used here~\eqref{equ:D_GHact} for $\Hgrp=\Grp$.
\end{proof}
Since our considerations are local, we may further assume that in Lemma~\ref{lm:Asu_sect} $\UNbh_{\un}=\UG$ and that $\Asect$ is defined on all of $\UX$, i.e. $\UFr=\UX$.
\begin{cor}\label{cor:tdop_ind_gImD}
If $\gel$,$\hel$, and $\gel\grgr\hel\in\UG$, then $\tdop{i}$ are also independent over the image $\ImDx{\gel,\hel\grx\fel}$ of the tangent linear map 
$$
  \Dmp{\Gact}{\gel,\hel\grx\fel}:T_{\gel}\Grp\to T_{(\gel\grgr\hel)\grx\fel}\XFr.
$$
\end{cor}
\begin{proof}
First notice that 
$\Dmpv{\Gact}{\gel,\hel\grx\fel}{\tgel}=(\gel\grgr\hel)\grtx\Dmpv{\Gact}{\un,\fel}{\tgel_1}$
for some $\tgel_1\in T_{\un}\Grp$.
Hence $\ImDx{\gel,\hel\grx\fel}=(\gel\grgr\hel)\grtx\ImD$.
Thus it suffices to show that $\tdop{i}$ are independent over $\gel\grtx\ImD$ for every $\gel\in\UG$.

Suppose that for some $\rrel\in\Rf$ and $\tgel\in T_{\un}\Grp$ we have
$$
\sprod{\rrel}{\tdop{}}=\gel\grtx\Dmpv{\Gact}{\un,\fel}{\tgel},
$$
or equivalently, $\sprod{\rrel}{\gel^{-1}\grtx\tdop{}}=\Dmpv{\Gact}{\un,\fel}{\tgel}$.
Then from~\eqref{equ:g_1_tdop} it follows that
$$ 
\sprod{\rrel}{\LsectRg{\gel}\cdot\tdop{}}= \sprod{\rrel}{\gel^{-1}\grtx\tdop{}} - \Dmpv{\Gact}{\un,\fel}{\tgel_1} = \Dmpv{\Gact}{\un,\fel}{\tgel-\tgel_1}
$$ 
for some $\tgel_1\in T_{\un}\Grp$.
Since $\LsectRg{\gel}$ is non-singular for $\gel\in\UG$ and $\tdop{}$ are independent over $\ImD$, it follows that $\LsectRg{\gel}\cdot\tdop{}$ are also independent over $\ImD$, whence $\rrel=0$.
\end{proof}
\begin{lem}\label{lm:fl_gf}
Let $\gel_1\in\UG$. Suppose that $\hrrel(t):I\to\Rf$ and $\hgel(t):I\to\UG$ are smooth paths such that 
$\hrrel(0)=0$, $\hgel(0)=\un$, and
\chEq\begin{equation}\label{equ:f_lt__gt_f}
\Ract{(\gel_1\grx\fel)}{\hrrel(t)}=\hgel(t)\grx\fel.
\end{equation}
Then $\hrrel(t)$ is in fact a constant path.
In particular, $\hrrel(0)=\hrrel(1)=0$.
\end{lem}
\begin{proof}
Differentiating~\eqref{equ:f_lt__gt_f} in $t$ we obtain:
$$ 
\Dmpv{\isot}{\gel_1\grx\fel,\hrrel(t)}{0,\hrrel'(t)} =
\sprod{\hrrel'(t)}{\tdop{}} = \Dmpv{\Gact}{\hgel(t),\fel}{\hgel'(t)}.
$$ 
The right term belongs to $\ImDx{\hgel(t),\fel}$
On the other hand, since $\tdop{i}$ are independent over $\ImDx{\hgel(t),\fel}$, 
we obtain that $\hrrel'(t)=0$ for all $t\in I$, whence $\hrrel(t)$ is a constant path.
\end{proof}
\begin{rem}
\rm The arguments of Lemma~\ref{lm:fl_gf} are used at the end of the proof of~\cite[Proposition~9.2.3]{Sergeraert}.
\end{rem}
\begin{lem}\label{lm:f_lphi_gf_l0}
{\rm\cite[Proposition~9.2.2]{Sergeraert}}
There is a neighborhood $\UNbh_{\un}\subset\UG$ of $\un\in\Grp$ such that $\Ract{(\gel\grx\fel)}{\rrel}\in\UX\cap\Orbit_{\fel}$ if and only if $\rrel=0$.
\end{lem}
\begin{proof}
In~\cite[Proposition~9.2.2]{Sergeraert} this statement was established via section $s_2$ (and only for $\gel=\un$). 
Our proof is similar, but is in opposite based on the section $\Asect$ which is analogue of $s_1$.

{\em Sufficiency.} Evidently, $\Ract{(\gel\grx\fel)}{0}=\gel\grx\fel\in\UX\cap\Orbit_{\fel}$.

{\em Necessity.}
Suppose that for each neighborhood $\UNbh_{0}\times\UNbh_{\un}$ of $(0,\un)$ in $\Rf\times\UG$ there exist $\rrel\in\UNbh_{0}$ and $\gel,\gel_1\in\UNbh_{\un}$ such that $\Ract{\gel_1\grx\fel}{\rrel}=\gel\grx\fel$.

We will show now that there are smooth paths $\hrrel(t):I\to\Rf$ and $\hgel(t):I\to\UG$ satisfying~\eqref{equ:f_lt__gt_f} and such that $\hgel(0)=\un$, $\hgel(1)=\gel$, $\hrrel(0)=0$, and $\hrrel(1)=\rrel$.
Then by Lemma~\ref{lm:fl_gf}, we will get $\rrel=\hrrel(1)=\hrrel(0)=0$.

Choose smooth paths $\gel(t):I\to\UG$ and $\rrel(t):I\to\Rf$ such
that $\gel(0)=\un$, $\gel(1)=\gel$, $\rrel(0)=0$, and $\rrel(1)=\rrel$.
Since $(\gel,\rrel)$ can be chosen arbitrary close to $(\un,0)$ we may also assume that $\ael_t=\Ract{\gel_1\grx\fel}{\rrel(t)}\in\UX$ for all $t\in I$.
Thus $\Asur\circ\Asect(\ael_t)=\ael_t$, i.e. \
$\ael_t = \Ract{\gel_1\grx\fel}{\rrel(t)} = \Ract{ \AsectG(\ael_t)\grx\fel }{ \AsectR(\ael_t) }$
whence
$$
\Ract{\gel_1\grx\fel}{(\rrel(t)- \AsectR(\ael_t))} = \AsectG(\ael_t)\grx\fel.
$$
So we may put $\hrrel(t)=\rrel(t)- \AsectR(\ael_t)$ and $\hgel(t)=\AsectG(\ael_t)$.
Then~\eqref{equ:f_lt__gt_f} holds true.
\end{proof}
\begin{cor}\label{cor:AR_gr_l__l}
If \ $\Ract{\gel\grx\fel}{\rrel}\,\in\UX$, then $\AsectR(\Ract{\gel\grx\fel}{\rrel})=\rrel$.
\end{cor}
\begin{proof}
Applying $\Asect$ to $\xel = \Ract{\gel\grx\fel}{\rrel}$ we get 
\chEq\begin{equation}\label{equ:gf_l__AG_AR}
\Ract{\gel\grx\fel}{\rrel} = 
\Ract{\, \AsectG(\xel)\grx\fel\, }{\, \AsectR(\xel)\, },
\end{equation}
whence by Lemma~\ref{lm:f_lphi_gf_l0} $\rrel=\AsectR(\Ract{\gel\grx\fel}{\rrel})$.
\end{proof}
\begin{cor}\label{cor:p_loc_sect}
The projection $p:\UG\to\Orbit_{\fel}$ has a smooth tame local section at $\fel$ defined by $\hel\mapsto\AsectG(\hel)$ for $\hel\in\Orbit_{\fel}$.
\end{cor}
\begin{proof}
If $\hel=\gel\grx\fel\in\Orbit_{\fel}$, then by Corollary~\ref{cor:AR_gr_l__l} we obtain that  $\AsectR(\gel\grx\fel)=0$, whence by~\eqref{equ:gf_l__AG_AR} $\hel=\AsectG(\hel)\grx\fel$.
\end{proof}
\begin{cor}\label{cor:Orf_Frechet}
There is a smooth tame embedding $\emb$ of $\UX$ into the tame \Frechet\ space $\ImD\times\Rf\times\Rf$ such that $\emb(\Orbit_{\fel}\cap\UX)$ is a neighborhood of $(0,0,0)$ in the closed linear subspace $\ImD\times0\times0$.
Hence $\Orbit_{\fel}$ admits the structure of a smooth tame \Frechet\ manifold.
\end{cor}
\begin{proof}
Let $\hel\in\UX$.
Then $\hel=\Ract{ \AsectG(\hel)\grx\fel }{\AsectR(\hel)} = 
\AsectG(\hel)\grx\fel + \sprod{\AsectR(\hel)}{\tdop{}}$.
Moreover, we can decompose $\AsectG(\hel)\grx\fel = \fel + (\AsectG(\hel)\grx\fel - \fel) $ via $\Lsect$, see~\eqref{equ:repr_sect_L}. This gives us the following representation:
$$
\hel = \fel + 
\overbrace{
\underbrace{  \Dmpv{\Gact}{\un,\fel}{ \Lsect(\AsectG(\hel)\grx\fel-\fel)} }_{\emb_{0}(\hel)}   + 
 \sprod{\! \underbrace{ \Lsect(\AsectG(\hel)\grx\fel-\fel) }_{\emb_{1}(\hel)} \!}{\!\tdop{}\!}  
}^{\AsectG(\hel)\grx\fel - \fel} +
 \sprod{\!\underbrace{ \AsectR(\hel)}_{\emb_{2}(\hel)}\!}{\!\tdop{}\!}.
$$
Thus we obtain a smooth tame mapping, see Figure~\ref{fig:orb_decomp}:
$$
\emb=(\emb_{0},\emb_{1},\emb_{2}):\UNbh \to \ImD\times\Rf\times\Rf.
$$%
\chFig\begin{figure}[ht]
\includegraphics[height=3.5cm]{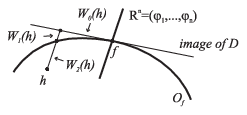}
\caption{}
\protect\label{fig:orb_decomp}
\end{figure}

It follows from the above decomposition that the following smooth tame mapping $\invemb: \ImD\times\Rf\times\Rf\to\UX$ defined by $$\invemb(\del,\muel,\rrel)=\fel+\del+\sprod{\muel+\rrel}{\tdop{}}$$
is a left inverse to $\emb$, i.e. $\invemb\circ\emb(\hel)=\hel$.
Hence $\emb$ is a tame smooth embedding.

It remains to show that $\emb_{0}:\Orbit_{\fel}\cap\UX\to\ImD$ is a homeomorphism onto a neighborhood of $0\in\ImD$. 
\begin{claim}
$\emb_{0}|_{\Orbit_{\fel}\cap\UX}$ is injective.
\end{claim}
\begin{proof}
Let $\hel_i=\gel_i\grx\fel\in\Orbit_{\fel}\cap\UX$, $i=1,2$.
Then by Corollary~\ref{cor:AR_gr_l__l}, $\emb_{2}(\gel_1\grx\fel)=\AsectR(\gel_2\grx\fel)=0$, whence 
$\gel_i\grx\fel = \fel + \emb_{0}(\hel_i) + \sprod{ \emb_{1}(\hel_i) }{\tdop{}}$.

Suppose that $\emb_{0}(\hel_1)=\emb_{0}(\hel_2)$. 
Then it follows that 
$$\gel_1\grx\fel = \gel_2\grx\fel + \sprod{ \emb_{1}(\hel_2)-\emb_{1}(\hel_1) }{\tdop{}}.$$
Moreover, from Lemma~\ref{lm:f_lphi_gf_l0} we also get $\emb_{1}(\hel_1)=\emb_{1}(\hel_2)$.
Thus $\emb_{i}(\hel_1)=\emb_{i}(\hel_2)$ for $i=0,1,2$, whence $\hel_1=\hel_2$.
\end{proof}
\begin{claim}
$\emb_{0}|_{\Orbit_{\fel}\cap\UX}$ is onto and its inverse $\invemb:\ImD\to\Orbit_{\fel}$  is given by the formula $\invemb(\del)=\AsectG(\fel+\del)\grx\fel$.
\end{claim}
\begin{proof}
Let $\del\in\ImD$. Then we have two representations:
$$
\begin{array}{rcccccc}
\AsectG(\fel+\del)\grx\fel & = & \fel & + & \emb_{0}(\AsectG(\fel+\del)\grx\fel) & + & 
\sprod{ \emb_{1}(\AsectG(\fel+\del)\grx\fel) }{\tdop{}}, \\ [1,5mm]
\AsectG(\fel+\del)\grx\fel & = & \fel & + & \del                                 & - & 
\sprod{ \AsectR(\fel+\del)\grx\fel }{\tdop{}}. 
\end{array}
$$
The second relation is just the application of $\Asect$ to $\fel+\del$.
Since $\tdop1,\ldots,\tdop{\fcod}$ are independent over $\ImD$ it follows that 
$\del=\emb_{0}(\AsectG(\fel+\del)\grx\fel)$ and 
$\AsectR(\fel+\del)\grx\fel=-\emb_{1}(\AsectG(\fel+\del)\grx\fel)$.
\end{proof}
This completes Corollary~\ref{cor:Orf_Frechet} and Theorem~\ref{th:Orb_struct}.
\end{proof}


\subsection{Examples of tame \Frechet\ manifolds}\label{sect:examples_tame_Frechet}
We will show here that $\DiffMcr$ and $\smoned$ (see Section~\ref{sect:intro}) are tame \Frechet\ manifolds.

\subsubsection{The space $\Cinf(\Mman,\Nman)$}
Let $\Mman$ and $\Nman$ be smooth finite-dimensional manifolds such that $\Mman$ is compact and $\partial\Nman=\emptyset$.
The following statement is well-known, e.g.~\cite{Hamilton, KrieglMichor}.
Nevertheless we briefly recall this construction.
\begin{lem}\label{lm:smmn_tame_Frechet_man}
{\rm e.g.~\cite[Example~I.4.1.2 and Theorem~II.2.3.1]{Hamilton}}
The space $\smmn$ with $C^{\infty}$ topology is a tame \Frechet\ manifold whose tangent space at a point $\mrsfunc\in\smmn$ is the space $\ssmftn$ of sections of the pullback of the tangent bundle $T^{*}\Nman$ under $\mrsfunc$.
\end{lem}
\begin{proof}
Fix some Riemannian metric $d$ on $\Nman$.
Then there exists a neighborhood $\WNbh\subset T\Nman$ of the zero section on which the exponential mapping $\exp:\WNbh\to\Nman\times\Nman$ is well-defined. Recall that $\, \exp(\xel,\xrel) = (x,v_{\xrel})\,$,
where $\xel\in\Nman$, $\xrel\in T_{\xel}\Nman$, and $v_{\xrel}\in\Nman$ is the end-point of the geodesic of length $\|\xrel\|$ starting at $\xel$ in the direction $\xrel$.
Decreasing $\WNbh$ if necessary, we can assume that $\exp$ is a diffeomorphism of $\WNbh$ onto a neighborhood of the diagonal $\Delta=\{(x,x)\,|\,x\in\Nman\}\subset\Nman\times\Nman$.

Notice also that there exists a $\eps>0$ such that if $a,b\in\Nman$ and $d(a,b)<\eps$, then $(a,b)\in\exp(\WNbh)\subset\Nman\times\Nman$ and there is a unique geodesic of length $d(a,b)$ connecting these points.

Let now $\mrsfunc\in\smmn$, $\Gamma_{\mrsfunc}=\{ (x,\mrsfunc(x)) \ | \ x\in\Mman\}\subset\Mman\times\Nman$
be the graph of $\mrsfunc$ in $\Mman\times\Nman$, and
$\feNbh=\mathop\cup\limits_{x\in\manif}x\times\anbh_{\eps}(\mrsfunc(x))$ a neighborhood of $\Gamma_{\mrsfunc}$ in $\Mman\times\Nman$,
where $\anbh_{\eps}(\mrsfunc(x))$ is an open $\eps$-neighborhood of $\mrsfunc(x)$ in $\Nman$.

Let $\Nbh_{\eps}(\mrsfunc)$ be the subset of $\smmn$ consisting of mappings whose graph is included in $\feNbh$.
Then $\Nbh_{\eps}(\mrsfunc)$ is an open neighborhood of $\mrsfunc$ in a strong $C^{0}$ Whitney topology. 
Suppose that $\gfunc\in\Nbh_{\eps}(\mrsfunc)$. 
Then for every $x\in\Mman$ we have that $(\mrsfunc(x),\gfunc(x))\in\exp(\WNbh)$ and the points $\mrsfunc(x)$ and $\gfunc(x)$ are connected with a unique geodesic in $\Nman$ of length $d(\mrsfunc(x),\gfunc(x))<\eps$.

It follows that $\exp^{-1}(\mrsfunc(x),\gfunc(x))=(\mrsfunc(x),\xrel_{\gfunc}(x))\in\WNbh$, where $\xrel_{\gfunc}(x) \in T_{\mrsfunc(x)}\Nman$.
Thus $\gfunc$ gives rise a mapping $\tgel:\Mman\to\Mman\times T\Nman$ defined by 
$\tgel(x)=(x,\mrsfunc(x),\xrel_{\gfunc}(x))$.
Such a mapping can be regarded as a section on the pullback $\mrsfunc^{*}T\Nman$.

Conversely, every smooth mapping $\thel:\Mman\to\Mman\times\WNbh\subset\Mman\times T\Nman$ of the form $\thel(x)=(x,\mrsfunc(x),\xrel(x))$, where $\xrel(x)\in T_{\mrsfunc(x)}\Nman$, yields a smooth mapping $\hel:\Mman\to\Nman$ defined by $\hel(x)=\exp(\mrsfunc(x),\xrel(x))$.

In other words, we can identify a neighborhood of $\mrsfunc$ in $\smmn$ with an open subset $\ZNbh$  of the tame \Frechet\ space $\ssmftn$ of sections of $\mrsfunc^{*}T\Nman$.
It follows that $\smmn$ is a tame \Frechet\ manifold whose tangent space at $\mrsfunc$ is $\ssmftn$, see for details~\cite[Example~I.4.1.2 and Theorem~II.2.3.1]{Hamilton}. 
\end{proof}

In particular, if $\Mman$ is closed, then the group of diffeomorphisms $\Diff(\Mman)$ is an open subset of $\Cinf(\Mman,\Mman)$ and therefore is a tame \Frechet\ manifold.  
The tangent space to $\Diff(\Mman)$ at $\id_{\Mman}$ is the space  
$$T_{\id_{\Mman}}\Diff(\Mman) = \Sinf(\Mman,\id_{\Mman}^{*}T\Mman)=\Sinf(\Mman,T\Mman)=\Gamma(\Mman)$$ of vector fields on $\Mman$.

\subsubsection{The group $\DMS$}
Let $\Mman$ be a compact manifold (possibly with boundary), $\Sigma$ a discrete (possibly empty) subset of $\Int\Mman$, and $\DMS$ of diffeomorphisms of $\Mman$ that preserve $\Sigma$.
Let also $\VMPD$ be the space of vector fields on $\Mman$ that vanish at $\Sigma$ and are tangent to $\partial\Mman$, i.e. $\tgel(\Sigma)=0$ and $\tgel(\partial\Mman)\subset T\partial\Mman$ for all $\tgel\in\VMPD$.
Evidently, $\VMPD$ is a closed subspace of the tame \Frechet\ space $\Gamma(\Mman)$ of all vector fields on $\Mman$, and therefore is a tame \Frechet\ space itself.
\begin{lem}
The group $\DMS$ is a tame \Frechet\ manifold. Its tangent space at $\id_{\Mman}$ is ${\VMPD}$.
\end{lem}
\begin{proof}
Choose a Riemannian metric $d$ on $\Mman$ in which $\partial\Mman$ is totally geodesic, i.e. consists of full geodesics and apply the construction of Lemma~\ref{lm:smmn_tame_Frechet_man}.
Then we get a mapping $\omega$ from some neighborhood $\Nbh$ of $\id_{\Mman}$ in $\Diff(\Mman)$ onto a neighborhood $\ZNbh$ of zero-vector field in $\Gamma(\Mman)$.

Suppose that $\gel\in\Nbh\cap\DMS$. 
We claim that $\tgel=\omega(\gel)\in\VMPD$.

Indeed, let $z\in\Sigma$.
Since $\Sigma$ is discrete and $\gel$ is close to $\id_{\Mman}$, we may assume that $\gel(z)=z$, whence the geodesic connecting $\gel(z)$ and $z$ is just a point and therefore $\tgel(z)=0$.

Further, if $z\in\partial\Mman$, then $\gel(z)\in\partial\Mman$ and since $\partial\Mman$ is totally geodesic, we see that the geodesic connecting $z$ and $\gel(z)$ is included in $\partial\Mman$. 
This implies that $\tgel(x)\in T\partial\Mman$, i.e. $\tgel$ is tangent to $\partial\Mman$.
Thus $\tgel\in\VMPD$.
It can also be shown that the image of $\omega(\Nbh\cap\DMS)$ is a neighborhood of zero-vector field in $\VMPD$.
We leave the details to the reader.
\end{proof}

\subsubsection{The space $\smoned$}
Let $\Mman$ be a compact manifold, $\Psp$ be either a real line $\RRR$ or a circle $\aCircle$, and
$\smoned$ be the space of smooth mapping $\Mman\to\Psp$ that take constant values on the connected components of $\partial\Mman$ (locally constant on $\partial\Mman$).
In particular, $\smrd$ is the space of smooth functions that are locally constant on $\partial\manif$.
It is a closed subspace of a tame \Frechet\ space $\smr$, whence $\smrd$ is a tame \Frechet\ space itself.

\begin{lem}
The space $\smoned$ is a tame \Frechet\ manifold. Its tangent space at each point $\mrsfunc\in\smoned$ is $\smrd$.
\end{lem}
\begin{proof}
Let $\mrsfunc\in\smoned$. Then by Lemma~\ref{lm:smmn_tame_Frechet_man} the tangent space to $\smone$ at $\mrsfunc$ is the space $\Sinf(\Mman,\mrsfunc^{*}T\Psp)$.
Since the tangent bundle $T\Psp$ is trivial (for both cases $\Psp=\RRR$ or $\aCircle$), we see that $\mrsfunc^{*}T\Psp\approx\Mman\times\RRR$ is also trivial, whence $T_{\mrsfunc}\smoned\approx\Sinf(\Mman,\mrsfunc^{*}T\Psp)\approx\smr$.
Notice that $\smrd$ is a closed subspace of $\smr$, and therefore is a tame \Frechet\ space itself.

In particular, we can identify a neighborhood $\Nbh$ of $\mrsfunc$ in $\smone$ with a neighborhood $\ZNbh$ of zero-function in $\smr$ via some homeomorphism $\omega:\Nbh\to\ZNbh$.
For the proof of this lemma it suffices to show that $\omega(\smoned\cap\Nbh)=\smrd\cap\ZNbh$.

Let $\gel\in\smoned\cap\Nbh$ and $\tgel=\omega(\gel)\in\ZNbh$.
Since $\mrsfunc$ and $\gfunc$ are locally constant on $\partial\Mman$, it follows that $(\mrsfunc(x),\gfunc(x))=(\mrsfunc(y), \gfunc(y))$ provided $x,y$ belong to same path-component of $\partial\Mman$. 
Therefore, the geodesic connecting $\mrsfunc(x)$ and $\gfunc(x)$ also connects $\mrsfunc(y)$ and $\gfunc(y)$ (this is a tautology).
Hence $\tgel(x)=\tgel(y)$. Thus $\tgel$ is locally constant on $\partial\Mman$, i.e. $\tgel\in\smrd$.

Conversely, if $\gel\in\Nbh$ and $\tgel=\omega(\gel)\in\smrd\cap\ZNbh$, then similar arguments show that $\gel$ is locally constant on $\partial\Mman$.
\end{proof}

\subsection{Proof of Theorem~\ref{th:loc-triv-fibering}}
\label{sect:proof_th_loctriv}
Let $\Mman$ be a smooth compact connected manifold, and $\bnd_1,\ldots,\bnd_b$ all of the connected components of $\partial\Mman$.
Let $\mrsfunc\in\smoned$ be a Morse mapping with critical points $z_1,\ldots,z_c$.

Instead of the action~\eqref{equ:action-DM} of $\DiffM$ on $\smoned$ it is more convenient to consider a {\em right\/} action $\Gact:\DiffM\times\smoned\to\smoned$ defined by: $\Gact(\difM,\mrsfunc)=\mrsfunc\circ\difM$. 
This action is smooth tame, has same orbits as~\eqref{equ:action-DM}, and differs from~\eqref{equ:action-DM} by the inversion of $\DiffM$. 
Nevertheless, Theorem~\ref{th:Orb_struct} can also be applied to this case.

Then the tangent mapping 
$\Dmp{\Gact}{\id_{\manif},\mrsfunc}:T_{\id_{\manif}}\DiffM\to T_{\mrsfunc}\smoned$ is in fact the linear mapping 
$\Dm{}:\VMd\to\smrd$ defined by the following formula:
$$ D(\txrel)=d\mrsfunc(\txrel).$$

{\bf Action of $\DiffM$.}
Due to Theorem~\ref{th:Orb_struct} it suffices to find $c+b$ smooth functions $\cfunc_1,\ldots,\cfunc_{c+b}$ and construct a linear tame mapping $\Lsect=(\LsectR,\LsectM):\smrd\to\RRR^{c+b}\times\VMd$ such that for each $\gfunc\in\smrd$ we will have:
$\gfunc=\sprod{ \LsectR(\gfunc) }{\cfunc} + d\mrsfunc( \LsectM(\gfunc))$.

Let $\Nbh=\{\nbh, \nbh_1,\ldots,\nbh_c, \anbh_1,\ldots\anbh_b\}$ be a covering of $\Mman$ such that $\nbh_i$ is an open neighborhood of $z_i$, $\anbh_j$ is an open neighborhood of $\bnd_j$, and 
$\overline{\nbh_i}\cap\overline{\nbh_i'} =
\overline{\nbh_i}\cap\overline{\anbh_j} =
\overline{\anbh_j}\cap\overline{\anbh_j'} = \varnothing,$
for $i\not=i'=1,\ldots,c$ and $j\not=j'=1,\ldots,b$.

Let $\mu_i (i=1,\ldots,c),\nu_j (j=1,\ldots,b):\Mman\to[0,1]$ be smooth functions that such that $\supp\mu_i\subset\nbh_i$, $\supp\nu_j\subset\anbh_j$, 
$\mu_i=1$ in a neighborhood of $z_i$, and $\nu_j=1$ in a neighborhood of $\bnd_i$.
Set $\theta = 1 - \sum_i \mu_i -\sum_j \nu_j$.
Then $\theta, \mu_i, \nu_j \in T_{\mrsfunc}\smoned=\smrd$ constitute a partition of unity subordinated to the covering $\Nbh$.

Let $\gfunc\in\smrd$. 
Then $\gfunc = \gfunc\theta + \sum_i \gfunc\mu_i + \sum_j \gfunc\nu_j$.
Choose some Riemannian metric on $\Mman$ and let $\grad\mrsfunc$ be the gradient of $\mrsfunc$ with respect to this metric.
We will show that $\mu_i$ and $\nu_j$ may stand for $\cfunc_k$.

1) Notice that $\supp\theta$ is distinct from the critical points of $\mrsfunc$ and $\partial\Mman$.
Let $\Fld=\gfunc\theta \frac{\grad\mrsfunc}{|\grad\mrsfunc|^2}$.
Then $\gfunc\theta = d\mrsfunc(\Fld)$.

2) Since $z_i$ is a non-degenerate critical point of $\mrsfunc$, we can assume (decreasing $\nbh_i$ is necessary) that in some local coordinates $(x_1,\ldots,x_m)$ near $z_i=0$ we have 
$\mrsfunc(x_1,\ldots,x_n)=\mrsfunc(0) + \mathop\sum\limits_{s=1}^{n} \eps_{is} \,x_s^2$,
where $\eps_{is}=\pm1$.

By Hadamard Lemma on $\nbh_i$ we have a representation 
$\gfunc(x) = \gfunc(z_i) + \sum_s x_s \tgfunc_{is}(x)$,
where $\tgfunc_{is}$ are some smooth functions linearly and tamely depending on $\gfunc$.

Then the following vector field $\Gfld_i=\frac{\mu_i}{2}(\eps_{i1}\tgfunc_{i1},\ldots,\eps_{im}\tgfunc_{im})$ belongs to $T_{\id_{\Mman}}\DiffM=\VMd$ and has a support in $\nbh_i$.

Hence 
$\gfunc\nu_i = \gfunc(z_i)\nu_i + d\mrsfunc(\Gfld_i)$
and $\gfunc(z_i)$ and $\Gfld_i$ linearly and tamely depend on $\gfunc$.

3) We can assume that $\anbh_j$ is a collar for $\bnd_j$, i.e. $\anbh_j$ is diffeomorphic with $\bnd_j\times[0,1)$ so that $\bnd_j$ corresponds to $\bnd_j\times0$ and $\mrsfunc(x,t) = \eps_j t+\mrsfunc(\bnd_j)$,
where $\eps_j=\pm1$ and $(x,t)\in\bnd_j\times[0,1)$.
Then $d\mrsfunc(x,t)=\eps_j dt$.

Let $\Hfld_j=\eps_j \nu_j(0,\gfunc(x,t)-\gfunc(\bnd_j))$ be a vector field on $\anbh_j$ (recall that $\gfunc$ is constant of $\bnd_j$).
Then $\gfunc(x,t)\nu_j = \gfunc(\bnd_j)\nu_j + d\mrsfunc(\Hfld_j)$.
Again $\gfunc(\bnd_j)$ and $\Hfld_j$ linearly and tamely depends on $\gfunc$.

Thus 
$\gfunc = d\mrsfunc\left(\Fld + \sum_i\Gfld_i + \sum_j \Hfld_j\right) + 
\sum_i \gfunc(z_i)\mu_i + \sum_j \gfunc(\bnd_j)\nu_j.$

It remains to note that $\mu_i,\nu_j$ are linearly independent, as they have disjoint supports.
Moreover, none of them can be represented in the form $d\mrsfunc(\Fld)$ for some vector field $\Fld\in\VMd$.
Indeed, if $\Fld\in\VMd$, then $d\mrsfunc(\Fld)=0$ at each $z_i$ (since $z_i$ is critical) and each $\bnd_j$ (since $\Fld$ is tangent to $\bnd_j$ and $\mrsfunc$ is constant on $\bnd_j$),
while $\mu_i(z_i)=\nu_j(\bnd_j)=1$.

{\bf Action of $\DiffMcr$.}
Let $\DiffMcr$ be the group of diffeomorphisms of $\Mman$ preserving the set of critical points of $\mrsfunc$.
We have to find $c\dimM+c+b$ functions $\cfunc_k$ and a tame linear map
$$\Lsect=(\LsectR,\LsectM):\smrd\to\RRR^{c\dimM+c+b}\times\VMdcr.$$ 
The proof is similar to the previous case and differs from it only at step 2).
We only indicate this difference.

2) Again on $\nbh_i$ we have a representation 
$$\gfunc(x) = \gfunc(z_i) + \sum_{s=1}^{m} \gfunc'(z_i) x_s + \sum_{s=1}^{m} \thfunc_{is}(x) x_s,$$
where $\thfunc_{is}$ are smooth functions that linearly and tamely depend on $\gfunc$ and such that $\thfunc_{is}(z_i)=0$.
Then a vector field 
$\Gfld_i=\frac{\mu_i}{2}(\eps_{is}\thfunc_{i1},\ldots,\eps_{is}\thfunc_{im})$ belongs to $T_{\id_{\Mman}}\DiffMcr=\VMdcr$, i.e. $\Gfld_i(z_i)=0$.

Let $\mu_{is}=x_s\mu_i$ for $s=1,\ldots,m$.
Then the functions $\mu_i, \mu_{is},\nu_j$ constitute a complementary basis in $\smrd$ to $d\mrsfunc(\VMdcr)$.
This completes Theorem~\ref{th:loc-triv-fibering}.\qed


\section*{Acknowledgements}
I am sincerely grateful to V.~V.~Sharko, D.~Bolotov, A.~Mozgova, M.~Pankov, E.~Polulyah, A.~Prishlyak, and I.~Vla\-sen\-ko for help, useful discussions and interest to this work.
I thank F.~Sergeraert for referring me to the paper by V.~Po\'enaru~\cite{Poenaru}.
I am indebted to the anonymous referee for careful reading this manuscript, valuable comments and suggestions which allow clarify the paper.

\end{document}